\documentclass[12pt,psamsfonts]{article}
\usepackage{amsmath}
\usepackage{amsfonts}
\usepackage{amssymb}
\usepackage{amsthm}
\usepackage{ragged2e}
\usepackage[mathscr]{eucal}
\usepackage{xspace}
\usepackage{enumerate}
\usefont{T1}{cmr}{m}{n}
\usepackage[russian,english]{babel}
\newtheoremstyle{my}%
 {4.0ex}{2.0ex}%
   {\itshape}%
     {0.0pt}%
       {\sffamily\mdseries}{\,:}%
          {1.0ex} { }%
\newtheoremstyle{my*}%
  {4.0ex}{2.0ex}%
     {\itshape}%
       {0.0pt}%
         {\sffamily\mdseries}{\,:}%
           {1.0ex} {#1 #3}%
\theoremstyle{my}
\newtheorem{theorem}{THEOREM}[section]
\newtheorem{lemma}{LEMMA}[section]
\newtheorem{definition}{DEFINITION}[section]

\newtheorem{remark}{REMARK}[section]

\theoremstyle{my*}
\newtheorem{nonumtheorem}{THEOREM}
\usepackage{titlesec}
\titleformat{name=\section}%
   {\raggedright\normalfont\small\bfseries\centering}%
     {\normalsize\thesection.}{0.5em}{}
\titlespacing{name=\section}{2.0ex}{3.0ex}{1.0ex}[0.0ex]
\titleformat{name=\subsection}%
   {\raggedright\normalfont\small\bfseries}%
     {\normalsize\thesubsection}{0.5em}{}
\DeclareSymbolFont{rsfs}{U}{rsfs}{m}{n}
\DeclareSymbolFontAlphabet{\mathscrn}{rsfs}

\DeclareMathOperator{\E}{\mathit{E}}

\DeclareMathOperator{\Mm}{\text{\raisebox{0.3ex}{\(\scriptscriptstyle\circ\)}}}
\DeclareMathOperator{\dist}{\textup{dist}}
\DeclareMathOperator{\Psii}{\Psi}
 \renewcommand{\theequation}{\thesection.\arabic{equation}}
 \numberwithin{equation}{section}

\begin{document}
\title{\textbf{ON H.WEYL AND H.MINKOWSKI POLYNOMIALS}}
\author{\textsf{Victor Katsnelson}
\thanks{Department of mathematics, the Weizmann Institute, Rehovot, 76100,
Israel.\newline
\hspace*{4.0ex}\texttt{victor.katsnelson@weizmann.ac.il,\ \
victorkatsnelson@gmail.com}}}

\date{06\,February\,2006}
\maketitle{}

\hfill
\begin{minipage}[l]{38.0ex}
 \footnotesize{ Erasm Darwin, the nephew of the great scientist Charles Darwin, believed
 that sometimes one should perform the most unusual experiments.
 They usually yield no results but when they do~...\,.
 So once he played trumped in front of tulips for the whole day. The experiment yielded
 no results.}
\end{minipage}

\vspace*{2.0ex} We introduce certain polynomials, so-called H.Weyl
and H.Minkowski polynomials, which have a geometric origin. The
location of roots of these polynomials is studied.
\section{H.\,WEYL AND H.MINKOWSKI POLYNOMIALS.}
Let \(\mathscr{M}\) be a smooth manifold of dimension \(n\):
\[\dim\mathscr{M}=n,\]
which is embedded injectively into the Euclidean space of a higher
dimension, say \(n+p\), \(p>0\). We identify \(\mathscr{M}\) with
the image of this embedding, so we consider \(\mathscr{M}\) as a
subset  of \(\mathbb{R}^{n+p}\):
\[\mathscr{M}\subset\mathbb{R}^{n+p}.\]
For \(x\in\mathscr{M}\), let \(\mathscr{N}_x\) be the normal
subspace to \(\mathscr{M}\) at the point \(x\). \(\mathscr{N}_x\)
is an affine subspace of the ambient space \(\mathbb{R}^{n+p}\),
\[\dim\mathscr{N}_x=p.\]
For \(t>0\), let
\begin{equation}%
\label{Ndisk}%
 D_x(t)=\{y\in\mathscr{N}_x:\,\dist(y,x)\leq{}t\},
\end{equation}%
 where
\(\dist(y,x)\) is the Euclidean distance between \(x\) and \(y\).
If the manifold \(\mathscr{M}\) is compact, and \(t>0\) is small
enough, then
\begin{equation}
\label{NInt}
D_{x_1}(t)\cap{}D_{x_2}(t)=\emptyset\quad\textup{for}\ \
x_1\in\mathscr{M}, x_2\in\mathscr{M}, x_1\not=x_2.
\end{equation}
\begin{definition}
\label{DTN} The set
\begin{equation}
\label{deftub}
\mathfrak{T}_{\,\,\mathscr{M}}^{\mathbb{R}^{n+p}}(t)
\stackrel{\textup{\tiny{}def}}{=}\bigcup_{x\in\mathscr{M}}D_x(t)
\end{equation}
is said to be \textsf{the tube neighborhood} of the manifold
\(\mathscr{M}\), or \textsf{the tube} around \(\mathscr{M}\). The
number \(t\) is said to be \textsf{the radius} of this tube.
\end{definition}
Is it clear that for manifolds \(\mathscr{M}\) \textit{without
boundary},
\begin{equation}
\label{DeTu}%
\mathfrak{T}_{\,\,\mathscr{M}}^{\mathbb{R}^{n+p}}(t)
=\{x\in\mathbb{R}^{n+p}:\,\textup{dist}(x,\mathscr{M})\leq t\},
\end{equation}
where \(\dist(x,\mathscr{M})\) is the Euclidean distance from
\(x\) to \(\mathscr{M}\). Thus, for manifolds without boundary,
the equality  \eqref{DeTu} could also be taken as a definition of
the tube \(\mathfrak{T}_\mathscr{M}(t)\). However, for manifolds
\(\mathscr{M}\) \textit{with boundary} the sets
\(\mathfrak{T}_{\,\,\,\mathscr{M}}^{\mathbb{R}^{n+p}}(t)\) defined
by \eqref{deftub} and \eqref{DeTu}
 do not coincide. In this, more general,
case the tube around \(\mathscr{M}\) should be defined by
\eqref{deftub}, but not by \eqref{DeTu}. Hermann Weyl,
\cite{Wey1}, obtained the following result, which is the starting
point of our work:\\

\textsf{THEOREM [H.Weyl].} \ %
 \textit{Let \(\mathscr{M}\) be a smooth compact manifold, without boundary or with boundary,
 of the dimension \(n\):
 \(\dim{}\mathscr{M}=n \),
which is embedded in the Euclidean space} \(\mathbb{R}^{n+p},\
p\geq{}1\).\\[-5.0ex]
\begin{enumerate}
\item[I.]%
\begin{itshape}If \(t>0\) is small enough%
\footnote{\label{Smt}If the condition \eqref{NInt} is satisfied.}, %
than the \((n+p)\)\,-\,dimensional volume \(\textup{Vol}_{n+p}\) %
of the tube
\(\mathfrak{T}_{\,\,\,\mathscr{M}}^{\mathbb{R}^{n+p}}(t)\) around
\(\mathscr{M}\), considered as a function of the radius \(t\) of
this tube%
 , is a polynomial of the form
\begin{equation}
\label{TuVo}
\textup{Vol}_{n+p}(\mathfrak{T}_{\,\,\,\mathscr{M}}^{\mathbb{R}^{n+p}}(t))=
\omega_p\,t^{p}\Big(
\sum\limits_{l=0}^{[\frac{n}{2}]}w_{2l,p}(\mathscr{M})\cdot{}t^{2l}\Big),
\end{equation}
where
\begin{equation}%
\label{Vpdb}%
 \omega_p=\frac{\pi^{p/2}}{\Gamma(\frac{p}{2}+1)}
\end{equation}%
is is the \(p\)-dimensional volume of the unit
\(p\)\,-\,dimensional ball.
\end{itshape}\\
\item[II.]
\begin{itshape}
 The coefficients \(w_{2l,p}(\mathscr{M})\) depend on \(p\) as:
 \begin{equation}
 \label{DeOnDi}
w_{2l,p}(\mathscr{M})=%
\frac{2^{-l}\,\Gamma(\frac{p}{2}+1)}{\Gamma(\frac{p}{2}+l+1)}%
\cdot{}k_{2l}(\mathscr{M})\,,
 \quad{} 0\leq{}l\leq{}\left[\textstyle\frac{n}{2}\right]\,,
 \end{equation}
where \emph{\textsf{the values
\(k_{2l}(\mathscr{M}),\,0\leq{}l\leq[\frac{n}{2}]\), may be
expressed only in terms of the intrinsic metric%
\footnote{That is the metric which is induced on
manifold \ \(\mathscr{M}\) from the ambient space \(R^{n+p}\).}
 of the manifold \ \(\mathscr{M}\).
}} In particular, the constant term \(
w_{0,p}(\mathscr{M})=k_{0}(\mathscr{M})\) is the \(n\)-dimensional
volume of \(\mathscr{M}\):
 \begin{equation}%
 \label{CTV}%
 k_0(\mathscr{M})=\textup{Vol}_n\,(\mathscr{M}).
 \end{equation}%
 \end{itshape}
\end{enumerate}
H.\,Weyl, \textup{\cite{Wey1}},  have expressed the coefficients
\(k_{2l}(\mathscr{M})\) as integrals of certain rather complicated
curvature functions of the manifold \(\mathscr{M}\).
\begin{remark}%
\label{EuC} In the case when \(\mathscr{M}\) is compact without
boundary and even dimensional, say \(n=2m\), the top coefficient
\(k_{2m}(\mathscr{M})\) is especially interesting:
\begin{equation}
\label{ToCo}%
 k_{2m}(\mathscr{M})=(2\pi)^{m}\chi(\mathscr{M}),
\end{equation}
where \(\chi(\mathscr{M})\) is the Euler characteristic of
\(\mathscr{M}\).
\end{remark}%

\begin{definition} %
\label{DEfWP}%
 Let \(\mathscr{M}\) be a smooth manifold, without boundary or with boundary,
  of the dimension \(n\):
 \(\dim{}\mathscr{M}=n \),
which is embedded in the Euclidean space \(\mathbb{R}^{n+p},\
p\geq{}1\), and \(\mathfrak{T}_{\,\,\,\mathscr{M}}^{\mathbb{R}^{n+p}}(t)\)
is the tube of the radius \(t\) around \(\mathscr{M}\), \eqref{DTN}.\\[1.0ex]
\hspace*{2.0ex} The polynomial \(W_{\,\mathscr{M}}^{\,p}(t)\)
which appears in the expression \eqref{TuVo} for the volume
\(\mathrm{Vol}_{n+p}\,\big(\mathfrak{T}_{\,\,\,\mathscr{M}}^{\mathbb{R}^{n+p}}(t)\big)\)
of this tube:
\begin{equation}%
\label{Soo}%
\mathrm{Vol}_{n+p}\,(\mathfrak{T}_{\,\,\,\mathscr{M}}^{\mathbb{R}^{n+p}}(t))=
\omega_pt^p\cdot{}W_{\,\mathscr{M}}^{\,p}(t)\quad\textup{ for
small positive \(t\),}
\end{equation} %
is said to be \textsf{the H. Weyl
polynomial of the index \(p\) for the manifold \(\mathscr{M}\).}
\end{definition}
\begin{remark}
The radius \(t\) of the tube is a positive number, so the formula
\eqref{Soo} is meaningful for positive \(t\) only. However the
polynomial \(W_{\mathscr{M}}^{\,p}\) is determined uniquely by  its
restriction on any fixed interval \([0,\varepsilon]\),
\(\varepsilon>0\), and we may and will consider this polynomial
for \textsf{every complex} \(t\).
\end{remark}

\begin{definition}
\label{DeWC}%
 Let \(\mathscr{M}\) be a smooth manifold of the dimension \(n\):
 \(\dim{}\mathscr{M}=n \),
which is embedded in the Euclidean space \(\mathbb{R}^{n+p},\
p\geq{}1\), and let \(W_{\mathscr{M}}^{\,p}\) be the Weyl polynomial
of \(\mathscr{M}\) \textup{(}defined by \eqref{NInt},\,\eqref{Soo}\,\textup{)}.\\[1.0ex]
\hspace*{2.0ex} The coefficients
\(k_{2l}(\mathscr{M}),\,0\leq{}l\leq{}[n/2]\)  which are
\textsf{defined} in terms of  the Weyl polynomial
\(W_{\,\mathscr{M}}^{\,p}\) by the equality
\begin{equation}%
\label{DefWC}%
W_{\,\mathscr{M}}^{\,p}(t)\stackrel{\textup{\tiny
def}}{=}\sum\limits_{l=0}^{[\frac{n}{2}]}%
\frac{2^{-l}\,\Gamma(\frac{p}{2}+1)}{\Gamma(\frac{p}{2}+l+1)}\,%
k_{2l}(\mathscr{M}) \cdot{}t^{2l}\,,
\end{equation}%
are said to be the Weyl coefficients of the manifold
\(\mathscr{M}\).
\end{definition}
\begin{remark}
\label{DecFo}%
 Often, the factor in \eqref{DefWC} appears in a
`decoded' form:
\begin{equation}
\label{DecFor}
\frac{2^{-l}\,\Gamma(\frac{p}{2}+1)}{\Gamma(\frac{p}{2}+l+1)}=
\frac{1}{(p+2)(p+4)\,\cdots\,\,(p+2l)}\,.
\end{equation}
\end{remark}
\begin{remark}
\label{RieMe} Defining the Weyl polynomials
\(W_{\,\mathscr{M}}^{\,p}\) of the manifold \(\mathscr{M}\) by
\eqref{Soo}, we assumed that \(\mathscr{M}\) is already embedded
into \(\mathbb{R}^{n+p}\). The tube around \(\mathscr{M}\) and its
volume are primary in this definition. So, in fact we defined the
notion of the Weyl polynomial not for the manifold \(\mathscr{M}\)
\textsf{itself} but for manifold \(\mathscr{M}\) which is already
\textsf{embedded} in an ambient space. Moreover, we assume
implicitly that from the very beginning the manifold
\(\mathscr{M}\) carries a `natural' Riemannian metric, and that
this `original' Riemannian metric coincides with the metric on
\(\mathscr{M}\) induced  from the ambient space
\(\mathbb{R}^{n+p}.\) (In other words, we assume that the
imbedding is isometrical.) However, in this approach the
`original' metric  does not play
 an `explicite' role in the definition \eqref{DTN}-\eqref{Soo}-\eqref{DefWC}
of the Weyl polynomial \(W_{\,\mathscr{M}}^{\,p}\) and the Weyl
coefficients \(k_{2l}(\mathscr{M})\).

 There is another approach to define the Weyl coefficients
 and the Weyl polynomials, which does not require an
actual embedding \(\mathscr{M}\) into the ambient space.
 Starting from the given Riemannian metric on
\(\mathscr{M}\), the Weyl coefficients \(k_{2l}(\mathscr{M})\) can
be introduced  formally, by means of the Hermann Weyl expressions
for \(k_{2l}(\mathscr{M})\) in terms of the given metric on
\(\mathscr{M}\). Then the Weyl polynomials
\(W_{\,\mathscr{M}}^{\,p}(t)\) can be defined by means of the
expression \eqref{DefWC}. In this approach, the intrinsic metric
of \(\mathscr{M}\) is primary, but not the tubes around
\(\mathscr{M}\) and their volumes.
\end{remark}
If the codimension \(p\) of
\(\mathscr{M}\)
 equals one\footnote{In other words,
\(\mathscr{M}\) is a hypersurface in \(\mathbb{R}^{n+1}\).},
\,\,\(\dim{}\mathscr{M}=n,\)  the
Weyl polynomial is of the form:
\begin{equation}
\label{PW}
\textup{Vol}_{n+1}(\mathfrak{T}_{\,\,\,\mathscr{M}}^{\mathbb{R}^{n+1}}(t))=
2t\cdot{}W_{\mathscr{M}}^{1}(t)\,,\quad
W_{\mathscr{M}}^{1}(t)=
\sum\limits_{l=0}^{[\frac{n}{2}]}w_{2l}(\mathscr{M})\cdot{}t^{2l},
\end{equation}
where
\begin{equation}%
\label{WKRel}%
w_{2l}(\mathscr{M})=
\frac{2^{-l}\Gamma(\frac{1}{2})}{\Gamma(\frac{1}{2}+l+1)}\,k_{2l}(\mathscr{M})\,,\quad
 {}0\leq{}l\leq{}[\textstyle\frac{n}{2}]\,.
\end{equation}%
In \eqref{PW} the `shortened' notation is used:
\(w_{2l}(\mathscr{M})\) instead of \(w_{2l,1}(\mathscr{M})\). The
factor \(2t\) is the one-dimensional volume of the one-dimensional
ball of radius \(t\), that is the length of the interval
\([-t,t]\).

If the \textit{hypersurface} \(\mathscr{M}\) is
orientable\,\footnote{The orientation of the hypersurface
\(\mathscr{M}\) can be specified by means of the continuous vector
field of unit normals on \(\mathscr{M}\). The half-tubes
\(\mathfrak{T}_\mathscr{M}^{+}(t)\) and
\(\mathfrak{T}_\mathscr{M}^{-}(t)\) are the parts of the tube
\(\mathfrak{T}_\mathscr{M}(t)\) corresponding to the `positive' и
`negative' directions of these normals.}, then  the tube
\(\mathfrak{T}_\mathscr{M}(t)\) can be decomposed into the union
of two \textit{half-tubes}, say,
\(\mathfrak{T}_\mathscr{M}^{+}(t)\) and
\(\mathfrak{T}_\mathscr{M}^{-}(t)\). The half-tubes
\(\mathfrak{T}_\mathscr{M}^{+}(t)\) and
\(\mathfrak{T}_\mathscr{M}^{-}(t)\) are the parts of the tube
\(\mathfrak{T}_\mathscr{M}(t)\) which are situated  on the
distinct sides of \(\mathscr{M}\).
 In particular, if the hypersurface \(\mathscr{M}\) is the boundary of a set
\(V:\,\mathscr{M}=\partial{}V\), then
\begin{equation}
\mathfrak{T}_\mathscr{M}^{+}(t)=\mathfrak{T}_\mathscr{M}(t)\setminus
V,\quad{}\mathfrak{T}_\mathscr{M}^{-}(t)=\mathfrak{T}_\mathscr{M}(t)\cap{}V\,.
\end{equation}
The \((n+1)-\,dimensional\) volumes
\(\textup{Vol}_{n+1}(\mathfrak{T}_\mathscr{M}^{\,+}(t))\) and
\(\textup{Vol}_{n+1}(\mathfrak{T}_\mathscr{M}^{\,-}(t))\) of the
half-tubes also are polynomials of \(t\). These polynomials are of
the form\,%
\footnote{The equalities \eqref{WhT}, \eqref{WPST} is one of the
results of the theory of tubes around manifolds. See
\cite{Gr},\,\cite{BeGo},\cite{AdTa}}\,:
\begin{equation}
\label{WhT}
\textup{Vol}_{n+1}(\mathfrak{T}_\mathscr{M}^{\,+}(t))=t\,W_\mathscr{M}^{\,+}(t),\quad
\textup{Vol}_{n+1}(\mathfrak{T}_\mathscr{M}^{\,-}(t))=t\,W_\mathscr{M}^{\,-}(t)\,,
\end{equation}
where:
\begin{subequations}%
 \label{WPST}%
\begin{alignat}{2}
W_\mathscr{M}^{\,+}(t)&=\sum\limits_{l=0}^{[\frac{n}{2}]}w_{2l}(\mathscr{M})\cdot{}t^{2l}&&+
t\,\sum\limits_{l=0}^{[\frac{n+1}{2}]-1}w_{2l+1}(\mathscr{M})\cdot{}t^{2l},\label{WPSTa}\\
W_\mathscr{M}^{\,-}(t)&=\sum\limits_{l=0}^{[\frac{n}{2}]}w_{2l}(\mathscr{M})\cdot{}t^{2l}&&-
t\,\sum\limits_{l=0}^{[\frac{n+1}{2}]-1}w_{2l+1}(\mathscr{M})\cdot{}t^{2l},\label{WPSTb}
\end{alignat}
\end{subequations}
and the coefficients \(w_{2l}(\mathscr{M})\) are the same that
in \eqref{PW}-\eqref{WKRel}. Unlike the coefficients \(w_{2l}(\mathscr{M})\),
the coefficients \(w_{2l+1}(\mathscr{M})\) depend not only on the
`intrinsic' metric of the manifold \(\mathscr{M}\), but also on
how \(\mathscr{M}\) is embedded to \(\mathbb{R}^{n+1}\). It is
remarkable that when the volumes of the half-tubes are summed:
\begin{equation*}
\label{IEPp}
2\,W_\mathscr{M}(t)=W_\mathscr{M}^{\,+}(t)+W_\mathscr{M}^{\,-}(t),
\end{equation*}
the dependence on the way of embedding disappears.
As it is seen from \eqref{WPST},
\(W_\mathscr{M}^{\,-}(t)=W_\mathscr{M}^{\,+}(-t)\), hence
\begin{equation}
\label{IEP}
2\,W_\mathscr{M}(t)=W_\mathscr{M}^{\,+}(t)+W_\mathscr{M}^{\,+}(-t).
\end{equation}
Remark also that the volumes of the half-tubes can be expressed only in the terms
of the polynomial \(W_\mathscr{M}^{\,+}\):
\begin{subequations}
\label{SeSa}
\begin{align}
\label{SeSa+}
\textup{Vol}_{n+1}(\mathfrak{T}_\mathscr{M}^{\,+}(t))&=t\,W_\mathscr{M}^{\,+}(\,\,t\,)
\ \ \textup{\,for small positive \(t\)\,.}\\[1.0ex]
\label{SeSa-}
\textup{Vol}_{n+1}(\mathfrak{T}_\mathscr{M}^{\,-}(t))&=t\,W_\mathscr{M}^{\,+}(-t)
\ \ \textup{for small positive \(t\)\,.}
\end{align}
\end{subequations}
The theory of the tubes around manifolds is presented in
\cite{Gr}, and to some extent in \cite{BeGo}, Chapter 6, and in
\cite{AdTa},\,Chapter 10. The comments of V.Arnold  \cite{Arn} to
the Russian translations of the paper \cite{Wey1} by H.Weyl are
very rich in content.

In the event that the hypersurface \(\mathscr{M}\) is the boundary
of a convex set \(V\): \(\mathscr{M}=\partial V\), the Weyl
polynomial \(W_{\mathscr{M}}^{1}\) can be expressed in terms of
polynomials considered in the theory of convex sets.

In the theory of convex sets the following fact, which was
discovered by Hermann Minkowski, \cite{Min1, Min2}, is of
principal importance: \textit{Let \(V_1\) and \(V_2\) be compact
convex sets in \(\mathbb{R}^n\). For positive numbers \(t_1,t_2\),
let us form the `linear combination' \(t_1V_1+t_2V_2\) of the sets
\(V_1\) и \(V_2\) (in the sense commonly accepted in the theory of
convex sets).  Then  the \(n\)-dimensional Euclidean volume
\(\textup{Vol}_n(t_1V_1+t_2V_2)\) of this linear combination,
considered as a function of the variables \(t_1,t_2\), is a
homogeneous polynomial of  degree~\(n.\) (It may be equal zero
identically.)} Choosing \(V\) as \(V_1\) and the unit ball \(B^n\)
of \(\mathbb{R}^n\) as \(V_2\), we conclude
:\\[1.5ex]
\textit{Let \(V\) be a compact convex set in \(\mathbb{R}^n\),
\(B^n\) be the unit ball of  \(\mathbb{R}^n\). Then
\(n\)-dimensional volume \(\textup{Vol}_n(V+tB^n)\), considered as
a function of the variable \(t\in[0,\infty)\), is a polynomial of
degree \(n\).}
\begin{definition}
\label{DeMiPo}
 Let \(V,\,V\subset{}\mathbb{R}^n,\) be a
compact convex set. The polynomial which expresses the
\(n\)-dimensional volume of the linear combination \(V+tB^n\) as a
function of the variable \(t\in[0,\infty)\) is said to be
\textsf{the Minkowski polynomial of the set \(V\)} and is denoted
by \(M_{\,V}^{\mathbb{R}{n}}(t)\):\\[-1.5ex]
\begin{equation}
\label{DMP} M_{V}^{\mathbb{R}^n}(t)=%
\textup{Vol}_n(V+tB^n)\,,\quad
(t\in[0,\infty)).
\end{equation}
The coefficient of Minkowski polynomial are denoted by
\(m_{\,k}^{\mathbb{R}^n}(V)\):
\begin{equation}
\label{MiP} M_{\,V}^{\mathbb{R}^{n}}(t)=\sum\limits_{0\leq k\leq
n}m_{\,k}^{\mathbb{R}^n}(V)t^k.
\end{equation}
If there is no need to emphasize that the ambient space is
\(\mathbb{R}^n\), the shortened notation \(M_V(t)\), \(m_k(V)\)
for the Minkowski polynomial and its coefficients will be used.
\end{definition}
Of course,
\[M_{\,V}^{\mathbb{R}^{n}}(t)=\textup{Vol}_n(\mathfrak{V}_{\,V}^{\,\mathbb{R}^{n}}(t)),\]
where \(\mathfrak{V}_{\,V}^{\mathbb{\,R}^{n}}(t))\) is
\(t\)-neighborhood of the set \(V\) with respect to
\(\mathbb{R}^n\):
\begin{equation}
\label{DeNeCo}%
\mathfrak{V}_{\,V}^{\,\mathbb{R}^{n}}(t)=\{x\in\mathbb{R}^n:\,\textup{dist}(x,{V})\leq
t\}.
\end{equation}
It is evident that
\begin{equation}
\label{EvEq}%
 m_0(V)=\text{Vol}_n(V), \ \ \text{and}\ \
m_n(V)=\text{Vol}_n(B^n).
\end{equation}
 If the boundary \(\partial V\) of a convex set \(V\) is smooth,
 then the \((n-1)\)-dimensional volume (`the area') of the boundary
\(\partial V\) can be expressed as
\begin{equation}
\label{Cau} m_{1}(V)=\text{Vol}_{\,n-1}(\partial V)\,.
\end{equation}
For a convex set \(V\), whose boundary \(\partial V\) may be
non-smooth, the formula \eqref{Cau} serves as a
\textit{definition} of the \textit{`area'} of \(\partial V\). (See
\cite{BoFe},\,\textbf{31}; \cite{Min1},\,\S\,24;
\cite{Web},\,\textbf{6.4}.)
Let us emphasize that the Minkowski polynomial is defined for an
 \textit{arbitrary} compact convex set \(V\), without any extra
 assumptions. The boundary of \(V\) may be non-smooth, and the
 interior of \(V\) may be empty. In particular, the Minkowski
 polynomial is defined for any convex polytope.
\begin{definition}
\label{solid} Let \(V,\,V\subset\mathbb{R}^n,\) be a convex set.
\(V\) is said to be \textsf{solid} if the interior of \(V\) is not
empty, and \textsf{non-solid} if the interior of \(V\) is not
empty.
\end{definition}
\begin{definition}
The \textsf{\(n\)\,-\,dimensional closed convex surface
\(\mathscr{M}\)} is the boundary \(\partial{}V\) of a solid
compact convex set \(V \):
\begin{equation}
\label{DeCoS}%
\mathscr{M}=\partial{}V,\quad{}V\subset\mathbb{R}^{n+1}\,.
\end{equation}
The set \(V\) is said to be the \textsf{generating set for the
surface \(\mathscr{M}\)}.
\end{definition}
 \begin{lemma}
\label{MPWP}
  \textit{If the closed \(n\)\,-\,dimensional convex surface \(\mathscr{M}\)
  is also a smooth manifold,
    then the Weyl polynomial \(W_{\mathscr{M}}^{\,1}\)
 of the surface \(\mathscr{M}\) and the Minkowski polynomial
  \(M_{\,\,V}^{\mathbb{R}^{n+1}}\)
 of its generating set \(V\) are related in the following way}:
 \begin{equation}
 \label{WMP}
 2t\,W_{\mathscr{M}}^{\,1}(t)=
 M_{\,\,V}^{\mathbb{R}^{n+1}}(t)-M_{\,\,V}^{\mathbb{R}^{n+1}}(-t).
\end{equation}
\end{lemma}
\noindent
 \textsf{PROOF OF LEMMA \ref{MPWP}.} We assign the positive orientation to the vector field
 of exterior normals on \(\partial V\). Let
 \(\mathfrak{T}_{\partial{}V}^{+}(t)\) is
 the `exterior' half-tube around \(\partial V\).
 For  positive \(t\),
\begin{equation*}
V+tB^{n+1}=V\cup{}\mathfrak{T}_{\partial{}V}^{+}(t),
\end{equation*}
Moreover the set \(V\) and \(\mathfrak{T}_{\partial{}V}^{+}(t)\)
do not  intersect. Therefore,
\[\textup{Vol}_{n+1}(V+tB^{n+1})=
\textup{Vol}_{n+1}(V)+\textup{Vol}_{n+1}(\mathfrak{T}_{\partial{}V}^{+}(t)).\]
Hence,
\begin{equation*}
M_V(t)=M_V(0)+t\,W_\mathscr{M}^{\,+}(t),\ \ \ %
{\mathscr{M}=\partial V}, %
\end{equation*}
where \(W_\mathscr{M}^{+}\) is a polynomial defined in
\eqref{WhT} (with \(n\) replaced by \(n+1\): now \(\dim
V=n+1\)). Then also
\begin{equation*}
M_V(-t)=M_V(0)-t\,W_\mathscr{M}^{\,+}(-t).
\end{equation*}
Thus,
\eqref{SeSa},
\begin{equation*}
M_V(t)-M_V(-t)=\textup{Vol}_{n+1}(\mathfrak{T}_{\partial{}V}^{+}(t))+
\textup{Vol}_{n+1}(\mathfrak{T}_{\partial{}V}^{-}(t))\,,
\end{equation*}
or
\begin{equation*}
M_V(t)-M_V(-t)=t\,(W_\mathscr{M}^{\,+}(t)+W_\mathscr{M}^{\,+}(-t)).
\end{equation*}
The equality \eqref{WMP} follows from the last equality and from
\eqref{IEP}.\hfill{Q.E.D.} 

 Since the Minkowski polynomial is defined for an
arbitrary compact convex set, the formula \eqref{WMP} can serve as
a \textit{definition} of the Weyl polynomial of an
\textit{arbitrary} closed convex surface, smooth or non-smooth.
Even more, we can define the Weyl polynomial for the `improper
convex surface \(\partial{}V\)', where \(V\) is a non-solid
compact convex set.
\begin{definition}
\label{Impr} Let \(V,\,V\subset\mathbb{R}^{n+1},\) be a compact
convex set. The boundary \(\partial{}V\) of the set \(V\) is said
to be the \textsf{boundary surface} of \(V\). The boundary surface
of \(V\) is said to be \textsf{proper} if \(V\) is solid, and
\textsf{improper} if \(V\) is non-solid.
\end{definition}

The following improper closed convex surface plays a role in what
follow:
\begin{definition}
\label{SqCyl} Let \(V\), \(V\subset\mathbb{R}^{n}\), be a compact
convex set which is solid \emph{with respect to
\(\mathbb{R}^{n}\).} We identify \(\mathbb{R}^{n}\) with it image
\(\mathbb{R}^{n}\times{}0\) by the `canonical' embedding%
\footnote{The point \(x\in\mathbb{R}^{n}\) is identified with the
point \((x,0)\in\mathbb{R}^{n+1}\).} %
 \(\mathbb{R}^{n}\) into \(\mathbb{R}^{n+1}\), and the set \(V\)
 with the set \(V\times{}0\) considered as a subset of
 \(\mathbb{R}^{n+1}\): \(V\times{}0\subset{}\mathbb{R}^{n+1}\).
 The set \(V\times{}0\),  considered as a subset of
 \(\mathbb{R}^{n+1}\), is said to be \textsf{the squeezed cylinder with
 the base \(V\).}
\end{definition}

\begin{remark}
\label{IntSqCy}%
The set \(V\times{}0\) can be interpreted as a ` cylinder of zero
hight' whose `lateral surface' is the Cartesian product
\(\partial{}V\times[0,0]\), and whose bases, lower and upper, are
the sets \(V\times{}(-0)\) and  \(V\times{}(+0)\)\textup{:}
\begin{equation}
\label{SuArIn}
\partial(V\times{}0)=\big((\partial{}V)\times{}[0,0]\big)\cup\big({V\times{}(-0)}\big)
\cup\big({V\times{}(+0)}\big)\,.
\end{equation}
In other words, the boundary surface \(\partial{}(V\times{}0)\)
 can be considered as `the doubly covered' set \(V\). In particular,
\begin{equation}%
\label{}%
\dim\partial{}(V\times{}0)=n\,.
\end{equation}%
and the number
\(\textup{Vol}_n(V\times(-0))+\textup{Vol}_n(V\times(+0))=2\,\textup{Vol}_n(V)\)
can be naturally interpreted as the `\(n\)-\,dimensional area' of
the \(n\)-\,dimensional convex surface (improper)
\(\partial(V\times{}0)\):
\begin{equation}
\label{SuAr}
\textup{Vol}_n(\partial(V\times{}0))=2\,\textup{Vol}_n(V)\,.
\end{equation}
\end{remark}

On the other hand, the equality \eqref{Cau}, in which the squeezed
cylinder \(V\times{}0\subset\mathbb{R}^{n+1}\) plays the role of
the set \(V\subset{}\mathbb{R}^{n}\), takes the form
\begin{equation}
\label{OOH}
\textup{Vol}_n(\partial(V\times{}0))=m_{\,1}^{\mathbb{R}^{n+1}}(V\times{}0)\,,
\end{equation}
where
\(m_{\,k}^{\mathbb{R}^{n+1}}(V\times{}0),\,\,k=0,1,\,\ldots\,,\,n+1,\)
are the coefficients of the Minkowski polynomial
\(M_{\,V\times{}0}^{\mathbb{R}^{n+1}}(t)\) of the squeezed
cylinder \(V\times{}0\) with respect to the ambient space
\(\mathbb{R}^{n+1}\). (See \eqref{MiP}.)

 In section \ref{ExInSp} we prove the following statement, which
 appears as Lemma \ref{IPRom} there:
\begin{lemma}
\label{IPR}%
Let \(V\) be a compact convex set in \(\mathbb{R}^n\), and
\begin{equation}%
\label{IRP1}%
 M_{\,V}^{\mathbb{R}^n}(t)=\sum\limits_{0\leq{}k\leq{}n}m_{\,k}^{\mathbb{R}^n}(V)t^k
\end{equation}%
be the Minkowski polynomial with respect to the ambient space
\(\mathbb{R}^n\). Then the Minkowski polynomial
\(M_{\,V\times{}0}^{\mathbb{R}^{n+1}}(t)\) with respect to the
ambient space \(\mathbb{R}^{n+1}\) is equal to:
\begin{equation}%
\label{MPn1}%
M_{\,V\times{}0^1}^{\mathbb{R}^{n+1}}(t)=%
t\!\!\sum\limits_{0\leq{}k\leq{}n}\frac{\Gamma(\frac{1}{2})
\Gamma(\frac{k}{2}+1)}{\Gamma(\frac{k+1}{2}+1)}\,m_{\,k}^{\mathbb{R}^n}(V)%
\, t^k\,.
\end{equation}%
\end{lemma}

So,%
 \[m_{\,0}^{\mathbb{R}^{n+1}}(V\times{}0)=0,\quad{}
m_{k+1}^{\mathbb{R}^{n+1}}(V\times{}0)=
\frac{\Gamma(\frac{1}{2})\Gamma(\frac{k}{2}+1)}{\Gamma(\frac{k+1}{2}+1)}\,
m_{\,k}^{\mathbb{R}^n}(V),\,k=0,\,\ldots,\,n\,.\] %
In particular, \(m_{\,1}^{\mathbb{R}^{n+1}}(V\times{}0)
=2m_{\,0}^{\mathbb{R}^n}(V).\) Since
\(m_{\,0}^{\mathbb{R}^n}(V)=\textup{Vol}_n(V)\), \eqref{EvEq},
\begin{equation}
\label{AgMP}
m_{\,1}^{\mathbb{R}^{n+1}}(V\times{}0)=2\,\textup{Vol}_n(V)\,.
\end{equation}

The equalities \eqref{SuAr}, \eqref{OOH}  and \eqref{AgMP} agree.

\begin{remark}%
\label{AprC1} Any non-solid compact convex set \(V\) can be
presented as the limit \textup{(}in the Hausdorff
\,metric\textup{)} of a monotonic\,%
\footnote{\label{monot}The monotonicity means that
\(V_{\varepsilon^\prime}\supset{}V_{\varepsilon^{\prime\prime}}\supset{}V\)
for \(\varepsilon^{\prime}>\varepsilon^{\prime\prime}>0\).}
 family \(\{V_{\varepsilon}\}_{\varepsilon>0}\)
of solid convex sets \(V_{\varepsilon}\) \textup{:}
\[V=\lim_{\varepsilon\to+0}V_{\varepsilon}.\]
Moreover, the approximating family
\(\{V_{\varepsilon}\}_{\varepsilon>0}\) of convex sets
 can be chosen  so that the boundary
\(\partial{}(V_{\varepsilon})\) of each set \(V_{\varepsilon}\) is
a smooth surface. Thus, the improper convex surface
\(\partial{}V\) may be presented as the limit of proper convex
smooth surfaces \(\partial{}(V_{\varepsilon})\) which shrink to
\(\partial{}V\):
\[\partial{}V=\lim_{\varepsilon\to+0}\partial{}(V_{\varepsilon}).\]
\end{remark}%
\begin{definition}
\label{DeWPDC} Let \(V,\,V\subset\mathbb{R}^{n+1},\) be an
arbitrary compact convex set. The \textsf{Weyl polynomial
\(W_{\partial{}V}^{\,1}(t)\)} of the convex surface
\(\mathscr{M}=\partial{}V\),
 proper or improper, is \textsf{defined} by the formula \eqref{WMP}.
 In other words, the Weyl polynomial \(t\,W_{\partial{}V}^{\,1}\)
 is defined as the odd part of the Minkowski polynomial
  \(M_{\,\,V}^{\mathbb{R}^{n+1}}\):
  \begin{equation}%
  \label{MPOP}
  t\cdot{}W_{\partial{}V}^{\,1}(t)={}^{\mathscr{O}}\!M_{\,\,V}^{\mathbb{R}^{n+1}}(t),
  \end{equation}%
  where the notions of the even part \({}^{\mathscr{E}}\!P\) and
  the odd part \({}^{\mathscr{O}}\!P\) of an arbitrary polynomial \(P\) are
  introduced in  \textup{Definition \ref{DefRIPa}} below.
\end{definition}
\begin{remark}
\label{agDe}%
 In the case when the set \(V\) is solid and its
boundary \(\partial{}V\) is smooth, both definitions,
\textup{Definition~\ref{DeWPDC}} \,and
\textup{Definition~\ref{DEfWP}} \,of the Weyl polynomial
\(W_{\,\partial{}V}^{\,1}\), are applicable to \(\partial{}V\). In
this case both definitions agree.
\end{remark}

\begin{remark}%
\label{why}%
\textup{Why may be useful to consider improper convex surfaces and
their Weyl polynomials?}\\
As it was remarked \textup{(Remark \ref{AprC1})}, every improper
convex surface \(\partial{}V\) is a limiting object for a family
of proper smooth convex surfaces \(\partial{}(V_{\varepsilon})\).
It turns out that the Weyl polynomial for this improper surface is
the limit of the Weyl polynomials for this `approximating' family
\(\{V_{\varepsilon}\}_{\varepsilon>0}\) of smooth proper
surfaces.\\
So the Weyl polynomials  for the improper surface \(\partial{}V\)
may be useful in the study of the limiting behavior of the family
of the Weyl polynomials \(\) for the proper surfaces
\(\partial{}(V_{\varepsilon})\) shrinking to the improper surface
\(\partial{}V\). In particular, see \textup{Theorem \ref{NMR}}
formulated in the end of \textup{Section \ref{FMR}}, and its proof
presented in the end of \textup{Section \ref{ExInSp}.}
\end{remark}%

Let \(\mathscr{M}\) be an \(n\)\,-\,dimensional closed convex
surface which is not assumed to be smooth, and \(V\) is the
generating convex set for \(\mathscr{M}\):
\(\mathscr{M}=\partial{}V\). Let \(M_{V}^{\mathbb{R}^{n+1}}\) be
the Minkowski polynomial for \(V\), defined by Definition
\ref{DeMiPo}. According to Definition \ref{DeWPDC}, the Weyl
polynomial \(W_{\mathscr{M}}^{\,1}\) is equal to
\begin{equation}%
\label{WPHC}%
W_{\mathscr{M}}^{\,1}(t)=
\sum\limits_{0\leq{}l\leq{}\left[\frac{n}{2}\right]}m_{2l+1}(V)t^{2l},
\end{equation}%
in other words,
\begin{equation}%
\label{IOW}%
w_{2l}(\mathscr{M})=m_{2l+1}(V),\quad
0\leq{}l\leq{}[\textstyle\frac{n}{2}],
\end{equation}%
 where \(w_{2l}(\mathscr{M})\) are the coefficients of the Weyl polynomial
\(W_{\mathscr{M}}^{1}\), \eqref{PW}, of the \(n\)-\,dimensional
surface \(\mathscr{M}\) with respect to the ambient space
\(\mathbb{R}^{n+1}\), and
 \(m_k(V),\,k=2l+1,\) are the coefficients of
the Minkowski polynomial~~\(M_{\,\,\,V}^{\mathbb{R}^{n+1}}\):
\begin{equation}
\label{MPGs}
M_{\,\,\,V}^{\mathbb{R}^{n+1}}(t)=\textup{Vol}_{n+1}(V+tB^{n+1}),\quad
M_{\,\,\,V}^{\mathbb{R}^{n+1}}(t)=
\sum\limits_{0\leq{}k\leq{}n+1}m_{k}(V)t^{k}\,.
\end{equation}
\begin{definition}
\label{DNWP}%
 Given a closed \(n\)-\,dimensional convex
surface \(\mathscr{M}\), proper or not,
\(\mathscr{M}=\partial{}V\), the numbers
 \(k_{2l}(\mathscr{M}),\,\,0\leq{}l\leq{}[\frac{n}{2}]\), are \textsf{defined} as
\begin{equation}%
\label{NWP}%
k_{2l}(\mathscr{M})=2^l\frac{\Gamma(l+\frac{1}{2}+1)}{\Gamma(\frac{1}{2}+1)}
\,m_{\,2l+1}^{\mathbb{R}^{n+1}}(V),
 \end{equation}
 where \( m_{\,k}^{\mathbb{R}^{n+1}}(V), \,k=2l+1,\,\) are the coefficients of the
 Minkowski polynomial \(M_{\,\,V}^{\mathbb{R}^{n+1}}\) for the
  generating set \(V\), \eqref{MPGs}.
 The numbers
 \(k_{2l}(\partial{}V),\,\,0\leq{}l\leq{}\left[\frac{n}{2}\right]\),
 are said to be \textsf{the Weyl coefficients for the surface
 \(\mathscr{M}\)}.
\end{definition}%
\begin{remark}
According to \textup{Lemma \ref{IPR}}, in the event that the
(improper) convex surface \(\mathscr{M}, \dim \mathscr{M}=n,\) is
the boundary of the squeezed cylinder (see Definition
\ref{SqCyl}), that is if \(\mathscr{M}=\partial{}(V\times{}0),\)
where \(V\subset\mathbb{R}^n\), the Weyl coefficients
\(k_{2l}(\mathscr{M}),\,\,0\leq{}l\leq{}[\frac{n}{2}]\), are:
\begin{equation}
\label{WCSqC}
k_{2l}(\mathscr{M})=2^{l+1}\,\Gamma(l+1)\,m_{2l}^{\mathbb{R}^{n}}(V)\,,
\end{equation}
where \( m_{\,k}^{{\mathbb{R}^{n}}}(V), \,k=2l,\,\) are the
coefficients of the
 Minkowski polynomial \(M_{\,\,V}^{\mathbb{R}^{n}}\) for
  the base   \(V\) of the squeezed cylinder \(\partial{}(V\times{}0)\).
\end{remark}
\begin{remark}
In the case when convex surface \(\mathscr{M}\),
\(\mathscr{M}=\partial{}V\), is smooth and `proper', that is the
set \(V\) generating the surface \(\mathscr{M}\) is solid, both
definitions, Definition~\ref{DNWP} \,and Definition~\ref{DeWC}
\,of the Weyl coefficients \(k_{2l}(\mathscr{M})\) are applicable.
According to \ \eqref{PW}-\eqref{WKRel} and \eqref{IOW}-\eqref{NWP}, in this case%
\footnote{Actually, the equalities \eqref{WKRel}, \eqref{IOW}
served as a motivation for Definition \ref{DNWP}.}
 both definitions agree.
\end{remark}
Note, that according to \eqref{Cau}, (see also Remark
\ref{IntSqCy}),
\begin{equation}
\label{WCA} k_0(\mathscr{M})=\textup{Vol}_{n}(\mathscr{M})
\end{equation}
for every \(n\)\,-\,dimensional closed convex surface
\(\mathscr{M}\).

\begin{lemma}
\label{Posi}%
{\ }\textup{\textsf{I}}. Let \(V\), \(V\subset\mathbb{R}^n\), be a
 solid (with respect to \(\mathbb{R}^n\)) compact convex set. Then
the coefficients \(m_k^{{\mathbb{R}}^n}(V),\,0\leq{}k\leq{}n,\) of
its Minkowski
polynomials%
\footnote{See \eqref{DMP}, \eqref{MiP}.} %
are strictly positive:
 \(m_k^{{\mathbb{R}}^n}(V)>0,\,\,\,0\leq{}k\leq{}n\,.\)\\
\hspace*{2.0ex} \textup{\textsf{II}}. Let \(\mathscr{M}\) be a
proper compact convex surface, \(\dim \mathscr{M}=n.\) Then all
its Weyl  coefficients \(k_{2l}(\mathscr{M})\) are strictly
positive\,\textup{:}\,\,\,
\(k_{2l}(\mathscr{M})>0,\,\,0\leq{}l\leq{}[\frac{n}{2}]\)\,.\\
\hspace*{1.5ex} \textup{\textsf{III}}. Let \(\mathscr{M}\) be the
boundary surface%
\footnote{See Definition \ref{SqCyl} and Remark \ref{IntSqCy}.}%
 of a squeezed cylinder whose base \(V\), \(\dim V=n,\)
is a compact convex set which is solid with respect to
\(\mathbb{R}^n\). Then all its Weyl coefficients
\(k_{2l}(\mathscr{M})\) are strictly positive\,\textup{:}\,\,\,
\(k_{2l}(\mathscr{M})>0,\,\,0\leq{}l\leq{}[\frac{n}{2}]\)\,.
\end{lemma}

The statement \textsf{I} of Lemma \ref{Posi} is a consequence of a
more general statement related to the monotonicity properties of
the mixed volumes. This will be discussed later, in Section
\ref{ChPrMiPo}. The statements \textsf{II} and \textsf{III} of
Lemma \ref{Posi} are consequences of the statement \textsf{I} and
\eqref{NWP}, \eqref{WCSqC}.

\begin{definition}
\label{DGWP} Given a closed \(n\)\,-\,dimensional convex surface
\(\mathscr{M}\),
 the \textsf{Weyl polynomial \(W_{\mathscr{M}}^{\,p}\)
  of the index \(p,\,\,p=1,\,2,\,3,\,\dots\,\,,\)
 for \(\mathscr{M}\)}
 is \textsf{defined} as
\begin{equation}
\label{UWP} %
W_{\mathscr{M}}^{\,p}(t)=
\sum\limits_{l=0}^{[\frac{n}{2}]}%
\frac{2^{-l}\,\Gamma(\frac{p}{2}+1)}{\Gamma(\frac{p}{2}+l+1)}%
k_{2l}(\mathscr{M}) \cdot{}t^{2l}\,,
\end{equation}
where the Weyl coefficients \(k_{2l}(\mathscr{M})\) are introduced
in \textup{Definition \ref{DNWP}}.
\end{definition}
Let us emphasize that in Definition \ref{DGWP} no assumption
concerning the smo\-othness of the surface  \(\mathscr{M}\) are
made. We already mentioned that the definitions of the Weyl
coefficients \(k_{2l}\) for smooth manifolds and for convex
surfaces agree. Therefore, Definitions \ref{DEfWP}\,-\,\ref{DeWC}:
\eqref{deftub}-\eqref{Soo}-\eqref{DefWC} of the Weyl polynomial
and the Weyl coefficients for a smooth manifold and Definition
\ref{DGWP} of the Weyl polynomials for a closed convex surface
agree if the convex surface is also a smooth manifold.

We also define the \(W_{\mathscr{M}}^{\,\infty}\) of the infinite
index.
\begin{definition}%
\label{deLWP}%
Let \(\mathscr{M},\,\dim{}\mathscr{M}=n\) be either a smooth
manifold, or a closed compact convex surface, and let
\(k_{2l}(\mathscr{M}),\,l=0,\,1,\,\ldots\,,\,[\frac{n}{2}]\),  be
the Weyl coefficients of \(\mathscr{M}\), defined by
\textup{Definition \ref{DeWC}} in the smooth case, and by
\textup{Definition \ref{DNWP}} in the convex case. The Weyl
polynomial of the infinite index \(W_{\mathscr{M}}^{\,\infty}\) is
defined as
\begin{equation}%
\label{DeWPInfInd}%
W_{\,\mathscr{M}}^{\infty}(t)=\sum\limits_{l=0}^{[\frac{n}{2}]}%
k_{2l}(\mathscr{M})\cdot{}t^{2l}.
\end{equation}%
\end{definition}%
\begin{remark}
\label{InfLimC}%
In view of \eqref{DecFor},
\[W_{\mathscr{M}}^{\,p}(\sqrt{p}t)=k_{0}(\mathscr{M})+\sum\limits_{l=1}^{[\frac{n}{2}]}%
\frac{p^l}{(p+2)(p+4)\,\cdots\,\,(p+2l)}\,
k_{2l}(\mathscr{M})\cdot{}t^{2l}\,.\] Therefore, the polynomial
\(W_{\,\mathscr{M}}^{\infty}(t)\) can be considered as a limiting
object for the family
\(\big\{W_{\mathscr{M}}^{\,p}(t)\big\}_{p=1,\,2,\,3,\,\dots\,}\)
of the Weyl polynomials of the index \(p\):
\begin{equation}
\label{LiRe}
W_{\,\mathscr{M}}^{\infty}(t)=\lim_{p\to\infty}W_{\,\mathscr{M}}^{p}(\sqrt{p}t)\,.
\end{equation}
\end{remark}

\begin{center}
\begin{minipage}{0.93\linewidth} \textsl{Thus, the sequence
\begin{math}
\label{SeWePo}\big\{W_{\mathscr{M}}^{\,p}\big\}_{p=1,\,2,\,3,\,\dots\,}
\end{math}
of the Weyl polynomials,
\begin{math}%
\deg{}W_{\mathscr{M}}^{p}= 2{\textstyle\left[\frac{n}{2}\right]}
\end{math}, %
as well as the `limiting' polynomial
\(W_{\,\mathscr{M}}^{\infty}\) are related  to any closed
\(n\)\,-\,dimensional convex surface \(\mathscr{M}\).}
\end{minipage}
\end{center}

\vspace{2.0ex}%
\textsf{ Weyl polynomials (and Minkowski polynomials
 in the convex case) reflect somehow intrinsic properties of the
 appropriate manifolds. On the other hand, there are known very
 distinguished and remarkable geometrical objects such as regular
 polytopes, compact matricial groups, spaces of constant
 curvatures, etc. Our belief is that the Weyl polynomials related
 to these geometric objects are of fundamental importance and
 possess interesting properties. These
 polynomials should be carefully studied. In particular, the
 following question is natural:\\
 \centerline{\textit{What can we say about roots of such polynomials?}}}\\
\section{FORMULATION OF MAIN RESULTS.\label{FMR}}
In this section we formulate the main results of this paper about
location of the roots of the Minkowski and Weyl  polynomials
related to convex sets and surfaces. \\[-5.0ex]
\paragraph{Dissipative and conservative polynomials.\label{DCPol}}
 We introduce two classes of polynomials: dissipative polynomials
 and conservative polynomials.
In many cases the Minkowski polynomials related to convex sets
 are dissipative, and the Weyl polynomials are conservative.

\begin{definition}
\label{DeHuPo}
The polynomial \(M\) is said to be \textsf{dissipative} if all roots
of \(M\) are situated in the open left half plane \(\{z:\textup{Re}\,z<0\}.\)
The dissipative polynomials are also called the \textsf{Hurwitz polynomials},
or the \textsf{stable polynomials}.
\end{definition}
\begin{definition}
\label{DeCoPo} The polynomial \(W\) is said to be
\textsf{conservative} if all roots of \(W\) are purely imaginary and
simple, in other words if all roots of \(W\) are contained in the
imaginary axis \(\{z:\textup{Re}\,z=0\}\), and each of them is of
multiplicity one.
\end{definition}

\begin{theorem} %
\label{H10H0}%
Given a closed compact convex surface \(\mathscr{M}\), \(\dim
\mathscr{M}=n\) \(\mathscr{M}=\partial{}V\),
 let \(W_{\mathscr{M}}^{\,1}\) be the Weyl polynomial of index
 \(1\) related to \(\mathscr{M}\), and let \(M_{\,\,V}^{\mathbb{R}^{n+1}}\) be the Minkowski
 polynomial related to the set \(V\).\\[0.8ex]
\hspace*{2.0ex}If the polynomial \(M_{\,\,V}^{\mathbb{R}^{n+1}}\)
 is dissipative, then the polynomial
\(W_{\mathscr{M}}^{1}\) is conservative.
\end{theorem}
The proof of Theorem \ref{H10H0} is based on the relation
\eqref{WMP}. Theorem \ref{H10H0} is derived from \eqref{WMP} using
Hermite-Biehler theorem. We do this in Section \ref{HBieT}.

From \eqref{LiRe} it follows that if for every \(p\) the
polynomial \(W_{\mathscr{M}}^{\,p}\) has only purely imaginary
roots, than all the roots of the polynomial
\(W_{\mathscr{M}}^{\,\infty}\) are purely imaginary as well. In
particular, \textit{all the roots of the polynomial
\(W_{\mathscr{M}}^{\,\infty}\) are purely imaginary if for every
\(p\) the polynomial \(W_{\mathscr{M}}^{\,p}\) is conservative}.

However, what is important for us that is the converse statement:
\begin{lemma}
\label{LWPo} If the polynomial \(W_{\mathscr{M}}^{\,\infty}\) is
conservative, then all the polynomials
\(W_{\mathscr{M}}^{\,p},\,p=1,\,2,\,3,\,\dots\,\,,\) are
conservative as well.
\end{lemma}

Lemma \ref{LWPo} is the consequence of some Laguerre result about
the multiplier sequences. Proof of Lemma \ref{LWPo} appeares in
the end of Section \ref{PEFG}.

Keeping in mind Lemma \ref{LWPo},  we will concentrate our efforts
on the study of the location of the roots of the Weyl polynomial
\(W_{\mathscr{M}}^{\,\infty}\) of the infinite index.

\paragraph{The case of low dimension.}
In this section we discuss the Minkowski polynomials of convex
sets \(V,\,V\subset{}\mathbb{R}^n,\) and the Weyl polynomials of
closed convex surfaces \(\mathscr{M}\), \(\dim V=n,\) for `small'
\(n\): \(n=2,\,3,\,4,\,5\).

\begin{theorem}
\label{LDC}%
 Let \(n\) be one of the numbers \(2,\,3,\,4\) or
\(5\), and let \(V,\,V\subset{}\mathbb{R}^n,\) be a solid compact
convex set. Then the Minkowski polynomial \(M_{V}^{\mathbb{R}^n}\)
is dissipative.
\end{theorem}
\begin{theorem}
\label{LDCW}%
 Let \(n\) be one of the numbers \(2,\,3,\,4\) or \(5\), and let
\(\mathscr{M}\) be closed proper\,%
\footnote{That is the generating set \(V\) is solid.} %
 convex surface of dimension
\(n\).\\[1.0ex]
\hspace*{1.0ex}Then: \\[-4.0ex]
\begin{enumerate}%
\item%
The Weyl polynomial \(W_{\mathscr{M}}^{\infty}\) of  infinite
index is conservative.
\item
For every \(p=1,\,2,\,3,\,\dots\,\,,\,\)  the Weyl polynomial
\(W_{\mathscr{M}}^{\,p}\)\,\ of index \(p\) %
 is conservative.
\end{enumerate}
\end{theorem}

Theorem \ref{LDC} and \ref{LDCW} are proved in section \ref{LoDi}.
Proving these theorems, we combine the Routh-Hurwitz criterion,
which express the property of a polynomial to be dissipative in
terms of its coefficients, and the Alexandrov-Fenchel
inequalities, which express the logarithmic convexity property for
the sequence of the cross-sectional measures of a convex set.
\paragraph{Selected 'regular' convex sets: balls, cubes, squeezed cylinders.}
For large \(n\), the statements analogous to Theorems \ref{LDC}
and \ref{LDCW} do not hold. If \(n\) is large enough, then there
exists such solid compact convex sets%
\footnote{Very flattened ellipsoids  can be taken as such V. See Theorem \ref{NMR}.} %
 \(V\), \(\dim V=n\), that
Min\-kow\-ski polynomials \(M_{\,V}^{\mathbb{R}^{n+1}}\) are not
dissipative, and the Weyl polynomials \(W^{p}_{\partial{}V}\) are
not conservative. However, for some `regular' convex sets \(V\),
like balls and cubes, the Weyl polynomials \(W^{p}_{\partial{}V}\)
are conservative, and the Minkowski polynomial are dissipative in
any dimension.

Let us present the collection of `regular' convex sets and their
boundary surfaces which we are dealing with further. Such sets and
surfaces will be considered for every \(n\).
\begin{itemize}%
\item[\(\Diamond\)]
The unit ball \(B^n\):
\begin{gather}
\label{Dub}%
 B^n=\{x=(x_1,\,\ldots\,,\,x_n)\in\mathbb{R}^{n}:
\sum\limits_{1\leq{}k\leq{}n}|x_k|^2\leq{}1\,\},\\
\textup{Vol}_n(B^n)=\frac{\pi^{n/2}}{\Gamma(\frac{n}{2}+1)}\,.
\end{gather}
\item[\(\Diamond\)]
The squeezed spherical cylinder \(B^{n}\times{}0\),
\(B^{n}\times{}0\subset\mathbb{R}^{n+1}\).
\item[\(\Diamond\)]
The unit sphere, \[
S^n=\{x=(x_1,\,\ldots\,,\,x_n,\,x_{n+1})\in\mathbb{R}^{n+1}:
\sum\limits_{1\leq{}k\leq{}n+1}|x_k|^2=1\,\},\]
 in other words, the boundary surface of the
unit ball: \(S^n=\partial{}B^{n+1}\,,\)
\begin{gather}
\label{Dus}%
\textup{Vol}_n(S^n)=(n+1)\,\textup{Vol}_{n+1}(B^{n+1})\,.
\end{gather}
\item[\(\Diamond\)]
 The boundary surface of the squeezed spherical cylinder
\(\partial{}(B^{n}\times{}0)\) :
\begin{equation}
\label{SCub}%
\textup{Vol}_n(\partial{}(B^{n}\times{}0))=2\,\textup{Vol}_n(B^n)\,.
\end{equation}
\item[\(\Diamond\)]
The unit cube \(Q^n\):
\begin{gather}
\label{Duc}%
 Q^n=\{x=(x_1,\,\ldots\,,\,x_n)\in\mathbb{R}^{n}:
\max\limits_{1\leq{}k\leq{}n}|x_k|\leq{}1\,\},\\
 \textup{Vol}_n(Q^n)=2^n\,.
\end{gather}
\item[\(\Diamond\)]
The squeezed cubic cylinder \(Q^{n}\times{}0\),
\(Q^{n}\times{}0\subset\mathbb{R}^{n+1}\).
\item[\(\Diamond\)]
The boundary surface \(\partial{}Q^{n+1}\) of the unit cube:
\begin{equation}
\label{Bsuc}%
\textup{Vol}_n(\partial{}Q^{n+1})=(n+1)\,\textup{Vol}_{n+1}(Q^{n+1}).
\end{equation}
\item[\(\Diamond\)]
The boundary surface of the squeezed cubic cylinder
\(\partial{}(Q^{n}\times{}0)\):
\begin{equation}
\label{SqCub}%
\textup{Vol}_n(\partial{}(Q^{n}\times{}0))=2\,\textup{Vol}_{n}(Q^{n})\,.
\end{equation}
\end{itemize}
\paragraph{The location of roots of the Minkowski and Weyl
polynomials\\
related to the `regular' convex sets.\\[1.0ex]} %
Let us state the main results about location of roots of the
Minkowski polynomials and the Weyl polynomials related to the
above mentioned `regular' convex sets and their surfaces.
\begin{theorem}%
\label{LoMP}%
For every \(n=1,\,2,\,3,\,\,\ldots\)\,\textup{:}\\[-5.0ex]
\begin{enumerate}%
\item%
The Minkowski polynomial \(M_{B^n}^{\mathbb{R}^n}\) related to the
ball \(B^n\) is dissipative, moreover all its roots are
negative\,%
\footnote{This part of the Theorem is trivial: \(M_{B^n}^{\mathbb{R}^n}(t)=(1+t)^n\)}. %
\item%
The Minkowski polynomial
\,\(M_{B^{n}\times{}0}^{\mathbb{R}^{n+1}}\) \,related to the
squeezed spherical cylinder \(B^{n}\times{}0\) is
of the form%
\footnote{\,\label{FaAp}The factors \(t\) appears
because the set \(B^{n}\times{}0\)
is not solid in \(\mathbb{R}^{n+1}\).} %
 \(M_{B^{n}\times{}0}^{\mathbb{R}^{n+1}}(t)=
t\cdot{}D_{B^{n}\times{}0}^{\mathbb{R}^{n+1}}(t)\), where the
polynomial \(D_{B^{n}\times{}0}^{\mathbb{R}^{n+1}}\) is
dissipative. If \(n\) is large enough, then the polynomial
\(M_{B^{n}\times{}0}^{\mathbb{R}^{n+1}}\) has non-real roots.
\item%
The Minkowski polynomial \(M_{Q^n}^{\mathbb{R}^n}\) related to
cube \(Q^n\) is dissipative, moreover all its roots are negative.
\item%
The Minkowski polynomial \(M_{Q^{n}\times{}0}^{\mathbb{R}^{n+1}}\)
related to the squeezed cubical cylinder \(Q^{n}\times{}0\) is of
the form\({}^{\ref{FaAp}}\)
\(M_{Q^{n}\times{}0}^{\mathbb{R}^{n+1}}(t)=
t\cdot{}D_{Q^{n}\times{}0}^{\mathbb{R}^{n+1}}(t)\), where the
polynomial \(D_{Q^{n}\times{}0}^{\mathbb{R}^{n+1}}\) is
dissipative, moreover all roots of the polynomial
\(D_{Q^{n}\times{}0}^{\mathbb{R}^{n+1}}\) are negative.
\end{enumerate}%
\end{theorem}
\begin{theorem}%
\label{LoWP}%
For every \(n=1,\,2,\,3,\,\,\ldots\)\,\textup{:}\\[-5.0ex]
\begin{enumerate}%
\item%
The Weyl polynomials \(W_{\partial{}B^{n+1}}^{\,\infty}(t)\) of
infinite index, as well as the Weyl polynomials
\(W_{\partial{}B^{n+1}}^{\,p}(t)\) of arbitrary finite index
\(p,\,p=1,\,2,\,\ldots\,\),
 related to
the boundary surface of the ball \(B^{n+1}\) are conservative.
\item%
The Weyl polynomials \(W_{\partial{}(B^{n}\times{}0)}^{\,p}\) of
order\,%
\footnote{\,The case \(p=3\) remains open.} %
 \(p=1,\,p=2\) and \(p=4\) related to the boundary surface
of the squeezed spherical cylinder \(B^{n}\times{}0\) are
conservative.
\item%
The Weyl polynomials \(W_{\partial{}Q^{n+1}}^{\,\infty}(t)\) of
infinite index, as well as the Weyl polynomials
\(W_{\partial{}Q^{n+1}}^{\,p}(t)\) of arbitrary finite index
\(p,\,p=1,\,2,\,\ldots\,\),
 related to
the boundary surface of the cube \(Q^{n+1}\) are conservative.
\item%
The Weyl polynomials
\(W_{\partial{}(Q^{n}\times{}0)}^{\,\infty}(t)\) of infinite
index, as well as the Weyl polynomials
\(W_{\partial{}(Q^{n}\times{}0)}^{\,p}(t)\) of arbitrary finite
index \(p,\,p=1,\,2,\,\ldots\,\),
 related to
the boundary surface of the squeezed cubic cylinder
\(Q^{n}\times{}0\) are conservative.
\end{enumerate}%
\end{theorem}
\begin{remark}
\label{EERW} The roots of the Weyl polynomial
\(W^1_{\partial{}B^{n+1}}\) can be found explicitly. Indeed
\[W^1_{\partial{}B^{n+1}}(it)=\textup{Vol}_{n+1}(B^{n+1})\frac{1}{2it}\big((1+it)^{n+1}-(1-it)^{n+1}\big)\,.\]
Changing variable \[t\to\varphi:\,1+it=|1+it|e^{i\varphi},
t=\tg{}\varphi\,,\ \ -\frac{\pi}{2}<\varphi<\frac{\pi}{2}\,,\] we
reduce the equation \(W^1_{\partial{}B^{n+1}}(it)=0\) to the
equation
\[\frac{\sin{}(n+1)\varphi}{\sin{}\varphi}=0,\ \ -\frac{\pi}{2}<\varphi<\frac{\pi}{2}\,.\]
The roots of the latter equation are:
\[\varphi_k=\frac{k\pi}{n+1},
\quad{}-\left[\frac{n}{2}\right]\leq{}k\leq{}\left[\frac{n}{2}\right],\
\ k\not=0\,.\]
 So, the roots \(t_k\) of the equation
\(W^{1}_{\partial{}B^{n+1}}(it)=0\) are
\[t_k=\tg{}\textstyle{\frac{k\pi}{n+1},\quad{}
-\left[\frac{n}{2}\right]\leq{}k\leq{}\left[\frac{n}{2}\right],\ \
k\not=0\,.}\] In particular, the polynomial \(W^{1}_{S^n}\) is
conservative.
\end{remark}
\textsf{Negative results}:
\begin{theorem}
\label{NRWPSq}%
 Given an integer \(p\), \(p\geq{}5\). If \(n\) is
large enough: \(n\geq{}N(p)\), then the Weyl polynomial
\(W^p_{\partial(B^n\times{}0)}\) is not conservative: some of its
roots do not belong to the imaginary axis.
\end{theorem}
 For an integer \(q:q\geq{}1\),  let
\(E_{n,\,q,\,\varepsilon}\) be the \(n+q\)-\,dimensional
ellipsoid:
\begin{subequations}
\label{Ell}
\begin{equation}
\label{Ell1}
E_{n,\,q,\,\varepsilon}=\{(x_1,\,x_2,\,\ldots\,,\,x_n,\,\ldots\,,\,x_{n+q})\in\mathbb{R}^{n+q}:
\sum\limits_{0\leq{}j\leq{}n+q}(x_j/a_j)^2\leq{}1\},%
\end{equation}%
where
\begin{equation}
\label{Ell2} a_j=1\ \ \textup{for}\ \ 1\leq{}j\leq{}n,\ \
a_j=\varepsilon\ \
\textup{for}\ \ n+1\leq{}j\leq{}n+q\,.%
\end{equation}%
\end{subequations}

\begin{theorem}\ \ \vspace*{-3.0ex}%
\label{NMR}%
\begin{enumerate}
\item Given an integer \(q:\,\ 5\leq{}q<\infty\). If \(n\) is large
enough:\,\(n\geq{}N(q)\), and \(\varepsilon\) is small enough:
\(0<\varepsilon\leq\varepsilon(n,\,q)\), then the Minkowski
polynomial \(M_{E_{n,\,q,\,\varepsilon}}^{\mathbb{R}^{n+q}}\) is
not dissipative: some of its roots are situated in the open
right-half plane.
\item
Given an integer \(p\) and an odd integer \(q\): \,
\(p\geq{}1,\,q\geq{}1,\,p+q\geq{}6\,\).  %
If \,\,\(n\) is large enough:\,\(n\geq{}N(p,q)\), and
\(\varepsilon\) is small enough:
\(0<\varepsilon\leq\varepsilon(n,\,p,\,q)\), then the Weyl
polynomial \, \(W^p_{E_{n,\,q,\,\varepsilon}}\) \,is not
conservative: some of its roots do not belong to the imaginary
axis.
\end{enumerate}
\end{theorem} Proof of Theorem \ref{NMR} is presented in
Section \ref{ExInSp}.
\section{THE EXPLICIT EXPRESSIONS\\ FOR THE MINKOWSKI AND WEYL
POLYNOMIALS\\ RELATED TO THE `REGULAR' CONVEX SETS.\label{EEWMP}}
Hereafter, we use the following identity for the
\(\Gamma\)-\,function:
\begin{equation}
\label{IfGa} \Gamma(\zeta+1/2)\,
\Gamma(\zeta+1)={\pi}^{1/2}2^{-2\zeta}\Gamma(2\zeta+1)\,,\ \forall
\zeta\in\mathbb{C}: 2\zeta\not=-1,\,-2,\,-3,\,\ldots\,\,.
\end{equation}
Let as present explicit expressions for the Minkowski polynomials
related to the `regular' convex sets: balls, cubes, squeezed
cylinders, as well as the expression for the Weyl polynomials
related to the boundary surfaces of these sets. The items related
to balls are marked by the symbol \(\bigodot\), the items related
to cubes are marked by the symbol \(\boxdot\).
\\[0.5ex]
\hspace*{0.0ex}\(\bigodot\)\hspace*{0.5ex}\textsf{The unit ball
\(B^n\)}.\\
Since \(B^n+tB^n=(1+t)B^n\) for \(t>0\), then, according to
\eqref{DeMiPo},
\begin{equation}
\label{EMPb}
M_{B^n}^{\mathbb{R}^n}(t)=%
\textup{Vol}_n(B^n)\cdot(1+t)^n\,,
\end{equation}
or
\begin{equation}
\label{EMPb1}
M_{B^n}^{\mathbb{R}^n}(t)=%
\textup{Vol}_n(B^n)\sum\limits_{0\leq{}k\leq{}n}\frac{n!}{(n-k)!}\cdot{}\frac{t^k}{k!}\,.
\end{equation}
\hspace*{2.0ex}Thus, the coefficients of the Minkowski polynomial
\(M_{B^n}^{\mathbb{R}^n}\) for the ball \(B^n\) are:
\begin{equation}
\label{MiCB}
m_{\,k}^{\mathbb{R}^n}(B^n)=\textup{Vol}_n(B^n)\cdot{}
\frac{n!}{(n-k)!}\cdot{}\frac{1}{k!}\,,\quad{}0\leq{}k\leq{}n\,.
\end{equation}
\hspace*{0.0ex}\(\bigodot\)\hspace*{0.5ex}\textsf{The squeezed
spherical
cylinder \(B^{n}\times{}0\)}.\\
The Minkowski polynomial for the squeezed spherical cylinder
\(B^{n}\times{}0\) is:
\begin{equation}
\label{MiSqB} M_{B^{n}\times{}0}^{\mathbb{R}^{n+1}}(t)=
\textup{Vol}_n(B^n)\cdot\!t\!\!\sum\limits_{0\leq{}k\leq{}n}\frac{n!}{(n-k)!}%
\frac{{\pi}^{1/2}\Gamma(\frac{k}{2}+1)}{\Gamma(\frac{k+1}{2}+1)}\frac{1}{k!}\,
t^k\,.
\end{equation}
The expression \eqref{MiSqB} is derived from \eqref{EMPb1} and
\eqref{IRP1}-\eqref{MPn1}. (See Lemma \ref{IPR}.)\\
\hspace*{1.5ex}Thus, the coefficients of the Minkowski polynomial
\(M_{B^n\times{}0}^{\mathbb{R}^{n+1}}\) for the squeezed spherical
cylinder \(B^{n}\times{}0\) are:
\begin{multline}
\label{CoMiSqB}%
 m_{\,0}^{\mathbb{R}^{n+1}}(B^n\times{}0)=0,\quad
m_{\,k+1}^{\mathbb{R}^{n+1}}(B^n\times{}0)=\\
=\textup{Vol}_n(B^n)\cdot{}\frac{n!}{(n-k)!}\cdot{}
\frac{{\pi}^{1/2}\Gamma(\frac{k}{2}+1)}%
{\Gamma(\frac{k+1}{2}+1)}\frac{1}{k!}\,,\quad{}0\leq{}k\leq{}n\,.
\end{multline}

\hspace*{0.0ex}\(\bigodot\)\hspace*{0.5ex}\textsf{The unit sphere
\(S^n=\partial{}B^{n+1}\)}.\\
According to \eqref{NWP} and \eqref{MiCB}, the Weyl coefficients
of the \(n\)-\,dimensional sphere \(S^n=\partial{}B^{n+1}\) are:
\begin{equation}
\label{WCSp} k_{2l}(\partial{}B^{n+1})=\textup{Vol}_{n}
(\partial{}B^{n+1})\cdot\frac{n!}{(n-2l)!}\cdot\frac{1}{l!}\frac{1}{2^l}\,,
\quad{}0\leq{}l\leq[{\textstyle\frac{n}{2}}]\,.
\end{equation}
Thus, the Weyl polynomials related to the \(n\)-\,dimensional
sphere are:
\begin{multline}
\label{WPpS} %
W_{\partial{}B^{n+1}}^{\,p}(t)=\textup{Vol}_{n}
(\partial{}B^{n+1})\cdot{}\\
\cdot{}\sum\limits_{l=0}^{[\frac{n}{2}]}%
\frac{n!}
{(n-2l)!}\cdot{}\frac{2^{-l}\,\Gamma(\frac{p}{2}+1)}{\Gamma(\frac{p}{2}+l+1)}
\cdot{}\frac{1}{l!}\cdot{}\Big(\frac{t^2}{2}\Big)^l,\ \
p=1,\,2,\,\ldots\,.
\end{multline}
\begin{equation}
 \label{WPiS}
 W_{\partial{}B^{n+1}}^{\,\infty}(t)=\textup{Vol}_{n}
(\partial{}B^{n+1})
\cdot{}\sum\limits_{l=0}^{[\frac{n}{2}]}%
\frac{n!}
{(n-2l)!}\cdot{}\frac{1}{l!}\cdot{}\Big(\frac{t^2}{2}\Big)^l\,\cdot
\end{equation}
\hspace*{0.0ex}\(\bigodot\)\hspace*{0.5ex}\textsf{The boundary
surface \(\partial(B^{n}\times{}0)\) of the
squeezed spherical cylinder \(B^{n}\times{}0\)}.\\
According to \eqref{WCSqC} and \eqref{MiCB}, the Weyl coefficients
of the \(n\)-\,dimensional improper surface
\(\partial(B^{n}\times{}0)\) are:
\begin{equation}
\label{WCSqB} k_{2l}(\partial(B^{n}\times{}0))=\textup{Vol}_{n}
(\partial(B^{n}\times{}0))\cdot\frac{n!}{(n-2l)!}\cdot\frac{\Gamma(1/2)}{\Gamma(l+1/2)}\,\frac{1}{2^l},
\quad{}0\leq{}l\leq[{\textstyle\frac{n}{2}}]\,.
\end{equation}
Thus, the Weyl polynomials related to the (improper) surface
\(\partial(B^{n}\times{}0)\) are:
\begin{multline}
\label{WPpSSqC} %
W_{\partial{}(B^{n}\times{}0)}^{\,p}(t)=\textup{Vol}_{n}
(\partial{}(B^{n}\times{}0))\cdot{}\\
\cdot{}\sum\limits_{l=0}^{[\frac{n}{2}]}%
\frac{n!}
{(n-2l)!}\cdot{}\frac{2^{-l}\,\Gamma(\frac{p}{2}+1)}{\Gamma(\frac{p}{2}+l+1)}
\cdot{}\frac{\Gamma(1/2)}{\Gamma(l+1/2)}\cdot{}\Big(\frac{t^2}{2}\Big)^l,\
\ p=1,\,2,\,\ldots\,.
\end{multline}
\begin{equation}
 \label{WPiSSqC}
 W_{\partial{}(B^{n+1}\times{}0)}^{\,\infty}(t)=\textup{Vol}_{n}
(\partial{}(B^{n+1}\times{}0))
\cdot{}\sum\limits_{l=0}^{[\frac{n}{2}]}%
\frac{n!}
{(n-2l)!}\cdot{}\frac{\Gamma(1/2)}{\Gamma(l+1/2)}
\cdot{}\Big(\frac{t^2}{2}\Big)^l\,\cdot
\end{equation}
\hspace*{0.0ex}\(\boxdot\)\hspace*{0.5ex}\textsf{The unit cube
\(Q^n\)}.\\
The Minkowski polynomial \(M_{Q^n}^{\mathbb{R}^n}\) is:
\begin{equation}
\label{MPUQ}%
 M_{Q^n}^{\mathbb{R}^n}(t)= \textup{Vol}_{n}%
(Q^{n})\cdot{}\!\!\!\!\sum\limits_{0\leq{}k\leq{}n}\frac{n!}{(n-k)!}
\frac{1}{\Gamma(\frac{k}{2}+1)k!}\Big(\frac{\sqrt{\pi}}{2}\Big)^k\,t^k\,.
\end{equation}
The expression \eqref{MPUQ} is obtained in the following way. The
\(n\)-\,dimensional cube \(Q^{n}\) is considered as the Cartesian
product of the one-dimensional cubes:
\begin{equation*}
\label{CsCP} Q^{n}=Q^1\times\,\cdots\,Q^1\,.
\end{equation*}
For \(n=1\), the Minkowski polynomial  is:
\(M_{Q^1}^{\mathbb{R}^1}(t)=2(1+t)\)\,. Then the rule is used how
to express the Minkowski polynomial of the Cartesian product in
terms of the Minkowski polynomials for the Cartesian factors. (See
details in Section \ref{MPCaPr}.)\\ %
\hspace*{2.0ex} Thus, the coefficients of the Minkowski polynomial
for the cube \(Q^n\) are:
\begin{equation}
\label{MiCQ}
m_{\,k}^{\mathbb{R}^n}(Q^n)=\textup{Vol}_n(Q^n)\cdot{}
\frac{n!}{(n-k)!}\cdot{}\frac{1}{\Gamma(\frac{k}{2}+1)k!}
\Big(\frac{\sqrt{\pi}}{2}\Big)^k\,,\quad{}0\leq{}k\leq{}n\,.
\end{equation}
\hspace*{0.0ex}\(\boxdot\)\hspace*{0.5ex}\textsf{The squeezed
cubic cylinder \(Q^n\times{}0\)}.\\
The Minkowski polynomial \(M_{Q^{n}\times{}0}^{\mathbb{R}^{n+1}}\)
is:
\begin{equation}
\label{MPUQC}%
 M_{Q^{n}\times{}0}^{\mathbb{R}^{n+1}}(t)= \textup{Vol}_{n}%
(Q^{n})\cdot{}t\!\!\!\!\sum\limits_{0\leq{}k\leq{}n}\frac{n!}{(n-k)!}
\frac{\Gamma(\frac{1}{2})}{\Gamma(\frac{k+1}{2}+1)k!}\Big(\frac{\sqrt{\pi}}{2}\Big)^k\,t^k\,.
\end{equation}
The expression \eqref{MPUQC} is derived from \eqref{MPUQ} and
\eqref{IRP1}-\eqref{MPn1}. (See Lemma \ref{IPR}.)\\ %
\hspace*{2.0ex}Thus, the coefficients of the Minkowski polynomial
\(M_{Q^n\times{}0}^{\mathbb{R}^{n+1}}\) for the squeezed cubic
cylinder are:
\begin{multline}
\label{CoMiSqQ} m_{\,0}^{\mathbb{R}^{n+1}}(Q^n\times{}0)=0,\quad
m_{\,k+1}^{\mathbb{R}^{n+1}}(Q^n\times{}0)=\\
=\textup{Vol}_n(Q^n)\cdot{}\frac{n!}{(n-k)!}\cdot{}
\frac{\Gamma(\frac{1}{2})}%
{\Gamma(\frac{k+1}{2}+1)}\frac{1}{k!}\Big(\frac{\sqrt{\pi}}{2}
\Big)^k\,,\quad{}0\leq{}k\leq{}n\,.
\end{multline}
\hspace*{0.0ex}\(\boxdot\)\hspace*{0.5ex}\textsf{The boundary
surface \(\partial{}Q^{n+1}\) of the unit cube \(Q^{n+1}\)}.\\
According to \eqref{NWP} and \eqref{MiCQ}, the Weyl coefficients
of the \(n\)-\,dimensional surface \(\partial{}Q^{n+1}\) are:
\begin{multline}
\label{WCQp} k_{2l}(\partial{}Q^{n+1})=\textup{Vol}_{n}
(\partial{}Q^{n+1})\cdot\frac{n!}{(n-2l)!}\cdot\\
\cdot\frac{1}{\Gamma(l+\frac{1}{2}+1)}\frac{1}{2^l\,l!}\,\Big(\frac{\sqrt{\pi}}{2}\Big)^{2l+1},
\quad{}0\leq{}l\leq[{\textstyle\frac{n}{2}}]\,.
\end{multline}
Taking into account the identity
\(\Gamma(l+1+{\textstyle\frac{1}{2}})\cdot\Gamma(l+1)={\pi}^{1/2}2^{-(2l+1)}\Gamma(2l+2)\),
which is the identity \eqref{IfGa} for \(\zeta=l+1/2\), we can
transform \eqref{WCQp} to the form
\begin{equation}
\label{aWCQp} k_{2l}(\partial{}Q^{n+1})=\textup{Vol}_{n}
(\partial{}Q^{n+1})\cdot\frac{n!}{(n-2l)!}
\cdot\frac{1}{(2l+1)!}\,\Big(\frac{\pi}{2}\Big)^{l},
\quad{}0\leq{}l\leq[{\textstyle\frac{n}{2}}]\,.
\end{equation}
Thus, the Weyl polynomials related to the \(n\)-\,dimensional
surface \(\partial{}Q^{n+1}\) are:
\begin{multline}
\label{WPpQ} %
W_{\partial{}Q^{n+1}}^{\,p}(t)=\textup{Vol}_{n}
(\partial{}Q^{n+1})\cdot{}\\
\cdot{}\sum\limits_{l=0}^{[\frac{n}{2}]}%
\frac{n!} {(n-2l)!}\cdot{}\frac{
2^{-l}\,\Gamma(\frac{p}{2}+1)}{\Gamma(\frac{p}{2}+l+1)}
\cdot{}\frac{1}{(2l+1)!}\cdot{}\Big({\frac{{\pi}t^2}{2}}\,\Big)^l,\
\ p=1,\,2,\,\ldots\,.
\end{multline}
\begin{equation}
 \label{WPiQ}
 W_{\partial{}Q^{n+1}}^{\,\infty}(t)=\textup{Vol}_{n}
(\partial{}Q^{n+1})
\cdot{}\sum\limits_{l=0}^{[\frac{n}{2}]}%
\frac{n!}
{(n-2l)!}\cdot{}\frac{1}{(2l+1)!}\cdot{}\Big(\frac{{\pi}t^2}{2}\Big)^l\,\cdot
\end{equation}
\hspace*{0.0ex}\(\boxdot\) \hspace*{0.5ex}\textsf{The boundary
surface \(\partial(Q^{n}\times{}0)\) of the
squeezed cubic cylinder \(Q^{n}\times{}0\)}.\\
According to \eqref{WCSqC} and \eqref{MiCQ}, the Weyl coefficients
of the surface (improper) \(\partial{}(Q^{n}\times{}0)\) are:
\begin{multline}
\label{WCSQp} k_{2l}(\partial{}Q^{n}\times{}0)=\textup{Vol}_{n}
(\partial{}Q^{n}\times{}0)\cdot\frac{n!}{(n-2l)!}\cdot\\
\cdot\frac{\sqrt{\pi}}{\Gamma(l+\frac{1}{2})}\frac{1}{l!\,2^{l}\,}\,\Big(\frac{\pi}{2}\Big)^{l},
\quad{}0\leq{}l\leq[{\textstyle\frac{n}{2}}]\,.
\end{multline}
Using the identity
\(\Gamma(l+1/2)\Gamma(l+1)=\sqrt{\pi}2^{-2l}\Gamma(2l+1)\), which
is the identity \eqref{IfGa} for \(\zeta=l\), the equality
\eqref{WCSQp} can be transformed to the form
\begin{equation}
\label{aWCSQp} k_{2l}(\partial{}Q^{n}\times{}0)=\textup{Vol}_{n}
(\partial{}Q^{n}\times{}0)\cdot\frac{n!}{(n-2l)!}\cdot
\frac{1}{(2l)!}\,\Big(\frac{\pi}{2}\Big)^{2l},
\quad{}0\leq{}l\leq[{\textstyle\frac{n}{2}}]\,.
\end{equation}
Thus, the Weyl polynomials related to the improper
\(n\)-\,dimensional surface \(\partial{}(Q^{n}\times{}0)\) are:
\begin{multline}
\label{WPpQSqC} %
W_{\partial{}(Q^{n}\times{}0)}^{\,p}(t)=\textup{Vol}_{n}
(\partial{}(Q^{n}\times{}0))\cdot{}\\
\cdot{}\sum\limits_{l=0}^{[\frac{n}{2}]}%
\frac{n!} {(n-2l)!}\cdot{}
\frac{2^{-l}\,\Gamma(\frac{p}{2}+1)}{\Gamma(\frac{p}{2}+l+1)}
\cdot{}\frac{1}{(2l)!}\cdot{}\Big({\frac{{\pi}t^2}{2}}\,\Big)^l,\
\ p=1,\,2,\,\ldots\,.
\end{multline}
\begin{equation}
 \label{WPiQSqC}
 W_{\partial{}(Q^{n}\times{}0)}^{\,\infty}(t)=\textup{Vol}_{n}
(\partial{}(Q^{n}\times{}0))
\cdot{}\sum\limits_{l=0}^{[\frac{n}{2}]}%
\frac{n!}
{(n-2l)!}\cdot{}\frac{1}{(2l)!}\cdot{}\Big(\frac{{\pi}t^2}{2}\Big)^l\,\cdot
\end{equation}

\section{WEYL AND MINKOWSKI POLYNOMIALS\\ OF `REGULAR' CONVEX
SETS\\
AS RENORMALIZED JENSEN POLYNOMIALS.\label{WMPJen}}%
To investigate directly a location of roots of the Minkowski
polynomials \(M_{B^n}^{\mathbb{R}^n}\),
\(M_{B^{n}\times{}0}^{\mathbb{R}^{n+1}}\),
\(M_{Q^n}^{\mathbb{R}^n}\),
\(M_{Q^{n}\times{}0}^{\mathbb{R}^{n+1}}\) \,and Weyl polynomials
\(W_{\partial{}B^{n+1}}^{\,p}\),
\(W_{\partial{}(B^{n}\times{}0)}^{\,p}\),
\(W_{\partial{}Q^{n+1}}^{\,\infty}\),
\(W_{\partial{}(Q^{n}\times{}0)}^{\,p}\) for a \textit{finite}
\(n\) is difficult. It turns out that it is much easier to
investigate first a location of roots of the entire functions
which are the limits of the (renormalized) Minkowski and Weyl
polynomials as \(n\to\infty\), and then to deduce properties of
roots of the original Minkowski  and Weyl polynomials from
properties of these limiting entire functions.
\paragraph{Jensen polynomials.}

From the explicit expressions \eqref{EMPb1}, \eqref{MiSqB},
\eqref{MPUQ}, \eqref{MPUQC} for the Minkowski polynomials and
\eqref{WPpS}, \eqref{WPiS},  \eqref{WPpSSqC}, \eqref{WPiSSqC},
\eqref{WPpQ}, \eqref{WPiQ}, \eqref{WPpQSqC}, \eqref{WPiQSqC} for
the Weyl polynomials we notice that each of this expressions
contains the factor \(\dfrac{n!}{(n-k)!}\), which is `a part' of
the binomial coefficient \(\binom{n}{k}\). The factorial ratio
can be presented as%
\begin{equation}%
\label{BinRat}%
 \frac{n!}{(n-k)!}=1\cdot{}
\Big(1-\frac{1}{n}\Big)\cdot{}\Big(1-\frac{2}{n}\Big)\cdot\,\,\cdots\,\cdot\,{}
\Big(1-\frac{k-1}{n}\Big)\cdot{}n^k\,,\quad{}\,1\leq{}k\leq{}n\,.
 \end{equation}
 \newpage
\begin{definition} \ \ %
\label{DefJP}%
\vspace*{-3.0ex}
\begin{enumerate}
\item
Given a formal power series \(f\):
\begin{equation}%
\label{FPS}%
 f(t)=\sum\limits_{0\leq{}l<\infty}a_{l}t^{l}\,.
\end{equation}
We associate with \(f\) the sequence of the polynomials
\(\mathscr{J}_n(f;t),\,n=1,\,2,\,3,\,\ldots\,\,:\)
\begin{equation}
\label{DJP}%
\mathscr{J}_n(f;t)=\sum\limits_{0\leq{}l\leq{}n}\frac{n!}{(n-l)!}\frac{1}{n^l}\cdot{}a_{l}t^{l},
\end{equation}
or, decoding the factor \(\frac{n!}{(n-l)!}\frac{1}{n^l}\),
\begin{equation}
\label{DJPd}%
\mathscr{J}_n(f;t)=a_0+\sum\limits_{1\leq{}l\leq{}n}1\big(1-{\textstyle{}\frac{1}{n}}\big)
\big(1-{\textstyle{}\frac{2}{n}}\big)\,\cdots\,\big(1-{\textstyle{}\frac{l-1}{n}}\big)\,
\cdot{}a_{l}t^{l}.
\end{equation}
The polynomials \(\mathscr{J}_n(f;t)\) are said to be \textsf{the
Jensen polynomials associated with the power series \(f\).}
\item
Given a \emph{function} \(f\) holomorphic in the disc
\(\{t:\,|t|<R\}\), where \(R\leq{}\infty\), we associate the
sequence of the Jensen polynomials with the Taylor series
\eqref{FPS} of the function \(f\) according the rule \eqref{DJP}.
We denote these
 polynomials by \(\mathscr{J}_n(f;t)\) as well and call them \textsf{the Jensen
 polynomials associated with the function \(f\).}
 \item
 The factors
 \begin{multline}
 \label{JeFa}
  j_{n,0}=1,\quad{}j_{n,k}=1\big(1-{\textstyle{}\frac{1}{n}}\big),
\big(1-{\textstyle{}\frac{2}{n}}\big)\,%
\cdots\,\big(1-{\textstyle{}\frac{k-1}{n}}\big),\,\,1\leq
  k\leq n,\\
j_{n,k}=0,\,\,k>n\,,
 \end{multline}
 are said to be \textsf{the Jensen multipliers}.\\ %
\hspace*{2.0ex}Thus, the Jensen polynomials associated with \(f\)
of the form \eqref{FPS}
 can be written as:
 \begin{equation}
\label{DJP1}%
\mathscr{J}_n(f;t)=\sum\limits_{0\leq{}<\infty}j_{n,l}\cdot{}a_{l}t^{l}\,.
\end{equation}
\end{enumerate}
\end{definition}
Since
\begin{math}%
\label{ReLRel}%
 j_{n,k}\to{}1 {}\ \textup{as \(k\) is fixed}, %
 \ n\to\infty\,,
\end{math}
 the following result is evident:
 \newpage
\begin{lemma}[The approximation property of Jensen polynomials.]{}\ \\%
\label{CJPL} %
Given the power series \eqref{FPS}, then:
\begin{enumerate}
\item
The sequence
 of the Jensen polynomials \(\mathscr{J}_n(f;t)\) converge to the series \(f\)
 coefficients-wise\,;
\item
If moreover the radius of convergence of the power series
\eqref{FPS} is positive, say equal to \(R,\,0<R\leq{}\infty\),
then the sequence
 of the Jensen polynomials \(\mathscr{J}_n(f;t)\) converge to the function
 which is the sum of this power series
 locally uniformly in the disc \(\{t:\,|t|<R\}\).
\end{enumerate}
\end{lemma}%

The approximation property in not specific for the polynomials
constructed from the \textit{Jensen multipliers} \(j_{n,k}\). This
property holds for \textit{any} multipliers \(j_{n,k}\) which
satisfy the conditions
\( j_{n,k}\to{}1 \,\textup{as \(k\) is fixed},%
\,n\to\infty\,,\) and are uniformly bounded:
\(\sup\limits_{k,n}|j_{n,k}|<\infty\,.\) What is much more
specific, that for some \(f\), the polynomial
\(\mathscr{J}_n(f;t)\) constructed from the \textit{Jensen
multipliers} \(j_{n,k}\) preserve the property of \(f\) to possess
only real roots. In particular:
\begin{nonumtheorem}
[{[}Jensen{]}] Let \(f\) be a \emph{polynomial} such that all its
roots are real. Then for each \(n\), all roots of the Jensen
polynomial \(\mathscr{J}_n(f,\,t)\) are real as well.
\end{nonumtheorem}
This result is a special case of Schur composition theorem
\cite{Schu1}. Actually, Jensen, \cite{Jen}, obtained a more
general result in which formulation \(f\) can be not only a
polynomial with real roots, but also an entire function which
belongs to the so called \emph{Laguerre-Polya class of entire
functions}. We return to this generalization later, is Section
\ref{LPEF}. Now we focus our attention on representation of the
Minkowski and Weyl polynomials as Jensen polynomials of certain
entire functions.

The relation \eqref{ReLRel} as well as the expressions
\eqref{EMPb1},\,\eqref{MiSqB},\,\eqref{MPUQ},\,\eqref{MPUQC} for
the Minkowski polynomials suggest us how the Minkowski polynomials
should be renormalized so that the renormalized polynomials tend
to a non-trivial limit as \(n\to\infty\).
\newpage
\paragraph{Entire functions which generate the Minkowski polynomials
for balls, cubes, spherical and cubic cylinders.} Let us introduce
the infinite power series:
\begin{subequations}
\label{LiEFMPS}
\begin{align}
 \mathcal{M}_{B^{\infty}}(t)&=\sum\limits_{0\leq{}k<\infty}
\frac{1}{k!}\,t^k\,;
\label{LiEFMPS1}\\%
\mathcal{M}_{B^{\infty}\times{}0}(t)&=\sum\limits_{0\leq{}k<\infty}
\frac{\Gamma(\frac{1}{2}+1)\Gamma(\frac{k}{2}+1)}{\Gamma(\frac{k+1}{2}+1)}\,
\frac{1}{k!}\,t^k\,;
\label{LiEFMPS2}\\%
\mathcal{M}_{Q^{\infty}}(t)&=\sum\limits_{0\leq{}k<\infty}
\frac{1}{\Gamma(\frac{k}{2}+1)k!}\Big(\frac{\sqrt{\pi}}{2}\Big)^k\,\,t^k\,;
\label{LiEFMPS3}\\%
\mathcal{M}_{Q^{\infty}\times{}0}(t)&=\sum\limits_{0\leq{}k<\infty}
\frac{\Gamma(\frac{1}{2}+1)}{\Gamma(\frac{k+1}{2}+1)k!}\Big(\frac{\sqrt{\pi}}{2}\Big)^k\,t^k\,.%
\label{LiEFMPS4}
\end{align}
\end{subequations}
The series \eqref{LiEFMPS} represent entire functions which grow
 not faster than  exponentially. More precisely, the
functions \(\mathcal{M}_{B^{\infty}}\) and
\(\mathcal{M}_{B^{\infty}\times{}0}\) grow exponentially: they are
of  order  \(1\) and normal type, the functions
\(\mathcal{M}_{Q^{\infty}}\) and
\(\mathcal{M}_{BQ^{\infty}\times{}0}\) grow subexponentially: they
are of order \(2/3\) and normal type.

With each of the entire functions \eqref{LiEFMPS} we associate the
sequence of polynomials which are the Jensen polynomials
associated with this entire function:
\begin{subequations}
\label{JLiEFMPS}
\begin{align}
 \mathcal{M}_{B^{n}}(t)&{=}
 \mathscr{J}_n(\mathcal{M}_{B^{\infty}};t)\,,
\label{JLiEFMPS1}\\%
\mathcal{M}_{B^{n}\times{}0}(t)&{=}\mathscr{J}_n(\mathcal{M}_{B^{\infty}\times{}0};t)\,,
\label{JLiEFMPS2}\\%
\mathcal{M}_{Q^{n}}(t)&{=}\mathscr{J}_n(\mathcal{M}_{Q^{\infty}};t)\,
\label{JLiEFMPS3}\\%
\mathcal{M}_{Q^{n}\times{}0}(t)&{=}
\mathscr{J}_n(\mathcal{M}_{Q^{\infty}\times{}0};t)\,,%
\label{JLiEFMPS4}
\end{align}
\end{subequations}
From the expressions \eqref{EMPb1}, \eqref{MiSqB}, \eqref{MPUQ},
\eqref{MPUQC} for the Minkowski polynomials it follows that they
are related to the above introduced polynomials \eqref{JLiEFMPS}
as:
\begin{subequations}
\label{RemOr}
\begin{alignat}{2}
M_{B^n}^{\mathbb{R}^n}(t)&=
\textup{Vol}_n(B^n)&\,&\mathcal{M}_{B^n}(nt)\,;
\label{RemOr1}\\%
M_{{B^n}\times{}0}^{\mathbb{R}^{n+1}}(t)&=
\textup{Vol}_n(B^{n})\omega_{1}t&\,&\mathcal{M}_{B^n\times{}0}(nt)\,;
\label{RemOr2}\\%
M_{Q^n}^{\mathbb{R}^n}(t)&=
\textup{Vol}_n(Q^n)&\,&\mathcal{M}_{Q^n}(nt)\,;
\label{RemOr3}\\%
M_{{Q^n}\times{}0}^{\mathbb{R}^{n+1}}(t)&=
\textup{Vol}_n(Q^n)\omega_{1}t&\,&\mathcal{M}_{Q^n\times{}0}(nt)\,;
\label{RemOr4}
\end{alignat}
\end{subequations}
The polynomials \(\mathcal{M}_{B^n}\),
\(\mathcal{M}_{B^n\times{}0}\), \(\mathcal{M}_{Q^n}\),
\(\mathcal{M}_{Q^n\times{}0}\) can be interpreted as
\textit{renormalized Minkowski polynomials}  respectively. We take
the equalities \eqref{RemOr} as the \textsf{definition} of the
renormalized Minkowski polynomials \(\mathcal{M}_{B^n}\),
\(\mathcal{M}_{B^n\times{}0}\), \(\mathcal{M}_{Q^n}\),
\(\mathcal{M}_{Q^n\times{}0}\) in terms of the `original'
Minkowski polynomials \(M_{B^n}^{\mathbb{R}^n}\),
\(M_{{B^n}\times{}0}^{\mathbb{R}^{n+1}}\),
\(M_{Q^n}^{\mathbb{R}^n}\), and
\(M_{{Q^n}\times{}0}^{\mathbb{R}^{n+1}}\),.

From the approximative property of Jensen polynomials and from
\eqref{JLiEFMPS} it follows that
\begin{multline}
\label{LiEFM} %
\mathcal{M}_{B^{n}}(t)\to\mathcal{M}_{B^{\infty}}(t),\ \ %
\mathcal{M}_{B^{n}\times{}0}(t)\to\mathcal{M}_{B^{\infty}\times{}0}(t),\ \ %
\mathcal{M}_{Q^{n}}(t)\to\mathcal{M}_{Q^{\infty}}(t),\\ %
\mathcal{M}_{Q^{n}\times{}0}(t)\to\mathcal{M}_{Q^{\infty}\times{}0}(t)\
\ \textup{as}\ \ n\to\infty\,.%
\end{multline}
This explains the notation \eqref{LiEFMPS}\,.

We summarize the above stated consideration as the following
\begin{theorem}%
\label{MPaJP}%
Let \(\{V^n\}\) be one of the four families of convex sets:
\(\{B^n\}\), \(\{B^n\times{}0\}\),  \(\{Q^n\}\),
\(\{Q^n\times{}0\}\). For each of these four families, there
exists the \textsf{single} entire function\,\footnote{The symbol
\(V^\infty\) means \(\{B^\infty\}\), \(\{B^\infty\times{}0\}\),
\(\{Q^\infty\}\) or
\(\{Q^\infty\times{}0\}\) respectively,} %
\(\mathcal{M}_{V^{\infty}}\)  such that in \textsf{every}
dimension \(n\), the renormalized Minkowski polynomials
\(\mathcal{M}_{V^{n}}\), defined  by  \eqref{RemOr}, are generated
by this entire function \(\mathcal{M}_{V^{\infty}}\) as the Jensen
polynomials \(\mathscr{J}_n(\mathcal{M}_{V^{\infty}})\): the
equalities \eqref{JLiEFMPS} hold.
\end{theorem}%
\vspace*{-6.0ex}
\paragraph{Entire functions which generate the Weyl polynomials
for the surfaces of balls, cubes, spherical and cubic cylinders.}
Let us introduce the infinite power series:
\begin{subequations}
\label{ReWP}
\begin{align}%
\mathcal{W}_{\partial{}B^{\infty}}^{\,p}(t)&=
\sum\limits_{l=0}^{\infty}%
\frac{2^{-l}\,\Gamma(\frac{p}{2}+1)}{\Gamma(\frac{p}{2}+l+1)}
\cdot{}\frac{1}{l!}\cdot{}\Big(-\frac{t^2}{2}\Big)^l,\ \
p=1,\,2,\,\ldots\,;
\label{ReWP1}\\
\mathcal{W}_{\partial{}B^{\infty}}^{\,\infty}(t)&=
\sum\limits_{l=0}^{\infty}%
\frac{1}{l!}\cdot{}\Big(-\frac{t^2}{2}\Big)^l\,;
 \label{ReWP2}\\
\mathcal{W}_{\partial{}(B^{\infty}\times{}0)}^{\,p}(t)&=
\sum\limits_{l=0}^{\infty}%
\frac{2^{-l}\,\Gamma(\frac{p}{2}+1)}{\Gamma(\frac{p}{2}+l+1)}
\cdot{}\frac{\Gamma(1/2)}{\Gamma(l+1/2)}\cdot{}\Big(-\frac{t^2}{2}\Big)^l,\
\  p=1,\,2,\,\ldots\,;%
 \label{ReWP3}\\ %
\mathcal{W}_{\partial{}(B^{\infty}\times{}0)}^{\,\infty}(t)&=
\sum\limits_{l=0}^{\infty}%
\frac{\Gamma(1/2)}{\Gamma(l+1/2)}\cdot{}\Big(-\frac{t^2}{2}\Big)^l\,; %
\label{ReWP4}  %
\end{align}
 \begin{align}
\mathcal{W}_{\partial{}Q^{\infty}}^{\,p}(t)&=
\sum\limits_{l=0}^{\infty}%
 \frac{2^{-l}\,\Gamma(\frac{p}{2}+1)}{\Gamma(\frac{p}{2}+l+1)}
\cdot{}\frac{1}{(2l+1)!}\cdot{}\Big(-{\frac{{\pi}t^2}{2}}\,\Big)^l,\
\  p=1,\,2,\,\ldots\,;%
\label{ReWP5}\\ %
\mathcal{W}_{\partial{}Q^{\infty}}^{\,\infty}(t)&=
\sum\limits_{l=0}^{\infty}%
\frac{1}{(2l+1)!}\cdot{}\Big(-{\frac{{\pi}t^2}{2}}\,\Big)^l\, ;
\label{ReWP6}\\ %
\mathcal{W}_{\partial{}(Q^{\infty}\times{}0)}^{\,p}(t)&=
\sum\limits_{l=0}^{\infty}%
\frac{2^{-l}\,\Gamma(\frac{p}{2}+1)}{\Gamma(\frac{p}{2}+l+1)}
\cdot{}\frac{1}{(2l)!}\cdot{}\Big(-{\frac{{\pi}t^2}{2}}\,\Big)^l,\
\ p=1,\,2,\,\ldots\,;%
\label{ReWP7}\\ %
\mathcal{W}_{\partial{}(Q^{\infty}\times{}0)}^{\,\infty}(t)&=
\sum\limits_{l=0}^{\infty}%
\frac{1}{(2l)!}\cdot{}\Big(-{\frac{{\pi}t^2}{2}}\,\Big)^l\,.
\label{ReWP8}
\end{align}%
\end{subequations}
The series \eqref{ReWP} represent entire functions. The functions
\eqref{ReWP2} and \eqref{ReWP4} are of order \(2\) and normal
type, the functions \eqref{ReWP1}, \eqref{ReWP3}, \eqref{ReWP6}
and \eqref{ReWP8} are of order \(1\) and normal type, the
functions \eqref{ReWP5} and \eqref{ReWP7} are of order \(2/3\) and
normal type.

With each of the entire functions \eqref{ReWP} we associate the
sequence of polynomials which are the Jensen polynomials
associated with this entire function:
\begin{subequations}
\label{JReWP}
\begin{alignat}{2}
 \mathcal{W}^p_{\partial{}B^{n+1}}(t)&{=}
 \mathscr{J}_{2[n/2]}(\mathcal{W}^p_{\partial{}B^{\infty}};t)\,,&\qquad&1\leq{}p\leq\infty\,;
\label{JReWP1}\\%
\mathcal{W}^p_{\partial{}(B^{n+1}\times{}0)}(t)&{=}
\mathscr{J}_{2[n/2]}(\mathcal{W}^p_{\partial{}(B^{\infty}\times{}0)};t)\,,&\qquad&1\leq{}p\leq\infty\,;
\label{JReWP2}\\%
\mathcal{W}^p_{\partial{}Q^{n}}(t)&{=}
\mathscr{J}_{2[n/2]}(\mathcal{W}^p_{\partial{}Q^{\infty}};t)\,,&\qquad&1\leq{}p\leq\infty\,;
\label{JReWP3}\\%
\mathcal{W}^p_{\partial{}(Q^{n}\times{}0)}(t)&{=}
\mathscr{J}_{2[n/2]}(\mathcal{W}^p_{\partial{}(Q^{\infty}\times{}0)};t)\,,&\qquad&1\leq{}p\leq\infty\,.
\label{JReWP4}
\end{alignat}
\end{subequations}
From the expressions \eqref{WPpS}, \eqref{WPiS}, \eqref{WPpSSqC},
\eqref{WPiSSqC},  \eqref{WPpQ}, \eqref{WPiQ}, \eqref{WPpQSqC},
\eqref{WPiQSqC},    for the Weyl polynomials it follows that they
are related to the above introduced polynomials \eqref{JReWP} as:
\begin{subequations}
\label{WRemOr}
\begin{alignat}{2}
W^p_{\partial{}B^{n+1}}(t)&=\phantom{(\times{})}
\textup{Vol}_n(\partial{}B^{n+1})&\,\cdot\,&\mathcal{W}^p_{\partial{}B^{n+1}}(int)\,;
\label{WRemOr1}\\%
W^p_{\partial{}(B^{n}\times{}0)}(t)&=
\textup{Vol}_n(\partial{}(B^{n}\times{}0))
\,&\,\cdot\,&\mathcal{W}^p_{\partial{}(B^{n}\times{}0)}(int)\,;
\label{WRemOr2}\\%
W^p_{\partial{}Q^{n+1}}(t)&=\phantom{(\times{})}
\textup{Vol}_n(\partial{}Q^{n+1})&\,\cdot\,&\mathcal{W}^p_{\partial{}Q^{n+1}}(int)\,;
\label{WRemOr3}\\%
W^p_{\partial{}(Q^{n}\times{}0)}(t)&=
\textup{Vol}_n(\partial{}(Q^{n}\times{}0))
\,&\,\cdot\,&\mathcal{W}^p_{\partial{}(Q^{n}\times{}0)}(int)\,;
\label{WRemOr4}
\end{alignat}
\end{subequations}
The equalities  \eqref{WRemOr} hold for all
\(n:\,1\leq{}n<\infty,\,\,p:\,1\leq{}p\leq{}\infty\).

The polynomials \(\mathcal{W}^p_{\partial{}B^{n+1}}\),
\(\mathcal{W}^p_{\partial(B^n\times{}0)}\),
\(\mathcal{W}^p_{\partial{}Q^{n+1}}\),
\(\mathcal{W}^p_{\partial{}(Q^n\times{}0)}\) can be interpreted as
\textit{renormalized Weyl polynomials}. We take the equalities
\eqref{WRemOr} as the \textsf{definition} of the renormalized Weyl
polynomials \(\mathcal{W}^p_{\partial{}B^{n+1}}\),
\(\mathcal{W}^p_{\partial(B^n\times{}0)}\),
\(\mathcal{W}^p_{\partial{}Q^{n+1}}\),
\(\mathcal{W}^p_{\partial{}(Q^n\times{}0)}\) in terms of the
`original' Minkowski polynomials \(W^p_{\partial{}B^{n+1}}\),
\(W^p_{\partial(B^n\times{}0)}\), \(W^p_{\partial{}Q^{n+1}}\),
\(W^p_{\partial{}(Q^n\times{}0)}\).

From the approximative property of Jensen polynomials and from
\eqref{JReWP} it follows that for every fixed
\(p,\,\,1\leq{}p\leq{}\infty\),
\begin{multline}
\label{LiEFMW} %
\mathcal{W}^p_{\partial{}B^{n+1}}(t)\to\mathcal{W}^p_{\partial{}B^{\infty}}(t),\ \ %
\mathcal{W}^p_{\partial{}(B^{n}\times{}0)}(t)
\to\mathcal{W}^p_{\partial{}(B^{\infty}\times{}0)}(t),\ \ \\ %
\mathcal{W}^p_{\partial{}Q^{n+1}}(t)\to\mathcal{W}^p_{\partial{}Q^{\infty}}(t), %
\mathcal{W}^p_{\partial{}(Q^{n}\times{}0)}(t)
\to\mathcal{W}^p_{\partial{}(Q^{\infty}\times{}0)}(t)\
\ \textup{as}\ \ n\to\infty\,.%
\end{multline}
This explains the notation \eqref{ReWP}\,.

We summarize the above stated consideration as the following
\begin{theorem}%
\label{WPaJP}%
Let \(\{\mathscr{M}^n\}\) be one of the four families of
\(n\)-dimensional convex surfaces: \(\{\partial{}B^{n+1}\}\),
\(\{\partial(B^n\times{}0)\}\), \(\{\partial{}Q^{n+1}\}\),
\(\{\partial(Q^n\times{}0)\}\). For each of these four families,
and for each \(p,\,1\leq{}p\leq{}\infty\), there exists the
\textsf{single} entire function\,\footnote{The symbol
\(\mathscr{M}^\infty\) means \(\{B^\infty\}\),
\(\{B^\infty\times{}0\}\), \(\{Q^\infty\}\) or
\(\{Q^\infty\times{}0\}\) respectively,} %
\(\mathcal{W}^p_{\,\mathscr{M}^{\infty}}\)  such that in
\textsf{every} dimension \(n\), the renormalized Weyl
polynomials \(\mathcal{W}^p_{\mathscr{M}^{n}}\), defined  by
\eqref{WRemOr}, are generated by this entire function
\(\mathcal{W}^p_{\mathscr{M}^{\infty}}\) as the Jensen polynomials
\(\mathscr{J}_2[n/2](\mathcal{W}^p_{\mathscr{M}^{\infty}})\).
\end{theorem}%

\section{ENTIRE FUNCTIONS OF THE HURWITZ\\
 AND OF THE LAGUERRE-POLYA CLASS.\\
MULTIPLIERS PRESERVING LOCATION OF ROOTS.
\label{LPEF}}%
\paragraph{Hurwitz class of entire functions.}
\begin{definition}
\label{DeHuC}
An  entire function \(H\) is said to
be in the \textsf{Hurwitz class}, written
\(H\in\mathscr{H}\), if
\begin{enumerate}
\item
\(H\not\equiv{}0\), and  roots of \(H\)
have negative real part: if \(H(\zeta)=0\),
then \(\textup{Re}\,\zeta<0\).
\item
The function \(H\) is of exponential type:
 \(\varlimsup\limits_{|z|\to\infty}\frac{\ln{}|H(z)|}{|z|}<\infty\),
 and its defect \(d_H\) is non-negative: \(d_H\geq{}0\), where
 \begin{equation}%
 \label{DefH}%
 2d_H=\textstyle{\varlimsup\limits_{r\to+\infty}\frac{\ln{}|H(r)|}{r}\,-
 \varlimsup\limits_{r\to+\infty}\frac{\ln{}|H(-r)|}{r}}\,.
 \end{equation}%
\end{enumerate}
\end{definition}

The following functions serve as examples of entire functions
of class \(\mathscr{H}\):\\
a). A dissipative polynomial \(P(t)\).\\ %
b). An exponential \(\exp\{\alpha{}t\}\), where
\(\textup{Re}\,\alpha\geq{}0\,\). \\ %
c). The product \(P(t)\cdot\exp\{\alpha{}t\}\): \(P(t)\) is a
dissipative polynomial, \(\textup{Re}\,\alpha\geq{}0\,\).

\textsf{ The significance of the Hurwitz class of entire functions stems from
the fact that function in this class\,%
\footnote{The full description of the class of entire functions which are the
 limits of dissipative polynomials can be found in \cite{Lev1},
 Chapter VIII, Theorem 4. This class (up to the change of variables \(z\to{}iz\))
 is denoted by \(P^{\ast}\) there.} %
are
the locally uniform limits in \(\mathbb{C}\) of
dissipative polynomials.}
\paragraph{Laguerre-Polya class of entire functions.}
\begin{definition}
\label{DFLPC} An  entire function \(E\) is said to
be in the \textsf{Laguerre-P{\'o}lya class}, written
\(E\in\mathscr{L}\text{-}\mathscr{P}\), if \(E\) is real
and
 can be expressed in the form
\begin{equation}
\label{DLP2}%
E(t)=ct^ne^{-{\beta}t^2+{\alpha}t}\prod\limits_{k=1}^{\infty}
\left(1+t\alpha_k\right)e^{-t\alpha_k},
\end{equation}
where
\(c\in\mathbb{R}\setminus{}0,\,\beta\geq{}0,\,\alpha\in\mathbb{R},\,
\alpha_k\in\mathbb{R}\), \(n\) is non-negative integer, and
\(\sum_{k=1}^{\infty}\alpha_k^{2}<\infty\).

Within the Laguerre-Polya class, those functions \(E\)
are said to be of \mbox{type \textup{I}}, written
\(E\in\mathscr{L}\text{-}\mathscr{P}\text{-}\textup{I}\)\,,
which are representable in the form
\begin{equation}
\label{DLP23}%
E(t)=ct^ne^{{\alpha}t}\prod\limits_{k=1}^{\infty}
\left(1+t\alpha_k\right),
\end{equation}
where \(c\in\mathbb{R}\setminus{}0,\,\,\alpha\geq{}0,\,
\alpha_k\geq{}0\), \(n\) is non-negative integer, and
\(\sum_{k=1}^{\infty}\alpha_k<\infty\).
\end{definition}
\hspace{3.0ex}%
\begin{minipage}[h]{0.9\linewidth}
\textsf{ The significance of the Laguerre-Polya class stems from
the fact that function in this class, \textit{and only these}, are
the locally uniform limits in \(\mathbb{C}\) of polynomials with
only real roots.} (See \cite{Lev1}, Chapter 8; \cite{Obr}, Chapter
II, Theorems 9.1,\,9.2,\,9.3.)
\end{minipage}

\begin{lemma}
\label{HLPI} An entire function \(E\) which is of type \textup{I}
in the Laguerre-Polya class also is the Hurwitz class:
\[\mathscr{L}\text{-}\mathscr{P}\text{-}\textup{I}\,\subset\,\mathscr{H}\,.\]
\end{lemma}
\textsf{PROOF.} The roots of the entire function \(E\) which admit
the representation \eqref{DLP23} are located at the points \(-(\alpha_k)^{-1}\),
thus is strictly negative. From the properties of the infinite product
\(\prod\limits_{k=1}^{\infty}\left(1+t\alpha_k\right)\)  %
with
 \(\sum_{k=1}^{\infty}|\alpha_k|<\infty\), it follows that
a function \(E\) which admit the representation \eqref{DLP23} is
of exponential type \(\alpha\), and
\(\varlimsup\limits_{r\to+\infty}\frac{\ln{}|H(\pm{}r)|}{r}=\pm{}\alpha\).
Thus, the defect \(d_H=\alpha\geq{}0\) since \(\alpha\geq{}0\).
 \hfill\framebox[0.45em]{ }

\paragraph{Multipliers preserving the reality of roots.}
\begin{definition}%
\label{PReZ}%
A sequence \(\{\gamma_k\}_{0\leq{}k<\infty}\) of real numbers is a
\textsf{multiplier sequence} if for every polynomial \(f\):
\[f(t)=\sum\limits_{0\leq{}k\leq{}n}a_kt^k\]
with only real roots, the polynomial
\[h(t)=\sum\limits_{0\leq{}k\leq{}n}\gamma_ka_kt^k\]
too has only real roots. (The degree \(n\) of the polynomial \(f\)
can be arbitrary.)
\end{definition}%
\begin{nonumtheorem} [{[Polya, Schur]}]
A sequence  \(\{\gamma_k\}_{0\leq{}k<\infty}\) of real numbers
which are not all roots is a multiplier sequence if and only if
the power series
\[\Psi(t)=\sum\limits_{0\leq{}k\leq{}\infty}\frac{\gamma_k}{k!}\,t^k\]
represents an entire function, and either the function \(\Psi(t)\)
or the function \(\Psi(-t)\) is in the Lagierre-Polya class of
type \textup{I}.
\end{nonumtheorem}
This result was obtained in \cite{PoSch}. The presentation of this
 and related results can be found in Chapter VIII of \cite{Lev1},
in Chapter II of \cite{Obr}, in \cite{RaSc} (Section 5 ), in
numerous papers by Th.\,Craven and G.\,Csordas.
\begin{theorem}{\textup{\textsf{[Jensen-Craven-Csordas-Williamson.]}}}
\label{JSW}
Let \(E(t)\) be an entire function belonging to the
Laguerre-Polya class \(\mathscr{L}\text{-}\mathscr{P}\), and
\(\{\mathscr{J}_n(E,\,t)\}_{n=1,\,2,\,3,\,\ldots}\) be
the sequence of the Jensen polynomials associated with the
function
\(\mathscr{E}\). \textup{(Definition \ref{DefJP}.)}
\begin{enumerate}
\item
Then for each \(n\), all roots of the polynomial
\(\mathscr{J}_n(E,\,t)\) are  real;\\
\item
If \(E(t)\) belongs to the subclass \(\mathscr{L}\text{-}\mathscr{P}\textup{-I}\)
of the Laguerre-Polya class \(\mathscr{L}\text{-}\mathscr{P}\), then
for each \(n\), all roots of the polynomial
\(\mathscr{J}_n(E,\,t)\) are  negative;\\
\item
If moreover \(E(t)\) is not of the form
\(E(t)=p(t)\,e^{\beta{}t}\), where \(p(t)\) is a polynomial, then
for each \(n\), all roots of the polynomial
\(\mathscr{J}_n(E,\,t)\) are simple\,.
\end{enumerate}
\end{theorem}
The statement 1 of the theorem was proved by
Jensen\,\footnote{Though Jensen himself did not introduce
explicitly the polynomials which are called `the Jensen
polynomials' now.}, \cite{Jen}. It is a special case of Theorem by
G.Polya and I.\,Schur corresponding to
\(\Psi(t)=\Big(1+\dfrac{t}{n}\Big)^n\). The refinement of the
statement 1 which is formulated as the statement 3 was done by
G.Csordas and J.\,Williamson in \cite{CsWi}, where the alternative
proof of the statement 1 also was done. In \cite{CsWi}, the main
Theorem formulated on p.\,263, which appeared as the statement 3
of Theorem \ref{JSW} of the present paper, was formulated not
accurately. The correction was done in \cite{CrCs3}, Section 4.1
there.
\begin{theorem}
\label{JDM}
Let \(H\) be an entire function belonging to the Hurwitz class \(\mathscr{H}\),
and \(\{\mathscr{J}_n(H,\,t)\}_{n=1,\,2,\,3,\,\ldots}\) be
the sequence of the Jensen polynomials associated with the
function
\(\mathscr{H}\). \textup{(Definition \ref{DefJP}.)}
Then for each \(n\), the polynomial \(\{\mathscr{J}_n(H,\,t)\)
is dissipative.
\end{theorem}
Theorem \ref{JDM} can be obtained as a consequence of Theorem
\ref{JSW} and Hermite-Bieler Theorem. Proof of Theorem \ref{JDM}
will be done in Section \ref{LPEF}.

\paragraph{Laguerre multipliers.}
\begin{nonumtheorem}[{[Laguerre]}]
Let an entire function \(E(t)\),
\begin{equation}
\label{EE}%
E(t)=\sum\limits_{0\leq{}l<\omega}\varepsilon_lt^l,\quad
\omega\leq{}\infty,
\end{equation}
be in the Laguerre-Polya class:
\(E\in\mathscr{L}\text{-}\mathscr{P}\), and let  an entire
function \(\psi\) be in the Laguerre-Polya class
\(\mathscr{L}\text{-}\mathscr{P}\) and moreover satisfy the
condition: all roots of \(\psi\) are negative.
\begin{enumerate}
\item
Then the power series
\begin{equation}
\label{PSPs}%
E_{\psi}=\sum\limits_{0\leq{}l<\omega}\varepsilon_l\psi(l)t^l
\end{equation}
converges for every \(t\), and its sum is an entire function of
the Laguerre-Polya class:
\(E_{\psi}\in\mathscr{L}\text{-}\mathscr{P}\).
\item
If moreover \(E(t)\) is of type \textup{I}:
\(E\in\mathscr{L}\text{-}\mathscr{P}\textup{-I}\), then the the
sum of power series \eqref{PSPs} also is an entire function of the
type \textup{I}:
\(E_{\psi}\in\mathscr{L}\text{-}\mathscr{P}\textup{-I}\).
\end{enumerate}
\end{nonumtheorem}
This theorem appeared by E.\,Laguerre, \cite{Lag1}, section 18,
p.117, or \cite{Lag}, p.\,202. Laguerre himself has formulated
this theorem for the function \(E\) which is a polynomial with
real roots. The extended formulation, where \(\E\) is a general
entire function from the class \(\mathscr{L}\text{-}\mathscr{P}\),
can be found in the paper \cite{PoSch}, p.\,112, or in its reprint
in \cite{Po}, p.123. In \cite{PoSch} the extended formulation is
attributed to Jensen, \cite{Jen}.

The presentation of the above mentioned results of Polya, Schur,
Laguerre, Jensen, as well as of many related results, can be found
in \cite{Obr}, Chapter II; \cite{Lev1}, Chapter VIII,\cite{RaSc};
\cite{RaSc}, Chapter 5, especially Sections 5.5, 5.6,  5.7; in
numerous papers of Th.\,Craven and G.\,Csordas (See for example
\cite{CrCs1}). See also \cite{PoSz}, Part five. The book of L.\,de
Branges \cite{deBr} is closely related to this circle of problems.
\newpage

\section{PROPERTIES OF ENTIRE FUNCTIONS\\
GENERATING MINKOWSKI AND WEYL POLYNOMIALS\\
OF `REGULAR' CONVEX SETS AND THEIR SURFACES. \label{PEFG}}%
\paragraph{Entire functions generating the Minkowski polynomials.}
\begin{theorem}\ \ %
\label{EFGWP}%
The entire functions \eqref{LiEFMPS} generating the renormalised
Minkowski polynomials of balls, cubes, squeezed spherical and
cubic cylinders, possesses the following properties:
\begin{enumerate}
\item
The function \(\mathcal{M}_{B^{\infty}}\) is of type \textup{I} of
the Laguerre-Polya class;
\item
The function \(\mathcal{M}_{B^{\infty}\times{}0}\) belongs to the
Hurwitz class. It has infinitely roots, all but finitely many its
roots are non-real;
\item
The function \(\mathcal{M}_{Q^{\infty}}\) is of type \textup{I} of
the Laguerre-Polya class;
\item
The function \(\mathcal{M}_{Q^{\infty}\times{}0}\) is of type
\textup{I} of the Laguerre-Polya class\,.
\end{enumerate}
\end{theorem}%
\begin{lemma}
\label{PrGa} The function \(\dfrac{1}{\Gamma(t+1)}\), where
\(\Gamma\) is the Euler Gamma function, is in Laguerre-Polya
class, and all its roots is negative.
\end{lemma}
Indeed,
\[\frac{1}{\Gamma(t+1)}=
e^{Ct}\prod\limits_{1\leq{}k<\infty}\left(1+\frac{t}{k}\right)e^{-\frac{t}{k}}\,,\]
(\(C\) is the Euler constant,
\(C\approx{}0.5772156\ldots\,\,\,.\))\hfill\framebox[0.45em]{ } %

\textsf{PROOF of Theorem \ref{EFGWP}}. The statement 1 is evident:
\(\mathcal{M}_{B^{\infty}}(t)=e^t\). \\%
\hspace*{3.0ex}To obtain Statement 3, we remark that the function
\(\mathcal{M}_{Q^{\infty}}\) is of the form \(E_{\psi}\),
\eqref{PSPs}, where \(E(t)=\exp\{\frac{\sqrt{\pi}}{2}t\}\), and
\(\psi(t)=\dfrac{1}{\Gamma(\frac{t}{2}+1)}\). Then we apply the
Laguerre theorem on multipliers to these \(E\) and \(\psi\). The
needed property of \(\psi\) is formulated as Lemma \ref{PrGa}.\\ %
\hspace*{3.0ex}The statement 4 can be obtained in the same way
that the statement 3. One need take
\(E(t)=\exp\{\frac{\sqrt{\pi}}{2}t\}\), and
\(\psi(t)=\dfrac{\Gamma(\frac{1}{2})}{\Gamma(\frac{t+1}{2}+1)}\).\\ %
\hspace*{3.0ex}Proof of the statement 2 is more complicated. From
\eqref{LiEFMPS2} it follows that
\[\mathcal{M}_{B^{\infty}\times{}0}(t)=
\sum\limits_{0\leq{}k<\infty}B({\textstyle{\frac{k}{2}+1,\frac{1}{2}}})\dfrac{1}{k!}\,t^k
=\sum\limits_{0\leq{}k<\infty}
\int\limits_0^1{\xi}^{\frac{k}{2}}(1-\xi)^{-\frac{1}{2}}\,d\xi\,%
\dfrac{1}{k!}\,t^k\,.\] Changing the order of summation and
integration and summarizing the exponential series, we obtain the
integral representation:
\begin{equation}
\label{IntRep}%
\mathcal{M}_{B^{\infty}\times{}0}(t)=
 2\int\limits_0^1{}(1-\xi^2)^{-\frac{1}{2}}\xi{}e^{\xi{}t}d\xi\,.
\end{equation}
The fact that the functions \(\mathcal{M}_{B^{\infty}\times{}0}\)
belongs to the Hurwitz class will be derived from the integral
representation \eqref{IntRep}. This will be done in %
Section~\ref{LocRoot}.\hfill\framebox[0.45em]{ }

\paragraph{Entire functions generating the Weyl polynomials.}
\begin{lemma}%
\label{PLaM}%
 Let \(E\):
\begin{equation}
\label{EEF}%
E(t)=\sum\limits_{0\leq{}l<\infty}a_l{}t^{2l}
\end{equation}
 be an even entire function of the class
\(\mathscr{L}\text{-}\mathscr{P}\), and let \(p>0\) be a number.
Then the function \(E_p(t)\) defined by the power series
\begin{equation}
\label{EEFs}%
E_p(t)\stackrel{\textup{\tiny
def}}{=}\sum\limits_{1\leq{}l<\infty}\frac{2^{-l}\Gamma\big(\frac{p}{2}+1\big)}
{\Gamma\big(l+\frac{p}{2}+1\big)}\cdot{}a_lt^{2l},
\end{equation}
belongs to the class \(\mathscr{L}\text{-}\mathscr{P}\) as well.
\end{lemma}
\textsf{PROOF.} %
Lemma \ref{PLaM} is the consequence of the Laguerre theorem on
multipliers. The function
\begin{equation}%
\label{LaMup}
\psi_p(t)=\frac{2^{-\frac{t}{2}}\Gamma\big(\frac{p}{2}+1\big)}
{\Gamma\big(\frac{t}{2}+\frac{p}{2}+1\big)}
\end{equation}%
is in the Laguerre-Polya class (see Lemma \ref{PrGa}), and its
roots are negative.\hfill\framebox[0.45em]{ }

\vspace{2.0ex}%
We point out, see \eqref{ReWP}, that the entire functions
\(\mathcal{W}_{\partial{}B^{\infty}}^{\,p}\),
\(\mathcal{W}_{\partial{}(B^{\infty}\times{}0)}^{\,p}\),
\(\mathcal{W}_{\partial{}Q^{\infty}}^{\,p}\),
\(\mathcal{W}_{\partial{}(Q^{\infty}\times{}0)}^{\,p}\), which
generate the Weyl polynomials of the finite index \(p\) for the
appropriate families of convex surfaces, can be obtained from the
entire functions
\(\mathcal{W}_{\partial{}B^{\infty}}^{\,\infty}\),
\(\mathcal{W}_{\partial{}(B^{\infty}\times{}0)}^{\,\infty}\),
\(\mathcal{W}_{\partial{}Q^{\infty}}^{\,\infty}\),
\(\mathcal{W}_{\partial{}(Q^{\infty}\times{}0)}^{\,\infty}\),
which generate the Weyl polynomials of the infinite index, by
means of the transformation of the form
\[\sum\limits_{0\leq{}k<\infty}a_kt^k\to\,
\sum\limits_{0\leq{}k<\infty}\psi_p(k)\,a_kt^k\,.\]
\begin{theorem}
\label{BLPC}%
\ \ %
\begin{enumerate}
\item %
 The functions \(\mathcal{W}_{\partial{}B^{\infty}}^{\,\infty}\),
\(\mathcal{W}_{\partial{}Q^{\infty}}^{\,\infty}\),
\(\mathcal{W}_{\partial{}(Q^{\infty}\times{}0)}^{\,\infty}\)
belong to the Laguerre-Polya class
\(\mathscr{L}\text{-}\mathscr{P}\).
\item
The function
\(\mathcal{W}_{\partial{}(B^{\infty}\times{}0)}^{\,\infty}\) does
not belong to the Laguerre-Polya class
\(\mathscr{L}\text{-}\mathscr{P}\): this function has infinitely
many non-real roots.
\end{enumerate}
\end{theorem}
\noindent%
\textsf{PROOF}. The statement 1 is evident in view of the explicit
expressions:
\begin{align}%
\label{EE1}%
\mathcal{W}_{\partial{}B^{\infty}}^{\,\infty}(t)&=\exp\{-t^2/2\}\,,\\
\label{EE2}%
 \mathcal{W}_{\partial{}Q^{\infty}}^{\,\infty}\,(t)&=
\frac{\sin\{(\pi/2)^{\frac{1}{2}}{}t\}}{(\pi/2)^{\frac{1}{2}}{}t}\,,\\
\label{EE3}%
 \mathcal{W}_{\partial{}(Q^{\infty}
\times{}0)}^{\,\infty}&=\cos\{(\pi/2)^{\frac{1}{2}}{}\,t\}\,.
\end{align}%
The function \(\mathcal{W}_{\partial{}(B^{\infty}
\times{}0)}^{\,\infty}\), which appears in \textup{Statement 2},
can not be expressed in terms of `elementary' functions, but it
can be expressed in terms of the Mittag-Leffler function
\(\mathscr{E}_{1,\,\frac{1}{2}}\):
\begin{equation}
\label{EE4}%
 \mathcal{W}_{\partial{}(B^{\infty}
\times{}0)}^{\,\infty}(t)=\sqrt{\pi}\mathscr{E}_{1,\,\frac{1}{2}}\left(-\frac{t^2}{2}\right),
\end{equation}
where
\begin{equation}%
\label{GML}%
\mathscr{E}_{\alpha,\beta}(z)=\sum\limits_{0\leq{}k<\infty}\frac{z^k}
{\Gamma(\alpha{}k+\beta)}\,.
\end{equation}
From \eqref{GML} the integral representation can be derived:
\begin{equation}%
\label{GMLI}%
\sqrt{\pi}\mathscr{E}_{1,\,\frac{1}{2}}(t)=1+
t\int\limits_0^1(1-\xi)^{-\frac{1}{2}}e^{t\xi}\,d\xi\,.
\end{equation}
The integral representation \eqref{GMLI} can be derived from the
Tailor series \eqref{GML} in the same way as the integral
representation \eqref{IntRep} was derived from the Taylor series
\eqref{LiEFMPS2}. From \eqref{GMLI} the following asymptotic can
be obtained:
\begin{equation}%
\label{AsML}%
\sqrt{\pi}\mathscr{E}_{1,\,\frac{1}{2}}(t)=
\begin{cases}%
\hspace*{2.0ex}\frac{1}{2t},& t\to{}-\infty,\\[1.0ex]
\sqrt{\pi{}t}\, e^t, &t\to{}+\infty.\,\\[1.0ex]
O(|t|),&t\to{}\pm{}i\infty.
\end{cases}%
\end{equation}%
From \eqref{AsML} it follows that the indicator diagram of the
entire function \(\mathscr{E}_{1,\,\frac{1}{2}}(t)\) of the
exponential type is the interval \([0,\,1]\,.\) Moreover, the
function \(\mathscr{E}_{1,\,\frac{1}{2}}(it)\) belongs to the
class \(C\), as this class was defined in \cite{Lev2}, Lecture 17.
From  Theorem of Cartwright-Levinson (Theorem 1 of the Lecture 17
from \cite{Lev2}) it follows that the function
\(\mathscr{E}_{1,\,\frac{1}{2}}(t)\) has infinitely many roots,
these roots have a positive density, and are located `near' the
rays \(\arg t=\frac{\pi}{2}\) and \(\arg t=-\frac{\pi}{2}\). From
this and from \eqref{EE4} it follows that the roots of the
function \( \mathcal{W}_{\partial{}(B^{\infty}
\times{}0)}^{\,\infty}(t)\) are located near four rays
\(\arg{}t=\frac{\pi}{4},\,\arg{}t=
\frac{3\pi}{4},\,\arg{}t=\frac{5\pi}{4},\,\arg{}t=\frac{7\pi}{4}\,.\)
In particular, infinitely many of the roots of the function \(
\mathcal{W}_{\partial{}(B^{\infty} \times{}0)}^{\,\infty}(t)\) are
non-real.
 \hfill\framebox[0.45em]{ }
 \begin{remark}
Much more precise results about the Mittag-Leffler function
\(\mathscr{E}_{\alpha,\beta}\) and distribution of its roots are
known. See, for example, \textup{\cite{EMOT},\,\,sec\-tion 18.1},
or \textup{\cite{Djr}}.
 \end{remark}

\begin{theorem}\ \\[-6.0ex] %
\label{ImCo}%
\begin{enumerate}
\item
For every \(p=1,\,2,\,\ldots\), the functions
\(\mathcal{W}_{\partial{}B^{\infty}}^{\,p}\),
\(\mathcal{W}_{\partial{}Q^{\infty}}^{\,p}\),
\(\mathcal{W}_{\partial{}(Q^{\infty}\times{}0)}^{\,p}\) belong to
the Laguerre-Polya class \(\mathscr{L}\text{-}\mathscr{P}\).
\item
If \(p\) is large enough, then the function
\(\mathcal{W}_{\partial{}(B^{\infty}\times{}0)}^{\,p}\) does not
belong the Laguerre-Polya class
\(\mathscr{L}\text{-}\mathscr{P}\): it has non-real roots.
\end{enumerate}
\end{theorem}
\noindent \textsf{PROOF}. The statement 1 of this theorem is a
consequence of the statement 1 of Theorem \ref{BLPC} and Lemma
\ref{PLaM}. The statement 2 of this theorem is a consequence of
the statement 2 of Theorem \ref{BLPC} and the approximational
property \eqref{LiRe}. \hfill\framebox[0.45em]{ }
\begin{remark}%
\label{RefToBes}%
The fact that the function
\(\mathcal{W}_{\partial{}B^{\infty}}^{\,p}\) belongs to the
Laguerre-Polya class \(\mathscr{L}\text{-}\mathscr{P}\), that is
all its roots are real, can be established without reference to
\textup{Lemma \ref{PLaM}}. The function
\(\mathcal{W}_{\partial{}B^{\infty}}^{\,p}\) can be expressed in
terms of Bessel functions \(J_{\nu}\). Recall that for arbitrary
\(\nu\),
\begin{equation}%
\label{Bess}%
J_{\nu}(t)=\left(\frac{t}{2}\right)^{\nu}\sum\limits_{0\leq{}l<\infty}%
\frac{(-1)^l(t^2/4)^l}{l!\,\Gamma(\nu+l+1)}\,.
\end{equation}%
Comparing \eqref{Bess} with \eqref{ReWP1}, we see that
\begin{equation}%
\label{CompBess}%
\mathcal{W}_{\partial{}B^{\infty}}^{\,p}(t)=\Gamma\left(\frac{p}{2}+1\right)
\left(\frac{t}{2}\right)^{-\frac{p}{2}}J_{\frac{p}{2}}(t)\,.
\end{equation}%
In particular, \\ %
 for\,%
\footnote{ Deriving \eqref{ExC1} from \eqref{CompBess}, we used
the formula
\(J_{\frac{1}{2}}(t)=\big(\frac{2}{\pi{}t}\big)^{\frac{1}{2}}\sin{}t\,.\)
\textup{(}Concerning this formula, see, for example,
\textup{\cite{WhWa}, section \textsf{17.24}.}\textup{)} However,
\eqref{ExC1} may be obtained directly from \eqref{ReWP1}.}
\(p=1\),
\begin{equation}%
\label{ExC1}%
\mathcal{W}^{\,1}_{\partial{}B^{\infty}}(t)=\frac{\sin{}t}{t},
\end{equation}
for
 \(p=2\),
\begin{equation}%
\label{ExC2}%
\mathcal{W}^{\,2}_{\partial{}B^{\infty}}(t)=2\,\frac{J_1(t)}{t}.
\end{equation}
It is known that for every \(\nu>-1\), all roots of the Bessel
function \(J_{\nu}(t)\) are real (This result is due to A.Hurwitz.
See, for example, \textup{\cite{Wat}, Chapter XV, Section 15.27}.)
\end{remark}%
The statement 2 of Theorem \ref{ImCo} may be strengthen
essentially.
\begin{theorem}\ \\[-6.0ex] %
\label{str}
\begin{enumerate}
\item
For \(p=1,\,2,\,4\), the function
\(\mathcal{W}_{\partial{}(B^{\infty}\times{}0)}^{\,p}\)
 belongs the Laguerre-Polya class
\(\mathscr{L}\text{-}\mathscr{P}\)\,;
\item
For \(p:\,5\leq{}p\leq{}\infty\), the function
\(\mathcal{W}_{\partial{}(B^{\infty}\times{}0)}^{\,p}\) does not
 belong the Laguerre-Polya class
\(\mathscr{L}\text{-}\mathscr{P}\)\,: it has infinitely many
non-real roots\,.
\end{enumerate}
\end{theorem}
\textsf{PROOF.} For every \(p\geq{}1\), the function
\(\mathcal{W}_{\partial{}(B^{\infty}\times{}0)}^{\,p}\) admits
 the integral representation
\begin{equation}
\label{IRWGF}
\mathcal{W}_{\partial{}(B^{\infty}\times{}0)}^{\,p}= %
p\int\limits_{0}^{1}(1-\xi^2)^{\frac{p}{2}-1}\xi\cos{}t\xi\,d\xi\,.
\end{equation}
This integral representation can be obtained from \eqref{ReWP3} in
the same way that the integral representation \eqref{IntRep} was
obtained from \eqref{LiEFMPS2}. Using the identity
\[\Gamma(l+1/2)\,\Gamma(l+1)=\Gamma(1/2)\,2^{-2l}\,\Gamma(2l+1),\]
we reduce  \eqref{ReWP3} to the form
\[\mathcal{W}_{\partial{}(B^{\infty}\times{}0)}^{\,p}=
\frac{p}{2}\sum\limits_{0\leq{}l<\infty}\text{B}(l+1,p/2)(-1)^l\frac{t^{2l}}{(2l)!}\,.\]
Then we use the integral representation for the function Beta,
change the order of summation and integration and summarize the
series using the Taylor expansion for \(\cos z\).
 For every \(p: \,1\leq{}p<\infty\), the function
\(\mathcal{W}_{\partial{}(B^{\infty}\times{}0)}^{\,p}\) can be
calculated asymptotically. This calculation may be done using the
integral representation \eqref{IRWGF}, or in other way. The
asymptotic expression for the function
\(\mathcal{W}_{\partial{}(B^{\infty}\times{}0)}^{\,p}\) is
presented in Section \ref{LocRoot}, see \eqref{AsyW},
\eqref{EstRema}. From this expression
it follows that:\\
1. For \(p>4\), infinitely many (actually all but finitely many)
roots of the
\(\mathcal{W}_{\partial{}(B^{\infty}\times{}0)}^{\,p}\) are
non-real. This is sufficiently for the negative result of the
statement 2 of Theorem \ref{str} to be obtained.
\\
2. For \(p\leq{}4\), all but finitely many roots of the function
\(\mathcal{W}_{\partial{}(B^{\infty}\times{}0)}^{\,p}\) are real
and simple. This alone is not sufficiently for  the result of the
statement 1 of to be obtained. The additional reasoning should be
invoked. For \(p=2\) and \(p=4\), the function
\(\mathcal{W}_{\partial{}(B^{\infty}\times{}0)}^{\,p}\) can be
calculated explicitly. The case \(p=3\) remains open. Proof of the
fact that for \(p=1,\,2,\,4\) all roots of the function
\(\mathcal{W}_{\partial{}(B^{\infty}\times{}0)}^{\,p}\) are real
will be done in Section \ref{LocRoot}. See Lemma \ref{WPBLP}.
\hfill\framebox[0.45em]{ }\\

\noindent%
\textsf{PROOF of THEOREM \ref{LoWP}.} According to Theorems
\ref{BLPC}, \ref{ImCo} and \ref{str} (Statements 1 of these
theorems), each of the functions
\(\mathcal{W}_{\partial{}B^{\infty}}^{\,p}\),
\(\mathcal{W}_{\partial{}Q^{\infty}}^{\,p}\),
\(\mathcal{W}_{\partial{}(Q^{\infty}\times{}0)}^{\,p}\) with
\(p:\,1\leq{}p\leq{}\infty\), and
\(\mathcal{W}_{\partial{}(B^{\infty}\times{}0)}^{\,p}\) with
\(p=1,\,2,\,4\), belongs to the class of Laguerre-Polya
\(\mathscr{L}\text{-}\mathscr{P}\). By Theorem of
Jensen-Csordas-Williamson, the Jensen polynomials associated with
each of these entire functions, has only simple real roots.
According to Theorem \ref{WPaJP}, the renormalized Weyl
polynomials \(\mathcal{W}^p_{\partial{}B^{n+1}}\),
\(\mathcal{W}^p_{\partial{}Q^{n+1}}\),
\(\mathcal{W}^p_{\partial{}(Q^{n}\times{}0)}\) with
\(p:\,1\leq{}p\leq{}\infty\), and
\(\mathcal{W}^p_{\partial{}(B^{n}\times{}0)}\) with
\(p=1,\,2,\,4\) have only simple real roots. In view of
renormalizing relations \eqref{WRemOr}, the Weyl polynomials
\(W^p_{\partial{}B^{n+1}}\), \(W^p_{\partial{}Q^{n+1}}\),
\(W^p_{\partial{}(Q^{n}\times{}0)}\) are conservative.
\hfill\framebox[0.45em]{ }

\noindent \textsf{PROOF of THEOREM \ref{NRWPSq}.} According to
Theorem \ref{str}, Statement 2, for \(p:\,5\leq{}p\leq{}\infty\),
each of the entire functions
\(\mathcal{W}_{\partial{}(B^{\infty}\times{}0)}^{\,p}\) has
infinitely many non-real roots. Since for fixed \(p\),
\(\mathscr{J}_n(\mathcal{W}_{\partial{}(B^{\infty}\times{}0)}^{\,p};t)
\to\mathcal{W}_{\partial{}(B^{\infty}\times{}0)}^{\,p}(t)\)
locally uniformly in \(\mathbb{C}\) as \(n\to\infty\), by Hurwitz
theorem, every polynomial
\(\mathscr{J}_n(\mathcal{W}_{\partial{}(B^{\infty}\times{}0)}^{\,p})\)
with \(p,\,n:\,\,p\geq{}5,\,\,n\geq{}N(p)\) has non-real roots. By
Theorem \ref{WPaJP}, the renormalized Weyl polynomials
\(\mathcal{W}^p_{\partial{}(B^{n}\times{}0)}\) have non-real
roots. In view of the renormalizing relations \eqref{JReWP}, the
Weyl polynomial \(W^p_{\partial{}(B^{n}\times{}0)}\) have roots
which do not belong to the imaginary axis.
\hfill\framebox[0.45em]{ } \vspace{2.0ex}

\noindent%
\textsf{PROOF of LEMMA \ref{LWPo}.} This lemma is a consequence of
lemma \ref{PLaM}. If the polynomial
\(W_{\mathscr{M}}^{\,\infty}(t)\) is conservative, then the
polynomial \(E(t)=W_{\mathscr{M}}^{\,\infty}(it)\) is a real
polynomial with only real simple roots. The function
\(E_p(t)=W_{\mathscr{M}}^{\,p}(it)\) is related with this \(E(t)\)
as well as the function \(E_p(t)\) appeared in \eqref{EEFs} is
related to \(E(t)\) from \eqref{PLaM}. By Lemma \ref{PLaM}, all
roots of \(E_p\) are real. Let us show that the roots are simple.
Consider the function \(E(t)+\varepsilon\), were \(\varepsilon\)
is a small real number, positive or negative. Since all roots of
the polynomial \(E(t)\) are real and simple, all roots of the
polynomial \(E(t)+\varepsilon\) are real if \(|\varepsilon|\) is
small enough. By Lemma \ref{PLaM}, all roots of the polynomial
\(E_p(t)+\varepsilon\) are real. But if the polynomial \(E_p(t)\)
have a multiple root,  this root splits into a group of simple
roots by the perturbation \(E_p(t)\to{}E_p(t)+\varepsilon\), and
by an appropriate choice of sign of \(\varepsilon\), some of roots
in this group will be non-real.\hfill\framebox[0.45em]{ }
\begin{remark}%
\label{SiSpC}%
We apply \textup{Lemma \ref{LWPo}} in special cases
\(n=2,\,3,\,4,\,5\) only. In these cases Lemma is quite
elementary. Actually only the cases \(n=4\) and \(n=5\) deserve
attention, the cases \(n=2\) and \(n=3\) are trivial.
The cases \(n=4\) and \(n=5\) are reduced to the following elementary statement:\\[1.0ex]
\hspace*{3.0ex}\begin{minipage}[t]{0.9\linewidth} \textsl{Let
\(k_0,\,k_2,\,k_4\) be positive numbers. Assume that the roots of
the polynomial \(Q(t)=k_0+k_2t+k_4t^2\) are negative and
different. Then for every \(p>0\), the roots of the polynomial
\[Q^p(t)=k_0+\frac{k_1}{(p+2)}t+\frac{k_4}{(p+2)(p+4)}t^2\]
are negative and different as well.}
\end{minipage}
\end{remark}%

Indeed, the conditions posed on polynomials \(Q\) and \(Q^p\) are
equivalent to the inequalities
\[k_1^2>k_0k_4\quad\textup{and}\quad{}\bigg(\frac{k_1}{p+2}\bigg)^2>k_0\frac{k_2}{(p+2)(p+4)}\,. \]
It is evident that the first of these inequalities implies the
second.\hfill\framebox[0.45em]{ }

\section{HERMITE-BIELER THEOREM AND ITS APPLICATION.
\label{HBieT}}%
In its traditional form, Hermite-Bieler theorem gives conditions under
which all roots of a polynomial belong to the upper half-plane
\(\{z:\,\textup{Im}\,z>0\}\). We need the version of this theorem
adopted to the left half-plane, and for the case of polynomials
with non-negative coefficients only. Before to present such a
reformulation of Hermite-Bieler theorem, we give several
definitions:
\begin{definition} %
\label{DefInlac}
Let \(S_1\) and \(S_2\) be two sets which are situated on the same straight
line\footnote{In our considerations the straight line \(L\) will be either the real axis
or the imaginary axis.}
 \(L\) of the complex plane:
\(S_1\subset{}L,\ S_2\subset{}L\,\), and moreover each of the sets
\(S_1,S_2\) consists of isolated points only. The sets \(S_1\) and
\(S_2\) \textsf{interlace} if between every two points of \(S_1\)
there is a point of \(S_2\), and between every two points of
\(S_2\) there is a point of \(S_1\).
\end{definition}
\begin{definition} %
\label{DefRIPa} %
Let \(P\) be a power series:
\begin{equation}
\label{CoPo} P(t)=\sum\limits_{0\leq{}k}p_kt^k,
\end{equation}
where \(t\) is a complex variable, and the coefficients \(p_k\)
are  complex numbers.

 The \textsf{real part} \ \({}^{\mathscr{R}\!}P\)
and the
\textsf{imaginary part} \ \({}^{\mathscr{I}\!}P\)  of \(P\) are defined as %
\begin{equation}%
\label{RIPa}
{}^{\mathscr{R}\!}P(t)=\frac{P(t)+\overline{P(\overline{t})}}{2},
\quad{}{}^{\mathscr{I}\!}P(t)=\frac{P(t)-\overline{P(\overline{t})}}{2i}\,,
\end{equation}%
where the overline bar is used as a notation for the complex
conjugation.

The  \textsf{even part} \({}^{\mathscr{E}\!}P\) and the
\textsf{odd part} \({}^{\mathscr{O}\!}P\)  of \(P\) are defined as %
\begin{equation}%
\label{EOIPa} {}^{\mathscr{E}\!}P(t)=\frac{P(t)+P(-t)}{2},
\quad{}{}^{\mathscr{O}\!}P(t)=\frac{P(t)-P(-t)}{2}\,,
\end{equation}%
In term of coefficients,
 \begin{subequations}
 \label{TeCo}
\begin{gather}%
\label{TeCo1}%
{}^{\mathscr{R}\!}P(t)=\sum\limits_{0\leq{}k}a_kt^k,\quad
{}^{\mathscr{I}\!}P(t)=\sum\limits_{0\leq{}k}b_kt^k, \\
\intertext{where}%
 a_k=\frac{p_k+\overline{p_k}}{2},
\quad{}b_k=\frac{p_k-\overline{p_k}}{2i}\,,
\label{TeCo2}%
\end{gather}%
\end{subequations}
and
\begin{equation}%
{}^{\mathscr{E}\!}P(t)=\sum\limits_{0\leq{}l}p_{2l}t^{2l},\quad
{}^{\mathscr{O}\!}P(t)=\sum\limits_{0\leq{}l}p_{2l+1}t^{2l+1}.
\label{EOCo}%
\end{equation}%
\end{definition}
\begin{nonumtheorem}[[Hermite-Bieler{]}] %
Let \(P\) be a polynomial, \(A={}^{\mathscr{R}\!}P\) and
\(B={}^{\mathscr{I}\!}P\) be the real and imaginary parts of
\(P\), i.e.
\[P(t)=A(t)+iB(t),\]
where \(A\) and \(B\) be a polynomials with real coefficients. In
order for all roots of \(P\) to be contained within the open upper
half-plane \(\{z:\,\textup{Im}\,z>0\}\), it is necessary and
sufficient that the following three conditions be satisfied:
\begin{enumerate}
\item
The roots of each of the polynomials \(A\) and \(B\) are all real
and simple.
\item
The sets \(\mathscr{Z}_A\) and \(\mathscr{Z}_B\) of the roots of
the polynomials \(A\) and \(B\) interlace.
\item
The inequality
\begin{equation}%
\label{ALN} %
B^{\prime}(0)A(0)-A^{\prime}(0)B(0)>0
\end{equation}%
 holds.
\end{enumerate}
\end{nonumtheorem}
Let us formulate a version of Hermite-Bieler Theorem for the left half-plane.
\begin{lemma} %
\label{VHBT}%
Let \(M\) be a polynomial with positive coefficients,
\[M(t)=\sum\limits_{0\leq{}k\leq{}n}m_{k}t^{k},
\quad{}m_k>0,\ 0\leq{}k\leq{}n\,,\] and let
\({}^{\mathscr{E}\!}M\) and \({}^{\mathscr{O}\!}M\) be the even
and the odd parts of \(M\). In order for the polynomial \(M\) be
dissipative it is necessary and sufficient that the following two
condition be satisfied:\vspace*{-2.6ex}
\begin{enumerate}%
\item
 The polynomials \({}^{\mathscr{E}\!}M\) and
\({}^{\mathscr{O}\!}M\) are conservative.
\item
The sets of roots of the polynomials  \({}^{\mathscr{E}\!}M\) and \({}^{\mathscr{O}\!}M\)
 interlace.
\end{enumerate}%
\end{lemma}%
\begin{lemma}
\label{CPtBC}%
 Let \(W\),
\begin{equation}
W(t)=w_{0}+w_{2}t^{2}+w_{4}t^{4}\,\cdots\,+w_{2m-2}t^{2m-2}+w_{2m}t^{2m}
\end{equation}
  be an even polynomial with positive coefficients \(w_{2l}\):
  \[w_{0}>0,\,w_{2}>0,\,\dots\,,\,w_{2m}>0\,.\]
  In order for the polynomial \(W\) to be conservative it is
  necessary and sufficient that the polynomial \(M=W+W^{\prime}\)
  to be dissipative, where \(W^{\prime}\) is the derivative of
  \(W\):
\begin{equation}
W^{\prime}(t)=2\cdot{}w_{2}t^{}+4\cdot{}w_{4}t^{3}\,\cdots\,
+(2m-2)\cdot{}w_{2m-2}t^{2m-3}+2m\cdot{}w_{2m}t^{2m-1}\,.
\end{equation}
\end{lemma}

 \textsf{PROOF of LEMMA \ref{VHBT}.} Let
\begin{equation}%
\label{LtoU}%
P(t)=M(it),\quad{}A(t)=({}^{\mathscr{E}\!}M)(it),
\quad{}B(t)=i\sp{-1}\cdot({}^{\mathscr{O}\!}M)(it),
\end{equation}%
so
\[P(t)=A(t)+iB(t)\,.\]
\(A\) and \(B\) are polynomials with real coefficients:
\[A(t)=\sum\limits_{0\leq{}l\leq\left[\frac{n}{2}\right]}(-1)^lm_{2l}t^{2l},\quad{}
B(t)=t\!\!\!\!\sum\limits_{0\leq{}l\leq\left[\frac{n-1}{2}\right]}\!\!\!\!(-1)^lm_{2l+1}t^{2l}\,.\]
Moreover,
\begin{equation}%
\label{ALN1}%
 B^{\prime}(0)A(0)-A^{\prime}(0)B(0)=m_0\,m_1\,.
\end{equation}%
From \eqref{LtoU} it is evident that
\begin{align*}
(\textup{\footnotesize{}All roots of }A\textup{ \footnotesize
{}are real and simple})&\Leftrightarrow(\textup{\footnotesize{}The
polynomial}{}^{\ \mathscr{E}\!}M\textup{ \footnotesize {}is conservative})\\
(\textup{\footnotesize{}All roots of }B\textup{ \footnotesize
{}are real and simple})&\Leftrightarrow(\textup{\footnotesize{}The
polynomial}{}^{\ \mathscr{O}\!}M\textup{ \footnotesize {}is conservative})\\
(\textup{\footnotesize{}All roots of }P\textup{ \footnotesize{}lie
in \(\{z\!:\textup{Im}\,z>0\}\)})&
\Leftrightarrow(\textup{\footnotesize{}The polynomial }M\textup{
\footnotesize{}is dissipative})
\end{align*}
and under condition that all roots of \(A\) and \(B\) are real,
\begin{equation*}
(\textup{\footnotesize{}The roots of \(A\) and \(B\)
interlace})\Leftrightarrow{}(\textup{\footnotesize{}The roots of \ %
\({}^{\mathscr{E}\!}M\) and \({}^{\mathscr{O}\!}M\) interlace})
\end{equation*}
Thus, Lemma \ref{VHBT} is an immediate consequence of
Hermite-Bieler Theorem it the above stated form. \textit{The
inequality \eqref{ALN} is ensured automatically by \eqref{ALN1}
since the
coefficients \(m_k\) are assumed to be positive.}\\ %
{\hspace*{0.ex} }\hfill Q.E.D.\\%
\textsf{PROOF of LEMMA \ref{CPtBC}.} It is clear that the
polynomials \(W\) and \(W^{\prime}\) are the even and the odd
parts of \(M=W+W^{\prime}\):
\[W={}^{\mathscr{E}}\!M,\quad{}W^{\prime}={}^{\mathscr{O}}\!M.\]
Let \(M\) be dissipative. Then, according to Lemma \ref{VHBT},
\(W\) is conservative. Conversely, let \(W\) be conservative.
According to Rolle theorem, the polynomial \(W^{\prime}\) is
conservative as well, and the sets of roots of \(W\) and
\(W^{\prime}\) interlace. By Lemma \ref{VHBT}, the polynomial
\(M\) is dissipative. {\ }\hfill Q.E.D.
\begin{remark}%
\label{NoNe}%
The claim of Lemma \ref{VHBT} remains true if to replace the
assumption posed on the coefficients \(m_k\) of of \(M\) with a
weaker assumption. It enough to assume that only the coefficients
\(m_0\) and \(m_1\) are strictly positive, whereas
 the other coefficients \(m_k,\,\, k=2,\,3,\,\dots\,,\,n,\) are real.
\end{remark}%
\textsf{PROOF of THEOREM \ref{H10H0}.} The relation \eqref{MPOP}
means that the polynomial \(tW_{\partial V}^{1}(t)\) is
the odd part of the Minkowski polynomial \(M_V\).
 Thus, we are in the situation of Lemma \ref{VHBT}.
 Since the polynomial
\(M_V\) is dissipative, the point \(z=0\) is not a root of \(M\),
that is \(m_0(V)\not=0\). According to \eqref{EvEq}, this means
that \(\textup{Vol}_n(V)\not=0\). Thus, the set \(V\) is solid. By
Proposition \ref{PCMP}, all the coefficients \(m_k(V)\) of the
polynomial \(M_V\) the are strictly positive. According to Lemma
\ref{VHBT}, the polynomial \({}^{\mathscr{O}\!}(M_V)\) is
conservative. Since \({}^{\mathscr{O}\!}(M_V)(0)=0\), the polynomial
\(t^{-1}\cdot{}{}^{\mathscr{O}\!}(M_V)(t)=W_{\partial V}^{1}(t)\) is conservative as well.%
{\hspace*{0.ex} }\hfill Q.E.D.%

\textsf{PROOF OF THEOREM \ref{JDM}.} In the course of the proof we
shall refer to some facts from the theory of entire functions
which usually are formulated  in literature for functions whose
roots are in the upper rather then in the left half-plane.
Therefore, it is convenient pass from the variable \(t\) to the
variable \(it\). Given a function \(H(t)\) of the Hurwitz class
\(\mathscr{H}\), let \(f(t)=H(it)\).  Then \(f\) is an entire
function of exponential type, all roots of \(f\) are in the upper
half-plane, and moreover, the defect \(d_f\) of \(f\) is
non-negative, where
\[2\,d_f=\varlimsup\limits_{r\to+\infty}f(-ir)-\varlimsup\limits_{r\to+\infty}f(-ir)\,.\]
(It is clear that \(d_f=d_H\), where \(d_H\) is the same as in
\eqref{DefH}.) Thus the function \(f\) is in the class \(P\) as
this class was defined in \cite{Lev1}, Chapter VII, Section 4. Let
\[f(t)=A(t)+iB(t)\,,\]
where \(A\) and \(B\) be real entire functions. Combining Lemma 1
from \cite{Lev1}, Chapter VII, Section 4 with Theorem 4 from
\cite{Lev1}, Chapter VII, Section 2, we obtain that the functions
\(A\) and \(B\) possess the properties:
\begin{enumerate}
\item
\(A\) and \(B\) are real entire functions of exponential type;
\item
\(A(0)B^{\prime}(0)-B(0)A^{\prime}(0)>0;\)
\item
For every \(\theta\in\mathbb{R}\), all roots of the linear
combination \(C_{\theta}\), where
\(C_{\theta}(t)=\cos\theta{}A(t)+\sin\theta{}B(t),\) are simple
and real.  (The entire functions \(A\) and \(B\) are \textit{a
real pair} in the terminology of N.G.Chebotarev, \cite{Cheb}.)
\end{enumerate}
According to Hadamard's factorization theorem, the entire function
\(C_{\theta}\) is in the Laguerre-Polya class. According to
Jensen-Csordas-Williamson Theorem (Theorem \ref{JSW}), for each
\(n\), all roots of the Jensen polynomial
\(C_{\theta,n}(t)=\mathscr{J}_n(C_{\theta};t)\) are real and
simple. Thus, the real polynomials \(A_n(t)=\mathscr{J}_n(A;t)\)
and \(B_n(t)=\mathscr{J}_n(B;t)\) possess the property:
\textit{For every \(\theta\in\mathbb{R}\), all roots of the linear
combination \(\cos\theta{}A_{n}(t)+\sin\theta{}B_{n}(t),\) are
real and simple}. (The polynomials \(A_n\) and \(B_n\) are real
pair as well.) From the property of the polynomials \(A_n\) and
\(B_n>0\) to be a real pair together with the property
\(A_n(0)B_n^{\prime}(0)-B_n(0)A_n^{\prime}(0)\) it follows that
all roots of the polynomial \(f_n(t)=A_n(t)+iB_n(t)\) are in the
upper half-plane. Thus, all roots of the polynomial
\(H_n(t)=f_n(-it)\) are in the left half-plane. In other words,
the polynomial \(H_n\) is a Hurwitz polynomial. On the other hand,
from the construction it follows that
\(H_n(t)=\mathscr{J}_n(H;t)\,.\) \hfill\framebox[0.45em]{ }
\section{PROPERTIES OF MINKOWSKI POLYNOMIALS
\newline
 OF A CONVEX SET. \label{ChPrMiPo}  }
 {\small\textsf{MOTION INVARIANCE}}: \textit{Let
\(V,\,V\subset \mathbb{R}^n\), be a compact convex set,
{\large\(\tau\)} be a motion\footnote{The motion of the space
\(\mathbb{R}^n\) is an affine transformation of \(\mathbb{R}^n\)
which preserves the Euclidean distance in \(\mathbb{R}^n\).} of
the space \(\mathbb{R}^n\)}, and {\large\(\tau\)}(V) be the image
of the set \(V\) under he motion {\large\(\tau\)}. Then
\(M_{{\mbox{\small\(\tau\)}}(V)}(t)=M_V(t)\).\\[1.0ex]
\noindent {\small\textsf{CONTINUITY}}: \textit{The correspondence
\(V\to{}M_V\) between compact convex sets \(V\) in
\(\mathbb{R}^n\) and their Minkowski polynomials \(M_V\) is
continuous \footnote{\label{topo}The set of compact convex sets in
\(\mathbb{R}^n\) is equipped by the Hausdorff metric, the set of
all polynomials is equipped by the topology of the locally uniform
convergence in \(\mathbb{C}\).}.}
\\[1.0ex]
A sketch of the proof of  the continuity property can be found in
\cite{BoFe},\,section \textsf{29}; \cite{BuZa},\,section
\textsf{19.2}; \cite{Schn},\,section \textsf{5.1}; \cite{Web},

\noindent {\small\textsf{MONOTONICITY}}:
\textit{%
Let \(V_1\) and \(V_2\) be compact convex sets in
\(\mathbb{R}^n\), and \(M_{V_1}, M_{V_2}\) be the appropriate
Minkowski polynomials.}
\textit{If \(V_1\subset{}V_2\), then the coefficients
\(m_k(V_1),\,m_k(V_2)\) of the polynomials \(M_{V_1}\),
\(M_{V_2}\), defined as in \eqref{MiP}, satisfy the inequalities}
\begin{equation}
\label{IMiP} m_k(V_1)\leq{}m_k(V_2),\quad{}0\leq k\leq n\,.
\end{equation}
\noindent%
\textsf{Explanation.} According to the definition of the
\textit{mixed volumes},
\begin{equation}
\label{MiX}%
 m_k(V)=\frac{n!}{(n-k)!\,k!}\,
\textup{Vol}(\underbrace{V,V,\,\dots\,,V}_{n-k};\,\underbrace{B^n,B^n,\,\dots\,,B^n}_{k}).
\end{equation}
Inequalities \eqref{IMiP} follow from the \textit{monotonicity of the
mixed volumes \eqref{MiX} with respect to \(V\)}. (Concerning the
monotonicity of the mixed volumes see, for example, \cite{BoFe},
section \textsf{29}; \cite{BuZa}, section \textsf{19.2};
\cite{Web}, Theorem \textsf{6.4.11}; \cite{Schn}, section
\textsf{5.1}, formula (5.1.23).)\hfill Q.E.D.\\
\begin{lemma}
\label{PCMP}
\begin{enumerate}
\begin{subequations}
\label{PosProp}
\item[\textup{a)}.]
 For any compact convex set
\(V,\,V\subset{}\mathbb{R}^n\), the coefficients \(m_k(V)\) of its
Minkowski polynomial, defined as in \eqref{MiP}, are non-negative:
\begin{equation}
\label{PPa} 0\leq{}m_k(V),\ \ 0\leq k\leq n.
\end{equation}
(According to \eqref{EvEq}, the coefficient \(m_n(V)\) is strictly
positive.)
\item[\textup{b)}.]
 If moreover the set \(V\) is solid, then
all coefficients \(m_k(V)\) are strictly  positive:
\begin{equation}
\label{PPb}%
 0<m_k(V),\ \ 0\leq k\leq n.
\end{equation}
\end{subequations}
The Weyl coefficients \(k_{2l}(\partial{}V),\,0\leq
l\leq\left[\frac{n-1}{2}\right]\), \,defined by \textup{Definition
\ref{DNWP}},
 are strictly positive as well.
\end{enumerate}
\end{lemma}

\textsf{PROOF.} Taking  \(V\) as \(V_2\) and an one-point subset
of \(V\) as \(V_1\) in \eqref{IMiP}, we obtain \eqref{PPa}. If the
set \(V\) is solid, then there exist a ball \(x_0+\rho{}B^n\) of
some positive radius \(\rho\): \(x_0+\rho{}B^n\subset{}V\). Taking
the ball \(x_0+\rho{}B^n\) as \(V_1\) and \(V\) as \(V_2\) in
\eqref{IMiP}, we obtain the inequalities
\(m_k(x_0+\rho{}B^n)\leq{}m_k(V),\,\,\,0\leq k\leq n.\) Moreover,
\(m_k(x_0+\rho{}B^n)=m_k(\rho{}B^n)=\rho^{n-k}m_k(B^n)=
\rho^{n-k}\frac{n!}{k!\,(n-k)!}\,\textup{Vol}_n(B^n)>0.\)
\hfill\framebox[0.45em]{ }

\vspace{2.0ex}
\begin{remark}%
\label{NSND2}%
The notion of the interior point of a set \(V\) depend on the
space in which \(V\) is embedded. The set \(V\),
\(V\subset\mathbb{R}^{n},\) which is non-solid with respect to the
`original' space \(\mathbb{R}^{n}\), is solid if \(V\) is
considered as be embedded in the space \(\mathbb{R}^{d}\) of the
`right' dimension \(d,\,d<n\). The dimension \(\dim{}V\) of the
set \(V\) should be taken as such \(d\).
\end{remark}
\begin{definition}
\label{DDCS} Let \(V,\,V\subseteq\mathbb{R}^n\), be a convex set.
The \textsf{dimension} \(\dim V\) \,of~\(V\) is the dimension of
the smallest affine subspace of \(\mathbb{R}^n\) which contains
\(V\).
\end{definition}
\begin{lemma}%
\label{Degdeg}%
Let \(V,\,V\subset\mathbb{R}^n\), be a compact convex set of the
dimension \(d\):
\begin{equation}
\label{DiCoS} \dim{}V=d,\quad{}0\leq{}d\leq{}n.
\end{equation}
Then
\begin{equation}
m_k^{\mathbb{R}^n}(V)=0\ \ \textup{for}\ \ 0\leq{}k<n-d;\ \
m_k^{\mathbb{R}^n}(V)>0\ \ \textup{for}\ \ n-d\leq{}k\leq{}n\,.
\end{equation}
\end{lemma}%
This lemma is a consequence of Lemma \ref{PCMP} and of the
following
\begin{lemma}
\label{CoRiDi}%
 Let \(V\), \(V\subset{}\mathbb{R}^n\), be a convex set of
 dimension \(d\), \(d\leq{}n\), and let
 \(M_V^{\mathbb{R}^n}(t)=\sum\limits_{0\leq{}k\leq{}n}m_k^{\mathbb{R}^n}(V)t^k\)
 and\,%
 \footnote{Defining the Minkowski polynomial
 \(M_V^{\mathbb{R}^d}\), we can assume that
 the
smallest affine subspace of \(\mathbb{R}^n\) which contains \(V\) is the space \(\mathbb{R}^d\).} %
  \(M_V^{\mathbb{R}^d}(t)=\sum\limits_{0\leq{}{\,k}\leq{}d}m_{\,k}^{\mathbb{R}^d}(V)t^k\)
 be the Minkovski polynomials of the set \(V\) with respect to the
 ambient spaces \(\mathbb{R}^n\) and \(\mathbb{R}^d\)
 respectively. Then
 \begin{equation}
M_V^{\mathbb{R}^n}(t)=t^{n-d}\cdot\hspace*{-1.0ex}\sum\limits_{0\leq{}{\,k}\leq{}d}
\,\frac{{\pi}^{\frac{n-d}{2}}%
\Gamma(\frac{k}{2}+1)}{\Gamma(\frac{k+n-d}{2}+1)}m_{\,k}^{\mathbb{R}^d}(V)t^k\,.
 \end{equation}
\end{lemma}
Lemma \ref{CoRiDi} appears in slightly different notation as
Theorem \ref{IRPp} in Section \ref{ExInSp}, where proof is
presented.

\begin{definition}
\label{CrSM}%
The mixed volumes appearing in \eqref{MiX} are said to be
\textsf{cross-sectional measures}  of the set \(V\) and are denote
 as \(v_{n-k}(V)\):
\begin{equation}
\label{CSM}
\textup{Vol}(\underbrace{V,V,\,\dots\,,V}_{n-k};\,\underbrace{B^n,B^n,\,\dots\,,B^n}_{k})=
v_{n-k}(V),\quad 0\leq k\leq n\,.
\end{equation}
\end{definition}
Thus, the coefficients of the Minkovski polynomials \(M_V\), which
appear in \eqref{MiP}, can be presented as %
\begin{equation}
\label{CrSFC}%
 m_k(V)=\binom{n}kv_{n-k}(V),\quad
\binom{n}k=\frac{n!}{k!\,(n-k)!} \text{ \small are binomial
coefficients,}%
\end{equation}
and the Minkowski polynomial itself can be presented as
\begin{equation}
\label{CrSFP}%
M_V(t)=\sum\limits_{0\leq k\leq n}\binom{n}kv_{n-k}(V)t^k,
\end{equation}
The following fact  will be used essentially in Section \ref{LoDi}: \\[2.5ex]
\noindent \textsf{ALEXANDROV\,--\,FENCHEL INEQUALITY.} \textit{Let
\(V\), \(V\subset\mathbb{R}^n\), be a compact convex set. Then its
cross-sectional measures \(v_k(V)\) satisfy the inequalities}
\vspace*{-1.5ex}%
\begin{equation}
\label{AFIn}%
v_k^2(V)\geq{}v_{k-1}(V)\,v_{k+1}(V),\ \ 1\leq k\leq n-1\,.
\end{equation}

\noindent%
A.D.\,Alexandrov published two proofs of this inequality in
\cite{Al1} and \cite{Al2}. The first of them, a combinatorial one,
is carried out for the convex polyhedra. The second proof is more
analytical. It uses the theory of self-adjoint elliptic operators
depending on parameter. This proof is carried out for smooth
convex bodies. To the general case, both proofs are generalized by
limit arguments. The first proof is developed in detail in the
textbook \cite{Le}. The second proof is reproduced in Busemann
\cite{Bus}. It has become customary to talk on `Alexandrov-Fenchel
inequality', because Fenchel \cite{Fen} also stated the inequality
and sketched  the proof. Its detailed exposition was never
published. At the end of 1978 independently Tessier in Paris and
A.G.Khovanski{\u\i} in Moscow obtained an algebraic-geometrical
proof of the Alexandrov-Fenchel inequality using the Hodge index
theorem. This proof is developed in \S 27 of the English
translation of \cite{BuZa} and was written by
A.G.\,Khovanski{\u\i}. (In the Russian original of \cite{BuZa} an
erroneous algebraic proof of the Alexandrov-Fenchel inequality was
included which has been excluded
 in the English translation.) Regarding the Alexandrov-Fenchel
 inequality see also \cite{BuZa}, \S\,20 and Section \textbf{6.3}
 of \cite{Schn}.
\begin{definition}%
\label{DefLogConc}%
A sequence \(\{p_k\}_{0\leq{}k\leq{}n}\) of non-negative numbers:
\begin{equation}%
\label{NoNeCo}%
p_k\geq{}0\,, \ \ \ 0\leq{}k\leq{}n,
\end{equation}%
 is said to be
\textsf{logarithmic concave}, if the following inequalities hold:
\begin{equation}
\label{LogCoIn}
p_k^2\geq{}p_{k-1}p_{k+1},\quad{}1\leq{}k\leq{}n-1\,.
\end{equation}
\end{definition}%
Thus, the Alexandrov-Fenchel inequalities can be formulated in the
form:\\[2.0ex]
\hspace*{3.0ex}
\begin{minipage}{0.9\linewidth}
 \textit{For any convex set \(V\), the sequence \(\{v_k(V)\}_{0\leq{}k\leq{}n}\) of
 its
cross sectional measures  is logarithmic concave.}
\end{minipage}

Under the extra condition \eqref{NoNeCo}, the logarithmic
concavity inequalities \eqref{LogCoIn} for the coefficients of the
polynomial
\begin{equation}
\label{TuPo}%
 P(t)=\sum\limits_{0\leq{}k\leq{}n}\binom{n}k\,p_kt^k,
\end{equation}
or for the coefficients of the entire function
\begin{equation}
\label{TuEf}
P(t)=\sum\limits_{0\leq{}k<\infty}\frac{p_k}{k!}\,t^{k},
\end{equation}
have been considered in connection with distribution roots of
\(P\). In this setting, such (and analogous) inequalities are
commonly known as \textsf{Tur\'an Inequalities}
(\textsf{Tur\'an-like Inequalities).} Concerning Tur\'an
inequalities see, for example, \cite{KaSc} and
\cite{CrCs2}.\\
\begin{remark}
\label{CoTiI}%
The Tur\'an inequalities \eqref{LogCoIn} for the coefficients of
the polynomial \eqref{TuPo} or entire function \eqref{TuEf} impose
some restrictions on location of roots of \(P\). However,
\textsf{these inequalities alone do not ensure} that all roots of
\(P\) are located in the left half-plane
\(\{z:\,\textup{Re}\,z<0\}.\)
\end{remark}
For example, given \(m\in\mathbb{N}\), let
\begin{equation}
\label{SE}%
p_k=1 \textup{ \ \small{}for \ }k=0, 1, \,\dots\,,\,m  \textup{ \
\small{}and \  }p_k=0\textup{ \ \small{}for \ }k>m.
\end{equation}
Such \(p_k\) satisfy the Tur\'an inequalities \eqref{LogCoIn}. The
function \eqref{TuEf} corresponding to these \(p_k\) is the
polynomial
\begin{equation}
\label{SES} %
P_m(t)=\sum\limits_{0\leq{}k\leq{}m}\frac{t^m}{m!}\,.
\end{equation}
This polynomial is a section of the exponential series. It is
known that already for \(m=5\) the polynomial \eqref{SES} has two
roots located in the half-plane \(\{z:\,\textup{Re}\,z>0\}.\) G.\,Szeg'o,
\cite{Sz}, studied the limiting distribution of roots of the
sequence of polynomials \(P_m\), \eqref{SES}, as \(m\to\infty\).
From his results on the limiting distributions of the roots it
follows that for large \(m\) the polynomial \(P_m\) not only has
roots in the half-plane \(\{z:\,\textup{Re}\,z>0\}\), but that the total
number of its roots located there has a positive density as
\(m\to\infty\). Regarding roots of sections of power series we
address to the book \textup{\cite{ESV}} and to the survey
\textup{\cite{Ost1}.}
 For \(m<n\), the polynomial
\eqref{TuPo} with \(p_k\) as in \eqref{SE} takes the form
\begin{equation}
\label{TBiP}%
P_{m,n}(t)=\sum\limits_{0\leq{}k\leq{}m}\binom{n}kt^k\,.
\end{equation}%
I.V.\,Ostrovskii, \textup{\cite{Ost2}}, studied the limiting
distribution of roots of the sequence of the polynomials
\(P_{m,n}\) as \(m,n\to\infty\),
\(m/n\to\alpha,\,\alpha\in(0,1).\) From his results it follows
that for large \(m,n\): \(n/m=O(1),\,n/(n-m)=O(1)\) the polynomial
\(P_{m,n}\) not only has roots in the half-plane
\(\textup{Re}\,z>0\), but that the total number of its roots
located there has a positive density as
\(m,n\to\infty,\,m/n\to\alpha\in(0,1)\).%

\vspace{3.0ex}
\section{ROUTH-HURWITZ CRITERION.\label{DeCr}}

Of course it is desirable to obtain an information about the
location of roots of the Weyl and Minkowski polynomials directly
from \textit{geometrical} considerations. At the moment we are not
able to do this. The only general tool from geometry which we can
use are the Alexandrov-Fenchel inequalities \eqref{AFIn} for
cross-sectional measures \(v_k(V)\) of convex sets. Therefore one
should express all polynomials which we investigate in terms of
this cross-sectional measures.

As it was explained in \eqref{CrSFP}, the expression of the
Minkowski polynomial \(M_V^{\mathbb{R}^n}\) for the convex set
\(V,\,V\subset\mathbb{R}^n\), in terms of the cross-sectional
measures \(v_k(V)\) is
\begin{equation}
\label{CrSFP1}%
M_V^{\mathbb{R}^n}(t)=\sum\limits_{0\leq k\leq
n}\binom{n}kv_{n-k}(V)t^k\,.
\end{equation}
\begin{lemma} %
\label{EWRC}%
 Let \(\mathscr{M}\) be a closed convex surface, \(\dim\mathscr{M}=n\), and let
 \(V\), \(V\subset\mathbb{R}^{n+1}\), be a generating convex set:
 \(\mathscr{M}=\partial{}V.\)

 Then the Weyl \(W_{\mathscr{M}}^{\infty}\) can be
 expressed as
\begin{equation}
\label{MaWePM}%
W_{\mathscr{M}}^{\infty}(t)=\sum\limits_{0\leq{}l\leq{}\left[\frac{n}{2}\right]}
\frac{(n+1)!}{2^ll!(n-2l)!}\,v_{n-2l}(V)\,t^{2l}\,,
\end{equation}
where \(v_k(V)\) are the cross-sectional measures of the
generating convex set \(V\).
\end{lemma}
 \textsf{PROOF}. The expression \eqref{MaWePM} is a consequence of
 \eqref{UWP}, \eqref{NWP} and \eqref{CrSFC}.
 \hfill\framebox[0.45em]{ }

\vspace{2.0ex}%
 To extract  an information about the location of
roots of the Minkowski polynomial \(M_V^{\mathbb{R}^n}\) from
\eqref{CrSFP1}, we may use the Alexandrov-Fenchel inequalities
\eqref{AFIn}. The Alexandrov-Fenchel inequalities relate the
cross-sectional measures \(v_k(V),\,v_{k-1}(V),\,v_{k+1}(V) \). To
extract such an information about the roots of the Weyl polynomial
\(W_{\mathscr{M}}^{\infty}\) from \eqref{MaWePM}, we need the
analogous inequalities which relate the cross-sectional measures
\(v_k(V\)),\, \(v_{k-2}(V)\),\, \(v_{k+2}(V)\).

\begin{lemma}
\label{CheO}
Let \(V,\,\,V\subset\mathbb{R}^{n+1}\), be a compact convex set.
Then its cross-sectional measures satisfy the inequalities
\begin{equation}
\label{AFlI}%
 v_k^2(V)\geq{}v_{k-2}(V)v_{k+2}(V), \quad
2\leq{}k\leq{n-1}.
\end{equation}
\end{lemma}
\textsf{PROOF}. We derive \eqref{AFlI} from \eqref{AFIn}. Rising
the inequality \eqref{AFIn} to square, we obtain that
\(v_{\,k}^4(V)\geq{}v_{k-1}^{\,\,2}(V)v_{k+1}^{\,\,2}(V).\)
Inequalities \eqref{AFIn} with \(k\) replaced with \(k-1\) and
\(k+1\) are: \[v_{k-1}^2(V)\geq{}v_{k-2}^{}(V)v_{k}^{}(V) %
\textup{ and } v_{k+1}^2(V)\geq{}v_{k}^{}(V)v_{k+2}^{}(V)\]
respectively. The inequality \eqref{AFlI} is the consequence of
the last three inequalities. {\ }\hfill Q.E.D.

The inequalities \eqref{AFIn} and \eqref{AFlI} for the
coefficients of the polynomials \eqref{CrSFP1} and \eqref{MaWePM}
respectively is one of two general tools which will be used in the
study of the location of roots of these polynomials. The second
general tool is the criteria which express the properties of
polynomials to be dissipative (or conservative respectively) in
terms of their coefficients. Such criteria are formulated as the
positivity of certain determinants constructed from the
coefficients of the tested polynomials.
\begin{nonumtheorem}[[Routh-Hurwitz\textup{]}]
Let
 \begin{equation}
 \label{RaHuPo}
 P(t)=a_0t^n+a_1t^{n-1}+\,\dots{}\,+a_{n-1}t+a_n
 \end{equation}
  be a
polynomial with strictly positive coefficients:
\begin{equation}
 \label{PoCo}
a_0>0,\,a_1>0,\,\dots\,,\,a_{n-1}>0,\,a_n>0\,.
\end{equation}
For the polynomial \(P\) to be dissipative, it is necessary and
sufficient that all the determinants
\(\Delta_k,\,k=1,\,2,\,\dots\,,\,n-1,\,n,\) be strictly positive:
\begin{equation}
\label{PRHD}%
\Delta_1>0,\,\Delta_2>0,\,\dots\,,\,\Delta_{n-1}>0,\,\Delta_n>0,
\end{equation}
where {\footnotesize
\begin{multline}
\label{RHD} \Delta_1=a_1,\quad\Delta_2=
\begin{vmatrix}a_1&a_3\\
a_0&a_2
\end{vmatrix},\quad
\Delta_3=
\begin{vmatrix}a_1&a_3&a_5\\
a_0&a_2&a_4\\
0&a_1&a_3
\end{vmatrix},\quad \\
\Delta_4=
\begin{vmatrix}a_1&a_3&a_5&a_7\\
a_0&a_2&a_4&a_6\\
0&a_1&a_3&a_5\\
0&a_0&a_2&a_4
\end{vmatrix},\quad\dots
\quad ,\ \Delta_n=
\begin{vmatrix}a_1&a_3&a_5&\dots&0\\
a_0&a_2&a_4&\dots&0\\
0&a_1&a_3&\dots&0\\
\hdotsfor{5} \\
\hdotsfor{4}&a_n%
\end{vmatrix}\cdot
\end{multline}
}
\end{nonumtheorem}
This result, as well as many relative results, can be found in
\cite{Gant}, Chapter XV. See also \cite{KrNa}.
\begin{remark}
\label{LieCh}%
 Actually, to prove that the polynomial \(P\), \eqref{RaHuPo}, of
 degree \(n\) with positive coefficients \(a_k,\,k=1,\,2,\,\dots\,,\,n,\) is dissipative,
  there is no need to inspect \textit{all} Hurwitz determinants
 \(\Delta_k,\,k=1,\,2,\,\dots\,,\,n,\) for positivity. It is
 enough to inspect either the determinants \(\Delta_k\) only with even
 \(k\), or the determinants \(\Delta_k\) only with even
 \(k\). \textup{(See \cite{Gant}, Chapter XV, \S 13.)}
\end{remark}
Applying the Routh-Hurwitz criterion to investigate whether the Minkowski
polynomial \(M_V^{\mathbb{R}^n}\) is dissipative, we should take, according to
\eqref{CrSFP1},
\begin{equation}%
 \label{RHuC}%
 a_k=\frac{n!}{k!(n-k)!}v_k(V)\,, \quad 0\leq{}k\leq{}n, \ \ a_k=0, \quad
 k>n\,.
 \end{equation}%
From the criterion of dissipativity, the criterion of
conservativity can be derived easily.
\begin{nonumtheorem}[[Criterion of conservativity\textup{]}]
Let
\begin{equation}
\label{CoPoA}
P(t)=a_0t^{2m}+a_{2}t^{2m-2}+\,\dots{}\,+a_{2m-2}t^{2}+a_{2m}
\end{equation}
be a polynomial with strictly positive coefficients
\(a_{2l},\,0\leq{}l\leq{}m\):
\begin{equation}
\label{PCFA}%
 a_0>0,\,a_2>0,\,\dots\,,\,a_{2m-2}>0,\,a_{2m}>0\,.
\end{equation}
For the polynomial \(P\) to be conservative, it is necessary and
sufficient that all the determinants
\(D_k,\,k=1,\,2,\,\dots\,,2m-1,\,2m,\) be strictly positive:
\begin{equation}
\label{PRHDc}%
D_1>0,\,D_2>0,\,D_3>0,\,\dots\,,\,D_{2m-1}>0,\,D_{2m}>0,
\end{equation}
where the determinants \(D_k\) are constructed from the
coefficients  of the polynomial \(P\) according to the following
rule. Determinant \(D_k\) are the determinant \(\Delta_k\),
\eqref{RHD}, whose entries \(a_{2l},\,\,0\leq{}l\leq{}m,\)  are
the coefficients of the polynomial \(P\), and
\(a_{2l+1}=(m-l)\,a_{2l}\), \(0\leq{}l\leq{}m-1\): {\footnotesize
\begin{multline*}
  D_1=\textstyle{}m\,a_0,\quad{}D_2=
\begin{vmatrix}\textstyle
m\,a_0&(m-1)a_2\\
a_0&a_2
\end{vmatrix},\quad
D_3=
\begin{vmatrix}ma_0&(m-1)a_2&(m-2)a_4\\
a_0&a_2&a_4\\
0&m\,a_0&(m-1)a_2
\end{vmatrix},
\end{multline*}
\begin{equation}
\label{RHCo} D_4=
\begin{vmatrix}m\,a_0&(m-1)\,a_2&(m-2)a_4&(m-4)a_6\\
a_0&a_2&a_4&a_6\\
0&m\,a_0&(m-1)\,a_2&(m-2)a_4\\
0&a_0&a_2&a_4
\end{vmatrix},\quad\dots
\end{equation}
\begin{equation*}
 D_{2m}=
\begin{vmatrix}m\,a_0&(m-1)a_2&(m-2)a_4&\dots&0\\
a_0&a_2&a_4&\dots&0\\
0&ma_0&(m-1)a_2&\dots&0\\
\hdotsfor{5} \\
\hdotsfor{4}&a_{2m}%
\end{vmatrix}\cdot
\end{equation*}
}
\end{nonumtheorem}
\textsf{PROOF.} This theorem is the immediate consequence of
Hermite-Bieler theorem and Lemma \ref{CPtBC}.{\ }\hfill{}Q.E.D.

Applying the Conservativity criterion to investigate whether the Weyl
polynomial \(W_{\mathscr{M}}^{\infty}\) is dissipative, we should take, according to
 \eqref{MaWePM} and \eqref{CoPoA},
\begin{multline}%
 \label{RHuCl}%
 a_{2l}=\frac{(n+1)!}{2^{m-l}(m-l)!(2l+n-2m)!}v_{2l+n-2m}(V)\,,\\
 \quad 0\leq{}l\leq{}m, \ \ a_{2l}=0, \quad
 l>m\,,\quad \textup{where }m=\textstyle{\left[\frac{n}{2}\right]}
 \end{multline}%

\section{THE CASE OF LOW DIMENSION:\\
PROOF OF THEOREMS \ref{LDC} AND \ref{LDCW}.\label{LoDi}}

\textsf{PROOF of THEOREM \ref{LDC}.} We apply the Routh-Hurwitz
criterion of dissipativity to the Minkowski polynomial
\(M_{\,\,V}^{{\mathbb{R}^n}}\). `Opening' the Hurwitz determinants
\(\Delta_k\), \eqref{RHD} , and taking into account that \(a_k=0\)
for \(k>n\), we obtain that for \(n\leq{}5\),
\begin{subequations}
\label{ODet}
\begin{gather}
\Delta_1=a_1\tag{\theequation.\footnotesize{1}}\label{ODet1}\\
\Delta_2=a_1a_2-a_0a_3\,,\tag{\theequation.\footnotesize{2}}\label{ODet2}\\
\Delta_3=a_1a_2a_3+a_0a_1a_5-a_0a_3^2-a_1^2a_4\,,\tag{\theequation.%
\footnotesize{3}}\label{ODet3}\\
\Delta_4=a_1a_2a_3a_4+a_0a_2a_3a_5+2a_0a_1a_4a_5-a_1^2a_4^2-a_0^2a_5^2-a_0a_3^2a_4-a_1a_2^2a_5,
\tag{\theequation.%
\footnotesize{4}}\label{ODet4}\\
\Delta_5=a_5\Delta_4\,,\tag{\theequation.%
\label{ODet5}\footnotesize{5}}
\end{gather}
\end{subequations}
where we should take \(a_k\) as in \eqref{RHuC}.

 According to Routh-Hurwitz criterion, we have to prove
 that \(\Delta_1>0,\ \Delta_2>0,\,\dots\,,\,\Delta_n>0\).
 The cases \(n=2,\,3,\,4,\,5\) will be considered separately.
 Since \(V\) is solid, \(v_k(V)>0,\,0\leq{}k\leq{}n\).
  (Corollary \ref{PCMP}.b  and \eqref{CrSFC}\,.)\\
  Thus, the determinant \(\Delta_1=\binom{n}{1}v_1(V)\) is always
  positive.\\[1.0ex]
  The cases \(n=2,\,3,\,4,\,5\) will be considered separately.
  To shorten notation, we right \(v_k\) instead
  \(v_k(V)\).\\[1.0ex]
  \hspace*{3.0ex}\(\boldsymbol{n=2.}\) In this case,
  \[a_0=v_0,\,a_1=2v_1,\,a_2=v_2,\]
  \[\Delta_2=2v_2v_1\,.\]
  The inequality \(\Delta_2>0\) is evident. Thus,
  the Minkowski polynomial \(M_{V}^{\mathbb{R}^2}\) is dissipative.
  \\[1.0ex]
  \hspace*{3.0ex}\(\boldsymbol{n=3.}\) In this case,
  \[a_0=v_0,\quad{}a_1=3v_1,\quad{}a_2=3v_2,\quad{}a_3=v_3\,,\quad{}a_k=0,\,k>3\,.\]
  Substituting these expressions for \(a_k\) into \eqref{ODet}, we
  obtain
  \[\Delta_2=9v_1v_2-v_0v_3,\quad{}\Delta_3=v_3\Delta_2\,.\]
  Thus, the property of \(M_V^{\mathbb{R}^3}\) to be dissipative
  is equivalent to the inequality
  \begin{subequations}
  \label{MP3}
  \begin{align}
 9v_1v_2&>v_0v_3,\tag{\ref{MP3}}
\end{align}
  \end{subequations}
  \hspace*{3.0ex}\(\boldsymbol{n=4.}\) In this case,
  \[a_0=v_0,\quad{}a_1=4v_1,\quad{}a_2=6v_2,\quad{}a_3=4v_3\,\quad{}a_4=v_4,\quad{}a_k=0,\,k>4\,.\]
 Substituting these expressions for \(a_k\) into \eqref{ODet}, we
  obtain
  \begin{multline*}
  \Delta_2=24v_1v_2-4v_0v_3,\quad{}\Delta_3=96v_1v_2v_3-16v_0v_3^2-16v_1^2v_4\,,\quad
\Delta_4=v_4\Delta_3\,.
  \end{multline*}
Thus, the property of \(M_V^{\mathbb{R}^4}\) to be dissipative
  is equivalent to the pair of inequalities
  \begin{subequations}
  \label{MP4}
  \begin{align}
  6v_1v_2&>v_0v_3,\tag{\ref{MP4}.\footnotesize{2}}\label{MP42} \\
  6v_1v_2v_3&>v_0v_3^2+v_1^2v_4\,.\tag{\ref{MP4}.\footnotesize{3}}
  \label{MP43}
  \end{align}
  \end{subequations}
  \hspace*{3.0ex}\(\boldsymbol{n=5.}\) In this case,
  \[a_0=v_0,\ a_1=5v_1,\ a_2=10v_2,\ a_3=10v_3,\,\ a_4=5v_4,\ %
  a_5=v_5,\,\,a_k=0,\,k>5\,.\]
 Substituting these expressions for \(a_k\) into \eqref{ODet}, we
  obtain
  \begin{multline*}
  \Delta_2=50v_1v_2-10v_0v_3,\quad{}\Delta_3=500v_1v_2v_3+5v_0v_1v_5-100v_0v_3^2-125v_1^2v_4\,,\\
\Delta_4=2500v_1v_2v_3v_4+100v_0v_2v_3v_5+50v_0v_1v_4v_5\\
-625v_1^2v_4^2-v_0^2v_5^2-500v_0v_3^2v_4 -500v_0v_1^2v_5\,,\ \
\Delta_5=v_5\Delta_4\,.
  \end{multline*}
Thus, the property of \(M_V^{\mathbb{R}^5}\) to be dissipative
  is equivalent to the triple of inequalities
  \begin{subequations}
  \label{MP5}
  \begin{align}
  5v_1v_2&>v_0v_3,\tag{\ref{MP5}.\footnotesize{2}}
  \label{MP52} \\
  100v_1v_2v_3+v_0v_1v_5&>20v_0v_3^2+25v_1^2v_4\,,\tag{\ref{MP5}.\footnotesize{3}}
  \label{MP53}
  \\
  2500v_1v_2v_3v_4+100v_0v_2v_3v_5+&50v_0v_2v_4v_5>{}\tag{\ref{MP5}.\footnotesize{4}}
  \label{MP54}\\
 {} >625v_1^2v_4^2&+500v_0v_3^2v_4+500v_1v_2^2v_5+v_0^2v_5^2\,.\notag
  \end{align}
  \end{subequations}
  As it is claimed in Lemma \ref{CoAFI} below,
  the inequalities \eqref{MP3}, \eqref{MP4},
  \eqref{MP5}, where \(v_k=v_k(V)\) are the cross-sectional
  measures of the solid compact set \(V\) of the appropriate
  dimension, are the consequences of the Alexandrov-Fenchel
  inequalities. This completes the proof. {\ }\hfill Q.E.D.
  \begin{remark} All this business works up to certain \(n\), but it
\textsf{\textsl{does not}} work for all \(n\). If \(n\) is large
enough, then the conditions \( v_k^2\geq{}v_{k-1}v_{k+1},\quad
1\leq{}k\leq{}n-1\,,\) posed on positive numbers \(v_k\) does not
imply the inequalities \(\Delta_k\geq{}0\) for all
\(k=1,\,\dots,\,n\), where \(\Delta_k\) is constructed from
\(a_k=\binom{n}{k}v_k\).
 Already by \(n=30\),  \(\Delta_5<0\) for certain
\(v_k\) satisfying these conditions.
 Moreover, as we will see later,
for \(n\)  large enough, there exist examples of such compact
convex sets \(V\subset\mathbb{R}^n\) for which the Minkowski
polynomial \(M_V\) is not dissipative, and the Weyl polynomial
\(W_{\partial{}V}^1\) is not conservative. In such examples the
sets \(V\) are although solid, but `almost degenerated'.
\end{remark}
  \begin{lemma}
  \label{CoAFI}
  Let \(v_k,\,0\leq{}k\leq{}n,\) be strictly positive numbers
  satisfying the inequalities
  \begin{equation}
  \label{CoFVk}
  v_k^2\geq{}v_{k-1}v_{k+1},\quad 1\leq{}k\leq{}n-1\,.
  \end{equation}
  Then:\\[1.0ex]
  \hspace*{2.0ex}\textup{1}. If \(n=3\), then the inequality \eqref{MP3} holds;\\[1.0ex]
  \hspace*{2.0ex}\textup{2}. If \(n=4\), then the inequalities \eqref{MP4} holds;\\[1.0ex]
  \hspace*{2.0ex}\textup{3}. If \(n=5\), then the inequalities \eqref{MP5} holds.
  \end{lemma}
  \textsf{PROOF of LEMMA \ref{CoAFI}.} Given \(k,\,1\leq{}k\leq{}n-2\),
   multiplying the inequalities
  \[v_k^2\geq{}v_{k-1}v_{k+1},\quad{}v_{k+1}^2\geq{}v_{k}v_{k+2},\]
  an then cancelling on \(v_kv_{k+1}\), we obtain the inequality
  \begin{equation}%
  \label{pain}%
  v_kv_{k+1}\geq{}v_{k-1}v_{k+2}\,,\quad 1\leq{}k\leq{}n-2\,.
  \end{equation}%
  In particular, for \(n=3,\,k=1\), as well as for \(n=4,\,k=1\) and \(n=5,\,k=1\).
  \[v_1v_2\geq{}v_0v_3\,.\]
This inequality implies the inequality \eqref{MP3}, \eqref{MP42}
and \eqref{MP52}. Multiplying the inequality
\(v_1v_2\geq{}v_0v_3\) with the positive number \(v_3\), we obtain
\begin{subequations}
\label{MP4C}
\begin{equation}
\label{MP4C1}
v_1v_2v_3\geq{}v_0v_3^2\,
\end{equation}
For \(n=4,\,k=2\), the inequality \eqref{pain} means
\[v_2v_3\geq{}v_1v_4.\]
Multiplying the inequality  \(v_2v_3\geq{}v_1v_4\) with  \(v_1\),
we get
\begin{equation}
\label{MP4C2}
v_1v_2v_3\geq{}v_1^2v_4.
\end{equation}
\end{subequations}
The inequalities \eqref{MP4C} imply the inequality \eqref{MP43} and \eqref{MP53} .\\
Multiplying the inequality \(v_2v_3\geq{}v_1v_4\) {\small{}(this
is \eqref{pain} for \(k=2, n=5\))} with \(v_1v_4\), we obtain
\begin{subequations}
\label{MP5C}
\begin{equation}
\label{MP5C1}
v_1v_2v_3v_4\geq{}v_1^2v_4^2\,.
\end{equation}
Multiplying the inequality \(v_1v_2\geq{}v_0v_3\) {\small{}(this
is \eqref{pain} for \(k=1, n=5\))}  with \(v_3v_4\), we obtain
\begin{equation}
\label{MP5C3}
v_1v_2v_3v_4\geq{}v_0v_3^2v_4\,.
\end{equation}
Multiplying the inequality \(v_3v_4\geq{}v_2v_5\) {\small{}(this
is \eqref{pain} for \(k=3, n=5\))} with the positive number
\(v_1v_2\), we obtain
\begin{equation}
\label{MP5C31} v_1v_2v_3v_4\geq{}v_1v_2^2v_5\,.
\end{equation}
Further, multiplying the inequalities \(v_1v_2\geq{}v_0v_3\)
{\small{}(this is \eqref{pain} for \(k=1, n=5\))} and
\(v_3v_4\geq{}v_2v_5\) {\small{}(this is \eqref{pain} for \(k=3,
n=5\))}, we obtain that \(v_1v_4\geq{}v_0v_5\). Rising it to
square, we obtain
\[v_1^2v_4^2\geq{}v_0^2v_5^2\,.\]
 Multiplying this inequality with the inequality \(v_2v_3\geq{}v_1v_4\)
 {\small{}(this
is \eqref{pain} for \(k=2, n=5\))}, we obtain that
\begin{equation}
\label{MP5C2}
v_1v_2v_3v_4\geq{}v_0^2v_5^2\,.
\end{equation}
\end{subequations}
Since \(2500>625+500+500+1\), the inequality \eqref{MP54} it the
consequence of the inequalities \eqref{MP5C}. (The number \(2500\)
is the coefficient before the monomial \(v_1v_2v_3v_4\) in the
left hand side of \eqref{MP54}, the numbers \(625\), \(500\),
\(500\), \(1\) are the coefficients before the monomials
\(v_1^2v_4^2\), \(v_0v_3^2v_4\), \(v_1v_2^2v_5\) and
\(v_0^2v_5^2\) in the right hand side of \eqref{MP54}
respectively.)
 {\ } \hfill Q.E.D.

\textsf{PROOF of THEOREM \ref{LDCW}.} We apply the
\textsf{Criterion of conservativity}, which was formulated in the
previous section, to the Weyl polynomial
\(W_{\mathscr{M}}^{\infty}\). Opening the determinants
\(\Delta_k\), \eqref{RHCo} , and taking into account that
\(a_{2l}=0\) for \(l>\left[\frac{n}{2}\right]\), we obtain that
for \(n\leq{}5\),
that is\,%
\footnote{Recall that \(n=\dim \mathscr{M},\,n+1=\dim V:
\mathscr{M}=\partial V.\)}
 for \(m=\left[\frac{n}{2}\right]\leq{}2\),
\begin{subequations}
\label{DDet}
\begin{gather}
D_1=ma_0\tag{\theequation.\footnotesize{1}}\label{DDet1}\\
D_2=a_0a_2\,,\tag{\theequation.\footnotesize{2}}\label{DDet2}\\
D_3=a_0\big((m-1)a_2^2-2ma_0a_4\big)\,,\tag{\theequation.%
\footnotesize{3}}\label{DDet3}\\
D_4=a_0a_4(a_2^2-4a_0a_4)\,.
\tag{\theequation.%
\footnotesize{4}}\label{DDet4}
\end{gather}
\end{subequations}
where we should take \(a_{2l}\) as in \eqref{RHuCl}.

According to the Conservativity criterion, we have to prove that
\(D_1>0\), \(D_2>0,\,\dots\,,\,D_{2m}>0\), where
 \(m=\left[\frac{n}{2}\right]\).\\
Since \(V\) is solid, \(v_k(V)>0,\,0\leq{}k\leq{}n+1\).
  (Corollary \ref{PCMP}.b  and \eqref{CrSFC}\,.)
Thus, the determinants \(D_1,\,D_2\) are always positive.\\[1.0ex]
Therefore, if \(n=2\), or if \(n=3\), that is if \(m=1\), the Weyl
polynomial \(W_{\mathscr{M}}^{\infty}\) is conservative. Of
course, this fact is evident without referring to the
conservativity criterion:\\[1.0ex]
In the case \(\boldsymbol{n=2}\), according to \eqref{RHuC} or \eqref{MaWePM},
\[W_{\mathscr{M}}^{\infty}(t)=3v_2+3v_0t^2\,.\]
In the case \(\boldsymbol{n=3}\), according to \eqref{RHuC} or \eqref{MaWePM},
\[W_{\mathscr{M}}^{\infty}(t)=v_3+3v_1t^2\,.\]
Evidently, in both cases, \(n=2\) or \(n=3\), the polynomial \(W_{\mathscr{M}}^{\infty}\)
is  conservative.\\[1.0ex]
In the cases \(n=4,\ n=5\), to what corresponds \(m=2\),
\[D_3=a_0(a_2^2-4a_0a_4),\quad D_4=a_4D_3\,.\]
According to \eqref{RHuCl}, we have to take in the cases: \\
In the case \(\boldsymbol{n=4}\)
\[a_0=15v_0,\quad{}a_2=30v_2,\quad{}a_4=5v_4\,.\]
Thus,
\[D_3=15v_0(900v_2^2-300v_0v_4)\,.\]
The conditions \(D_3>0,\,D_4>0\) take the form
\begin{equation}
\label{WePo4}
3v_2^2>v_0v_4\,.
\end{equation}
In the case \(\boldsymbol{n=5}\)
\[a_0=90v_1,\quad{}a_2=60v_3,\quad{}a_4=6v_5\,.\]
Thus,
\[D_3=90v_1(900v_2^2-300v_0v_4)\,.\]
The conditions \(D_3>0,\,D_4>0\) take the form
\begin{equation}
\label{WePo5}
5v_3^2>3v_1v_5\,.
\end{equation}
So, in the cases \(n=4\) and \(n=5\) the property of the Weyl
polynomial \(W_{\mathscr{M}}^{\infty}\) be conservative is
equivalent to the inequality  \eqref{WePo4} and \eqref{WePo5}
respectively, where \(v_k=v_k(V)\) are the cross-sectional
  measures of the solid compact set \(V\) generating the surface
  \(\mathscr{M}:\ \mathscr{M}=\partial{}V.\) In its turn, the inequalities
\eqref{WePo4} and \eqref{WePo5} are evident consequences of the inequalities
\(v_2^2\geq{}v_0v_4\) and \(v_3^2\geq{}v_1v_5\) respectively.
 The latter inequalities are special cases
of the inequalities \eqref{AFlI}. (See Lemma \ref{CheO}.) Thus, in
the cases \(n=4\) and \(n=5\) the Weyl
\(W_{\mathscr{M}}^{\infty}\) polynomial of infinite index is
conservative. By Lemma \ref{LWPo}, all Weyl polynomials
\(W_{\mathscr{M}}^{\,p},\,p=1,\,2,\,3,\,\dots\,\,,\) are
 conservative as well.\hfill\framebox[0.45em]{ }

\section{EXTENDING OF THE AMBIENT
SPACE.\label{ExInSp}}
\paragraph{Adjoint convex sets.}
Let \(V\) be a compact convex set, \(V\subset\mathbb{R}^n\).
However we may consider the space \(\mathbb{R}^n\) as a subspace
of the space \(\mathbb{R}^{n+q}\) of a higher dimension
 \(q=1,\,2,\,3,\,\dots\,\). The
embedding \(\mathbb{R}^n\) to   \(\mathbb{R}^{n+q}\) is standard:
\begin{equation*}%
\label{Emb}%
\mathbb{R}^{n}\hookrightarrow \mathbb{R}^{n+q}:\ \ %
(\xi_1,\,\dots\,,\,\xi_n)\to(\xi_1,\,\dots\,,\,\xi_n;%
\underbrace{\,0,\,\dots\,,\,0}_q\,).
\end{equation*}%
 Thus, the set \(V\), which originally was considered as a
subset of \(\mathbb{R}^n\), may also be considered as a subset of
\(\mathbb{R}^{n+q}\). In other words, we identify the set
\(V\subset\mathbb{R}^n\) with the set \(V\times{}0_{q}\), which is
the Cartesian product of the set \(V\) and the zero point
\(0^{q}\) of the space \(\mathbb{R}^q\):
\(V\times{}0^{q}\subset\mathbb{R}^{(n+q)}.\)
\begin{definition}%
\label{DeAdM}
 Given a compact convex set \(V\), \(V\subset{}\mathbb{R}^n\),
and a number \(q\),
 \(q=0,\,1,\,2,\,3,\,\dots\,,\,\,\) the \textsf{\(q\)-th adjoint to \(V\)
set \(V^{(q)}\)} is defined as:
\begin{equation}
\label{AdS}
V^{(q)}\stackrel{\textup{\tiny{}def}}{=}V\times{}0^{q},\quad
V^{(q)}\subset\mathbb{R}^{n+q}\,,
\end{equation}
where \(0^q\) is the zero point of the space \(\mathbb{R}^{q}\),
and the space \(\mathbb{R}^{n+q}\) is considered as the Cartesian product:
\(\mathbb{R}^{n+q}=\mathbb{R}^{n}\times\mathbb{R}^{q}\).\\
The Minkowski polynomial
\(M_{\,V\times{}0^q}^{\,\mathbb{R}^{n+q}}\) of the \(q\)-th
adjoint set \(V^{(q)}\),
\begin{equation}
\label{defAMP}
M_{\,V\times{}0^q}^{\mathbb{R}^{n+q}}=\textup{Vol}_{n+q}(V\times{}0^q+tB^{n+q}),
\end{equation}
 is said to be the \textsf{\(q\)-th adjoint
Minkowski polynomial for the set~\(V\)}.
\end{definition}
For \(q=0\), the set \(V^{(0)}\) coincides with \(V\), and the polynomial
 \(M_{\,V\times{}0^0}^{\mathbb{R}^{n+0}}\) coincides with  \(M_{\,\,V}^{\mathbb{R}^{n}}\).
 For \(q=1\), the set \(V^{(1)}\) is what we are called the
 squeezed cylinder with the base \(V\).
\paragraph{Minkowski polynomials for adjoint sets.}
Let us answer the natural question:\\[1.0ex]
\centerline{\textit{How the polynomials
\(M_{\,V}^{\mathbb{R}^{n}}(t)\) and
\(M_{\,V\times{}0^q}^{\mathbb{R}^{n+q}}(t)\) are related}?}\\ 

 The answer  this question will be done by an inductive reasoning.
 Lemma~\ref{IPRom} below provides the step of the induction.
\begin{lemma}
\label{IPRom}%
Let \(V\) be a compact convex set in \(\mathbb{R}^n\), and
\begin{equation}%
\label{IRPom1}%
 M_{\,V}^{\mathbb{R}^n}(t)=\sum\limits_{0\leq{}k\leq{}n}m_{\,k}^{\mathbb{R}^n}(V)t^k
\end{equation}%
be the Minkowski polynomial with respect to the ambient space
\(\mathbb{R}^n\). Then the Minkowski polynomial
\(M_{\,V\times{}0}^{\mathbb{R}^{n+1}}(t)\) is equal to:
\begin{equation}%
\label{MPn1om}%
M_{\,V\times{}0^1}^{\mathbb{R}^{n+1}}(t)=%
t\!\!\sum\limits_{0\leq{}k\leq{}n}\frac{{\pi}^{1/2}
\Gamma(\frac{k}{2}+1)}{\Gamma(\frac{k+1}{2}+1)}\,m_{\,k}^{\mathbb{R}^n}(V)%
\, t^k\,.
\end{equation}%
\end{lemma}
The following theorem completes the inductional reasoning:
\begin{theorem}
\label{IRPp}%
Let \(V\) be a compact convex set in \(\mathbb{R}^n\), and
\begin{equation}%
\label{MPn}%
M_{\,V}^{\mathbb{R}^n}(t)=
\sum\limits_{0\leq{}k\leq{}n}m_{\,k}^{\mathbb{R}^n}(V)t^k
\end{equation}%
 be the Minkowski polynomial of the set
\(V\). Then the \(q\)-th adjoint Minkowski polynomial
\(M_{\,V\times{}0^q}^{\,\mathbb{R}^{n+q}}(t)\) of the set \(V\) is
equal to:
\begin{gather}%
\label{MPn2}%
M_{\,V\times{}0^q}^{\,\mathbb{R}^{n+q}}(t)=%
\sum\limits_{0\leq{}k\leq{}n}m_{\,k}^{\mathbb{R}^n}(V)\,\gamma_{\,\,k}^{(q)}\,\,%
t^{k+q}\,,\\
\intertext{where}
\label{MuSe}
\gamma_{\,k}^{(q)}={\pi}^{q/2}\,\frac{\Gamma(\frac{k}{2}+1)}{\Gamma(\frac{k+q}{2}+1)}\,,
\quad k=0,\,1,\,2,\,\dots\,\,;\ \ q=0,\,1,\,2,\,\dots.
\end{gather}
\end{theorem}
A sketch of proof of this theorem can be found in \cite{Had},
Chapter VI, Section 6.1.9. A detailed proof is presented below.
\begin{remark}%
\label{FnTp}%
Theorem \ref{IRPp} means that the sequence of the coefficients\\
\(\{m_{\,k}^{\mathbb{R}^{n+q}}(V\times{}0^q)\}_{0\leq k\leq n+q}\)
of the polynomial \(M_{\,V\times{}0^q}^{\mathbb{R}^{n+q}}\):
\begin{equation}%
M_{\,V\times{}0^q}^{\,\mathbb{R}^{n+q}}(t)=
\sum\limits_{0\leq{}k\leq{}n+q}m_{\,k}^{\mathbb{R}^{n+q}}(V\times{}0^q)\,t^k
\end{equation}%
 are obtained from the sequence of the coefficients
\(\{m_{k}(V)\}_{0\leq k\leq n}\) of the polynomial \(M_{\,V}^{\mathbb{R}^{n}}\),
\textup{\eqref{MPn}}, by means of shift and multiplication:
\begin{subequations}
\label{SiMu}
\begin{alignat}{2}%
m_{\,\,k}^{\mathbb{R}^{n+q}}(V\times{}0^q)&=0, & & 0\leq k<q;
\label{SiMua1} \\ %
m_{\,{k+q}}^{\mathbb{R}^{n+q}}(V\times{}0^q)&=
m_{k}^{\mathbb{R}^{n}}(V)\,\gamma_{\,k}^{(q)},& \ \ &0\leq k\leq
n\,.\label{SiMua2}
\end{alignat}%
\end{subequations}
\end{remark}
\begin{remark} %
\label{CompPol}%
According to \textup{Theorem \ref{IRPp}},
the transformation which maps the polynomial \(M_{\,V}^{{\mathbb{R}^{n}}}\)
into the polynomial \(M_{\,V\times{}0^q}^{{\mathbb{R}^{n+q}}}\) is essentially of the form
\begin{equation}
\label{MTP} %
\sum\limits_{0\leq k\leq n}m_kt^k\to \sum\limits_{0\leq k\leq
n}\gamma_km_kt^k,
\end{equation}
where \(\gamma_k\) is a certain sequence of \textsf{multipliers}.
\textup{(}The factor \(t^q\) before the sum in \eqref{MPn2}  is
not essential\textup{)}. The transformations of the form
\textup{\eqref{MTP}} were already discussed is Section \ref{LPEF}.
There such transformations were considered in relation with
location of roots of polynomials and entire functions.
\end{remark}%
\begin{lemma}%
\label{NonMult}%
 For any \(q\), \(q=1,\,2,\,3,\,\ldots\), the
sequence \(\{\gamma_{\,k}^{(q)}\}_{k=0,\,1,\,2,\,\ldots}\) is not
a multiplier sequence in the sense of \textup{Definition
\ref{PReZ}}.
\end{lemma}%
\textsf{PROOF.} In Section  \ref{LocRoot} we explain that the
entire function
\begin{equation}
\label{Multq}
\mu_q(t)=\displaystyle\sum\limits_{0\leq{}k<\infty}\frac{\gamma_{\,k}^{(q)}}{k!}t^k
\end{equation}
has infinitely many non-real roots.
 The entire function \(\mu_q(t)\),
\eqref{Multq}, appears as the function
\(\mathcal{M}_{B^n\times{}0^q}(t)\) in Section \ref{LPEF}. (Up to
a constant factor which is not essential for study the roots.)
 According to the Polya-Schur
Theorem, which was formulated in Section \ref{LPEF}, the sequence
\(\{\gamma_{\,k}^{(q)}\}_{k=0,\,1,\,2,\,\ldots}\) is not a
multiplier sequence.
\hfill\framebox[0.45em]{ }\\

\begin{remark}
In Section \ref{LPEF} we study the function \(\mu_q(t)\) in much
more details that it is needed to prove \textup{Lemma
\ref{NonMult}}. The study of section \ref{LPEF} is aimed to
clarify for which \(q\) the roots of the function in question are
located in the left half plane. The question whether there are
non-real roots is much more rough. This question may be answered
from very general considerations.
 The function \(\mu_q\) admits the
integral representation:
\begin{equation}%
\label{IRMuF}%
\mu_q(t) =
q\omega_q\int\limits_{0}^{1}(1-\xi^2)^{\frac{q}{2}-1}\xi\,{}e^{\xi{}t}\,d\xi\,.
\end{equation}
(Expanding the exponential \(e^{\xi{}t}\) into the Taylor series,
we see that the Taylor coefficients of the function in the right
hand side of \eqref{IRMuF} are the numbers
\(\frac{\gamma_{\,k}^{(q)}}{k!}\).) From \eqref{IRMuF} it follows
that the function \(\mu_q(t)\) is an entire function of
exponential type, and that its indicator diagram is the interval
\([0,\,1]\). Moreover,
\(\sup_{-\infty<t<\infty}|\mu_q(it)|<\infty\). In particular, the
function \(\mu_q(it)\) belongs to the class of entire functions
which is denoted by \(C\) in \cite{Lev2}, Lecture 17. From
\textup{ Theorem of Cartwright-Levinson (Theorem 1 of the Lecture
17 from \cite{Lev2})} it follows that the function \(\mu_q(t)\)
has infinitely many roots, these roots have  positive density, and
the `majority' of these roots is located `near' the rays \(\arg
t=\frac{\pi}{2}\) and \(\arg t=-\frac{\pi}{2}\). In particular,
the function \(\mu_q(t)\) has infinitely many non-real roots. (We
already used this reasoning proving \textup{Statement 2 of Theorem
\ref{BLPC}}.) \hfill\framebox[0.45em]{ }
\end{remark}

\noindent%
\textsf{PROOF of LEMMA \ref{IPRom}.} Let
\((x,s)\in{\mathbb{R}^{n+1}}\), where \(x\in\mathbb{R}^n\), and
\(s\in\mathbb{R}\). Then by Pythagorean theorem,
\[\textup{dist}_{\mathbb{R}^{n+1}}^2((x,s),V\times{0})=
\textup{dist}_{\mathbb{R}^n}^2(x,V)+s^2.\] Therefore, the
equivalence holds:
\begin{equation}
\label{eqd}%
\Big( \textup{dist}_{\mathbb{R}^{n+1}}((x,s),V\times{}0)\leq
t\Big) \Longleftrightarrow
\Big(\textup{dist}_{\mathbb{R}^{n}}(x,V)\leq\sqrt{t^2-s^2}\Big)
\end{equation}
 Let
\begin{equation}
\label{EV} \mathfrak{T}_{V\times{}0^1}^{\mathbb{R}^{n+1}}(t)=
\{(x,s)\in\mathbb{R}^{n+1}:\,\textup{dist}_{\mathbb{R}^{n+1}}((x,s),{V\times{}0^1})\leq
t\}.
\end{equation}
be the \(t\)-neighborhood of the set \(V\times{}0^1\) with respect
to the ambient space \(\mathbb{R}^{n+1}\). Thus,
\begin{equation}%
\label{Vnp1}%
\textup{Vol}_{n+1}(\mathfrak{T}_{V\times{}0^1}^{\mathbb{R}^{n+1}}(t))=
M_{V\times{}0^1}^{\mathbb{R}^{n+1}}(t).
\end{equation}%
For fixed \(s\in\mathbb{R}\), let \(\mathfrak{S}(s)\) be the
`horizontal section' of the set
\(\mathfrak{T}_{V\times{}0^1}^{\mathbb{R}^{n+1}}(t)\) on the
`vertical level' \(s\):
\[\mathfrak{S}(s)=
\{x\in\mathbb{R}^n:\,(x,s)\in\mathfrak{T}_{\,\,V}^{\mathbb{R}^{n+1}}(t)\}.\]
It is clear that
\begin{equation}%
\label{Kav}
\textup{Vol}_{n+1}(\mathfrak{T}_{V\times{}0^1}^{\mathbb{R}^{n+1}}(t))=
\int\textup{Vol}_n(\mathfrak{S}(s))ds\,.
\end{equation}%
 The equivalence \eqref{eqd} means that
\begin{equation*}%
\mathfrak{S}(s)=\mathfrak{T}_{V}^{\mathbb{R}^{n}}(\sqrt{t^2-s^2})=
\{x\in\mathbb{R}^n:\,\textup{dist}_{\mathbb{R}^{n}}(x,V)\leq{}\sqrt{t^2-s^2}\}.
\end{equation*}
Thus,
\begin{equation}%
\label{VoSe}%
 \textup{Vol}_n(\mathfrak{S}(s))=M_V^{\mathbb{R}^n}(\sqrt{t^2-s^2}).
\end{equation}%
From \eqref{Kav} and \eqref{VoSe} it follows that
\begin{equation*}%
\label{fin}
M_{V\times{}0^1}^{\mathbb{R}^{n+1}}(t)=\int\limits_{-t}^tM_V^{\mathbb{R}^{n}}(\sqrt{t^2-s^2})ds\,.
\end{equation*}%
Changing variable \(s\to{}ts^{1/2}\), we obtain
\begin{equation*}
M_{\,\,V}^{\mathbb{R}^{n+1}}(t)=
t\,\int\limits_{0}^1M_V^{\mathbb{R}^{n}}(t(1-s)^{1/2})s^{-1/2}ds.
\end{equation*}
Substituting the expression \eqref{IRP1} for
\(M_{V}^{\mathbb{R}^{n}}\) into the last formula, we obtain
\begin{equation*}
M_{V\times{}0^1}^{\mathbb{R}^{n+1}}(t)=t\sum\limits_{0\leq k\leq
n}m_k(V)\,t^k\int\limits_0^1(1-s)^{k/2}s^{-1/2}ds.
\end{equation*}
According to Euler,
\[\int\limits_0^1(1-s)^{k/2}s^{-1/2}ds=\textup{B}{\textstyle\big(\frac{k}{2}+1,\frac{1}{2}\big)}=
\frac{\Gamma\big(\frac{1}{2}\big)\,
\Gamma\big(\frac{k}{2}+1\big)}{\Gamma\big(\frac{k+1}{2}+1\big)}=
\pi^{1/2}\frac{
\Gamma\big(\frac{k}{2}+1\big)}{\Gamma\big(\frac{k+1}{2}+1\big)}\cdot\]
Thus, \eqref{MPn1om} holds.
 \hfill\framebox[0.45em]{ }

\noindent%
 \textsf{PROOF of THEOREM \ref{IRPp}.} For \(q=0\), the statement of the Theorem is
self-evident. Let us show how to pass from \(q\) to \(q+1\). Since
 \(V\times{}0^{q+1}=(V\times{}0^{q})\times{}0^1\), and
  \(\mathbb{R}^{n+q+1}=\mathbb{R}^{n+q}\times{}\mathbb{R}^{1}\),
 we can apply Lemma \ref{IPRom} to the convex set \(V\times{}0^{q}\)
 whose Minkowski polynomial is \eqref{MPn2} by the induction assumption.
 The induction assumption can be formulated as
\begin{subequations}
\label{SiMup}
\begin{alignat}{2}%
m_{\,\,k}^{\mathbb{R}^{n+q}}(V\times{}0^q)&=0, & & 0\leq k<q;
\label{SiMup1} \\ %
m_{\,{k}}^{\mathbb{R}^{n+q}}(V\times{}0^q)&=
m_{k-q}^{\mathbb{R}^{n}}(V)\,\gamma_{\,k-q}^{(q)},& \ \ &q\leq
k\leq q+n\,.\label{SiMup2}
\end{alignat}%
\end{subequations}
 By Lemma \ref{IPR},
 \begin{subequations}
\label{SiMuInd}
 \begin{align}
 \label{SiMuaInd1}
 m_{\,{k}}^{\mathbb{R}^{(n+q)+1}}((V\times{}0^q)\times{}0^1)&=0,\ \ k=0\,;\\
 \label{SiMuaInd2}
m_{\,{k}}^{\mathbb{R}^{(n+q)+1}}((V\times{}0^q)\times{}0^1)&=
m_{\,{k-1}}^{\mathbb{R}^{(n+q)}}((V\times{}0^q))\cdot{}\gamma_{\,k-1}^{(1)}\,,\
\ 1\leq{}k\leq{}n+q+1\,.
 \end{align}
 \end{subequations}
 In view of the identity
 \[\gamma_{k}^{(q)}\cdot{}\gamma_{k+q}^{(1)}=\gamma_{k}^{(q+1)}\,,\]
 \eqref{SiMuInd} takes the form \eqref{SiMup} with \(q\) replaced
 by \(q+1\)\,.
 \hfill\framebox[0.45em]{ }
\begin{remark}
From \eqref{Vpdb} and \eqref{MuSe} it follows that
\begin{equation}
\label{IdFG}%
 \gamma_{\,k}^{(q)}=\frac{\omega_{k+q}}{\omega_{k}}\,.
\end{equation}
Thus, the equalities \eqref{SiMua2} can be rewritten as
\begin{equation}
\label{IntrCo}%
\frac{\,m_{\,{k+q}}^{\mathbb{R}^{n+q}}(V\times{}0^q)}{\omega_{k+q}}=
\frac{\,m_{\,{k}}^{\mathbb{R}^{n}}(V)\,}{\omega_{k}},\quad
q=0,\,1,\,2,\,\ldots\,.
\end{equation}
The equality \eqref{IntrCo} holds for \(k=0,\,1,\,\ldots\,,\,n\).
For other \(k\), the value \(m_k^{\mathbb{R}^n}(V)\) is not yet
defined. Let us agree that
\begin{equation}
\label{Agre}%
\omega_k=1 \textup{\ \ for\ \ }k<0,\ \ m_k^{\mathbb{R}^n}(V)=0
\textup{\ \ for\ \ }k<0\textup{\ \ and for\ \ }k>n\,.
\end{equation}
Under this agreement, the equality \eqref{IntrCo} holds for
\emph{every} \(k\in\mathbb{Z}\)\textup{:} for \(k>n\) or for
\(k<-q\)  \eqref{IntrCo} is trivial, for \(-q\leq{}k\leq{}-1\) it
coincides with \eqref{SiMua1}, for \(0\leq{}k\leq{}n\) -- with
\eqref{SiMua2}.
\end{remark}
\paragraph{The Minkowski polynomials for the \(q\)-th adjoint to the ball \(B^n\).}
In particular, applying Theorem \ref{IRPp} to the case \(V=B^n\),
\(B^n\subset\mathbb{R}^n\), we obtain:
\begin{equation}
\label{ExpMinAdq}%
M_{B^n\times{}0^q}^{\mathbb{R}^{n+q}}(t)=
\omega_n\omega_q\,t^q\mathcal{M}_{B^n\times{}0^q}(nt)\,,
\end{equation}
where the normalized Minkowski polynomial
\(\mathcal{M}_{B^n\times{}0^q}\) is defined as
\begin{equation}
\label{NMPq}%
\mathcal{M}_{B^n\times{}0^q}(t)=
\sum\limits_{0\leq{}k\leq{}n}\frac{n!}{(n-k)!n^k}%
\frac{\Gamma{(\frac{q}{2}+1)}\Gamma(\frac{k}{2}+1)}{\Gamma(\frac{k+q}{2}+1)}\frac{t^k}{k!}\,.
\end{equation}
The polynomial \(\mathcal{M}_{B^n\times{}0^q}\) is the Jensen
polynomial associated with the entire functions
\(M^{B^n\times{}0^q}\):
\begin{equation}
\label{JPMPa}
\mathcal{M}_{B^n\times{}0^q}(t)=\mathscr{J}_n(M_{B^{\infty}\times{}0^q};t),
\end{equation}
where
\begin{equation}
\label{JPMEb}
\mathcal{M}_{B^{\infty}\times{}0^q}(t)=\sum\limits_{0\leq{}k<\infty}
\frac{{\Gamma(\frac{q}{2}+1)}\Gamma(\frac{k}{2}+1)}{\Gamma(\frac{k+q}{2}+1)}\frac{t^k}{k!}\,.
\end{equation}
Comparing with \eqref{IRMuF}, we obtain Comparing with
\eqref{IRMuF}, we obtain
\begin{equation}
\label{NMPq1}%
\mathcal{M}_{B^{\infty}\times{}0^q}(t)=
q\int\limits_{0}^{1}(1-\xi^2)^{\frac{q}{2}-1}\xi\,{}e^{\xi{}t}\,d\xi\,.
\end{equation}
For every \(q=0,\,1,\,2,\,\ldots\,\), the function
\(\mathcal{M}_{B^{\infty}\times{}0^q}\) is an entire function of
the exponential type one.
\begin{lemma}\ \
\label{HPEGM}
\begin{enumerate}
\item %
For \(q=0,\,1,\,2,\,4\), the entire function
\(\mathcal{M}_{B^{\infty}\times{}0^q}\) is in the Hurwitz class
\(\mathscr{H}\);
\item
For \(q\geq{}5\),  the entire function
\(\mathcal{M}_{B^{\infty}\times{}0^q}\) is not in the Hurwitz
class: it has infinitely many roots in the open right half plane
\(\{z:\,\,\textup{Im} z >0\}.\).
\end{enumerate}
\end{lemma}
Proof of this Lemma is presented in Section \ref{LocRoot}.
Statement 2 is a consequence of the asymptotic calculation of the
function \(\mathcal{M}_{B^{\infty}\times{}0^q}\). (See Lemma
\ref{NonHq}.)

 For \(q=0\), the
function \(\mathcal{M}_{B^{\infty}\times{}0^0}=e^t\), thus it is
of type \textup{I} in the Laguerre-Polya class:
\(\mathcal{M}_{B^{\infty}\times{}0^0}\in
\mathscr{L}\text{-}\mathscr{P}\text{-}\textup{I}\). For \(q=2\)
and \(q=4\) the functions \(\mathcal{M}_{B^{\infty}\times{}0^q}\)
can be calculated explicitly and investigated by elementary
methods. The case \(q=1\) is more involved. The case \(q=3\)
remains open.\hfill\framebox[0.45em]{ }

\noindent%
 \textsf{PROOF of Statement 1 of THEOREM \ref{NMR}.} Let \(q\geq{5}\) be given.
 According to statement 2 of Lemma \ref{HPEGM}, the function
\(\mathcal{M}_{B^{\infty}\times{}0^q}\) has infinitely many roots
in the open right half-plane. In view of the approximation
property of the Jensen polynomials (Lemma \ref{CJPL}), for
\(n\geq{}n(q)\) some  roots of the Jensen polynomial
\(\mathscr{J}_n(M_{B^{\infty}\times{}0^q};t)\) are located in the
open right half-plane. In view of \eqref{ExpMinAdq}  and
\eqref{JPMPa}, some roots of the Minkovski polynomial
\(M_{B^n\times{}0^q}\) of the (non-solid) convex set
\(B^n\times{}0^q\), \(B^n\times{}0^q\subset{}\mathbb{R}^{n+q}\),
are located in the open right half-plane. Fix \(n:\,n\geq{}n(q)\).
Consider the ellipsoids \(E_{n,\,q,\,\varepsilon}\) defined in
\eqref{Ell}, \(E_{n,\,q,\,\varepsilon}\subset\mathbb{R}^{n+q}\).
For \(\varepsilon>0\), the ellipsoid \(E_{n,\,q,\,\varepsilon}\)
is a solid convex set with respect to the ambient space
\(\mathbb{R}^{n+q}\). The family of the convex sets
\(\{E_{n,\,q,\,\varepsilon}\}_{\varepsilon>0}\) is monotonic,
 (See Remark \ref{AprC1} and footnote \ref{monot}), and
 \begin{equation}
 \label{LimEll}
 \lim_{\varepsilon\to+0}E_{n,\,q,\,\varepsilon}=B^n\times{}0^q\,.
 \end{equation}
It is known that the Minkowski polynomials \(M_V(t)\) depends on
the set \(V\) continuously: see Section \ref{ChPrMiPo} and
footnote \ref{topo}. Therefore,
\begin{equation}%
\label{LRLMP}%
\lim_{\varepsilon\to{}0}M_{E_{n,\,q,\,\varepsilon}}^{\mathbb{R}^{n+q}}(t)=
M_{B^n\times{}0^q}^{\mathbb{R}^{n+q}}(t)
\end{equation}%
locally uniformly in \(\mathbb{C}\). Hence, there exists
\(\varepsilon(q,n), \ \varepsilon(q,n)>0\) such that the Minkowski
polynomial \(M_{E_{n,\,q,\,\varepsilon}}^{\mathbb{R}^{n+q}}\) has
roots located in the open right half-plane. %
\mbox{\hspace*{0.1ex}} \hfill\framebox[0.45em]{ }

\paragraph{The Weyl polynomials for the surfaces of the adjoint convex sets.}
 Passing to define the so-called adjoint Weyl polynomials\
\(W_{\,V\times{}0^q}^{\,\,p}\), we do this  following
\textup{Definition \ref{DeWPDC}} as a sample.
\begin{definition}
\label{DaWP}
Given a convex compact set \(V\), \(V\subset{}\mathbb{R}^n\), %
and a number \(q\),
 \(q= 0,\,1,\,2,\,3,\,\dots\,,\,\,\), the \(q\)-th adjoint
  Weyl polynomial \(W_{\,\partial(V\times{}0^q)}^{\,\,1}\) of the index \(1\) for the
convex surface \(\partial(V\times{}0^q)\) is defined by means of
the odd part of the \(q\)-th adjoint Minkowskii polynomial
\(M_{\,V\times{}0^q}^{\mathbb{R}^{n+q}}\):
\begin{equation}%
\label{fdAwp}%
2tW_{\partial(V\times{}0^q)}^{\,\,1}(t)
\stackrel{\textup{\tiny{}def}}{=}M_{\,V\times{}0^q}^{\mathbb{R}^{n+q}}(t)-
M_{\,V\times{}0^q}^{\mathbb{R}^{n+q}}(-t)\,,
\end{equation}%
where \(M_{\,V^{(q)}}^{\mathbb{R}^{n+q}}\) is the \(q\)-th adjoint
Minkowskii polynomial which was introduced in \textup{Definition
\ref{DeAdM}}. In more detail\,%
\footnote{According to the agreement \eqref{Agre},
\(m_{2l+1}^{\mathbb{R}^{n+q}}(V\times{}0^q)=0\) for \(2l+1<0\) or
\(2l+1>n+q\,\).}\,,%
\begin{equation}%
\label{fdAwp1}%
W_{\partial(V\times{}0^q)}^{\,\,1}(t)=
\sum\limits_{l\in\mathbb{Z}}m_{2l+1}^{\mathbb{R}^{n+q}}(V\times{}0^q)t^{2l}\,.
\end{equation}%
\end{definition}
From \eqref{fdAwp1} we may define the Weyl coefficients
\(k_{2l}(\partial{}(V\times{}0^q))\) according to Definition
\ref{DNWP}:
\begin{equation}%
\label{WCASu}%
k_{2l}(\partial(V\times{}0^q))=
m_{2l+1}^{\mathbb{R}^{n+q}}(V\times{}0^q)\frac{(2\pi)^l\omega_1}{\omega_{1+2l}}\,.
\end{equation}%
 Then we define the Weyl
polynomials \(W_{\partial{}(V\times{}0^q)}^p\) with higher \(p\)
according to\,%
\footnote{Remark that
\(\displaystyle{}\frac{2^{-l}\,\Gamma(\frac{p}{2}+1)}{\Gamma(\frac{p}{2}+l+1)}=
(2\pi)^{-l}\frac{\omega_{2l+p}}{\omega_{p}}\,.\)}
 Definition \ref{DGWP}:
 \begin{definition}
\begin{equation}%
\label{WPACSpq}%
W_{\partial{}(V\times{}0^q)}^p(t)\stackrel{ \textup{\tiny def}}{=}
\sum\limits_{l\in\mathbb{Z}}k_{2l}(\partial{}(V\times{}0^q))
(2\pi)^{-l}\frac{\omega_{2l+p}}{\omega_{p}}t^{2l}\,.
\end{equation}%
\end{definition}%
Thus,
\begin{equation}%
\label{WPACS}%
\omega_{p}t^pW_{\partial{}(V\times{}0^q)}^p(t)=
\sum\limits_{l\in\mathbb{Z}}\frac{m_{2l+1}^{\mathbb{R}^{n+q}}(V\times{}0^q)}{\omega_{2l+1}}\,
\omega_{1}\omega_{2l+p}\,t^{2l+p}\,.
\end{equation}%

Let us clarify how the Weyl polynomials for the convex surfaces
\(\partial{}V\) and \(\partial{}(V\times{}0^q)\) are related. Here
we also have to distinguish the cases even and odd \(q\).
\begin{lemma}%
\label{TINVQ}%
Let \(V,\,\,V\subset\mathbb{R}^n,\) be a solid compact convex set,
and let \(p>0,\,q>0\) be integers. Then
 \begin{enumerate}
 \item%
 For even \(q\)
 \begin{equation}
\label{CeQ}%
\omega_pt^p\cdot{}W_{\partial(V\times{}0^q)}^{\,p}(t)=
\omega_{p+q}t^{p+q}\cdot{}W_{\partial{}V}^{p+q}(t)\,;
\end{equation}
 \item%
 For odd \(q\)
\begin{equation}
\label{CoQ}
 \omega_pt^p\cdot{}W_{\partial(V\times{}0^q)}^{\,p}(t)=
\omega_{p+q-1}t^{p+q-1}W_{\hspace*{0.5ex}\partial{}(V\times{}0^1)}^{p+q-1}(t)\,.
\end{equation}
\end{enumerate}
\end{lemma}%
\textsf{PROOF of LEMMA \ref{TINVQ}.} We distinguish cases of
even and odd \(q\).\\
\hspace*{3.0ex}1. \(q\) is even. The equality \eqref{IntrCo} with
\(k=2l+1-q\) takes the form
\[\frac{\,m_{\,{2l+1}}^{\mathbb{R}^{n+q}}(V\times{}0^q)}{\omega_{2l+1}}=
\frac{\,m_{\,{2l+1-q}}^{\mathbb{R}^{n}}(V)\,}{\omega_{2l+1-q}}.\]
From this and \eqref{WPACS} it follows that
\begin{equation*}%
\omega_{p}t^pW_{\partial{}(V\times{}0^q)}^p(t)=
\sum\limits_{l\in\mathbb{Z}}\frac{m_{2l+1-q}^{\mathbb{R}^{n}}(V)}{\omega_{2l+1-q}}\,
\omega_{1}\omega_{2l+p}\,t^{2l+p}\,.
\end{equation*}%
Changing the summation variable: \(l\to{}l+\frac{q}{2}\), we
obtain
\begin{equation*}%
\omega_{p}t^pW_{\partial{}(V\times{}0^q)}^p(t)=
\sum\limits_{l\in\mathbb{Z}}\frac{m_{2l+1}^{\mathbb{R}^{n}}(V)}{\omega_{2l+1}}\,
\omega_{1}\omega_{2l+p+q}\,t^{2l+p+q}\,.
\end{equation*}%
The expression in the right hand side of the last equality has the
same structure that the expression in the right hand side of
\eqref{WPACS}, with \(V\times{}0^q\) replaced to \(V\), \(p\)
replaced by \(p+q\), \(q\) replaced by \(0\). So, \eqref{CeQ} is
proved.
\\
\hspace*{3.0ex}2. \(q\) is odd. The equality \eqref{IntrCo}
implies the equality
\[\frac{\,m_{\,{2l+1}}^{\mathbb{R}^{n+q}}(V\times{}0^q)}{\omega_{2l+1}}=
\frac{\,m_{\,{2l+1-(q-1)}}^{\mathbb{R}^{n+1}}(V\times{}0^1)\,}{\omega_{2l+1-(q-1)}}.\]
From this and \eqref{WPACS} it follows that
\begin{equation*}%
\omega_{p}t^pW_{\partial{}(V\times{}0^q)}^p(t)=
\sum\limits_{l\in\mathbb{Z}}\frac{m_{2l+1-
(q-1)}^{\mathbb{R}^{n+1}}(V\times{}0^1)}{\omega_{2l+1-(q-1)}}\,
\omega_{1}\omega_{2l+p}\,t^{2l+p}\,.
\end{equation*}%
Changing the summation variable: \(l\to{}l+\frac{q-1}{2}\), we
obtain
\begin{equation*}%
\omega_{p}t^pW_{\partial{}(V\times{}0^q)}^p(t)=
\sum\limits_{l\in\mathbb{Z}}\frac{m_{2l+1}^{\mathbb{R}^{n+1}}(V\times{}0^1)}{\omega_{2l+1}}\,
\omega_{1}\omega_{2l+p+q-1}\,t^{2l+p+q-1}\,.
\end{equation*}%
The expression in the right hand side of the last equality has the
same structure that the expression in the right hand side of
\eqref{WPACS}, with \(V\times{}0^q\) replaced to \(V\times{}0^1\),
\(p\) replaced by \(p+q-1\), \(q\) replaced by \(1\). So,
\eqref{CoQ} is proved. \hfill\framebox[0.45em]{ }

The meaning of Lemma \ref{TINVQ} lies in the following. Studying
the location of roots of Weyl polynomials related to convex
surfaces there is no need to consider the boundary surfaces
\(\partial{}(V\times{}0^q)\) of  \(q\)-th adjoint convex sets
\(V\times{}0^q\) for arbitrary large \(q\). It is enough to
restrict the consideration to the cases \(q=0\) and \(q=1\) only,
that is to the case of the set \(V\) itself and
to the case of the squeezed cylinder with the base \(V\).\\[1.0ex]

\noindent%
 \textsf{PROOF of Statement 2 of THEOREM \ref{NMR}.}
 By Statement 2 of Theorem \ref{str}, the entire function
 \(\mathcal{W}_{\partial{}(B^{\infty}\times{}0)}^{\,p+q-1}\) has
 infinitely many non-real roots which. (We have assumed that \(p+q-1\geq{}5.\)) If \(n\) is large enough,
the Jensen polynomial
\(\mathcal{W}^{p+q-1}_{\partial{}(B^{n+1}\times{}0)}(t){=}
\mathscr{J}_{2[n/2]}(\mathcal{W}^{p+q-1}_{\partial{}(B^{\infty}\times{}0)};t)\)
also has non-real roots. According to \eqref{JReWP2} and
\eqref{WRemOr2}, the Weyl polynomial
\(W^{p+q-1}_{\partial{}(B^{n}\times{}0)}(t)\) has roots which do
not belong to the imaginary axis. By Statement 2 of Lemma
\ref{TINVQ},
\[W_{\partial{}(B^{n}\times{}0^q)}^{\,p}
={\textstyle\frac{\omega_{p+q-1}}{\omega_p}t^{q-1}}W_{\partial{}(B^{n}\times{}0)}^{\,p+q-1}\,.\]
Thus, the Weyl polynomial
\(W_{\partial{}(B^{n}\times{}0^q)}^{\,p}\) has roots which do not
belong to the imaginary axis. For fixed \(\,q,\,n\) and a positive
\(\varepsilon\), consider the ellipsoid
\(E_{n,\,q,\,\varepsilon}\) defined by \eqref{Ell}. Since
\(E_{n,\,q,\,\varepsilon}\to{}B^{n}\times{}0^q\) as
\(\varepsilon\to+0\), also \(W_{E_{n,\,q,\,\varepsilon}}^p
\to{}W_{\partial{}(B^{n}\times{}0^q)}^{\,p}\) as
\(\varepsilon\to+0\). Hence, if \(\varepsilon\) is small enough:
\(0<\varepsilon\leq{}\varepsilon(n,p,q)\), the polynomial
\(W_{E_{n,\,q,\,\varepsilon}}\) has roots which do not belong to
the imaginary axis.
 \hfill\framebox[0.45em]{ }

\section{THE MINKOWSKI POLYNOMIAL OF THE CARTESIAN \newline
 PRODUCT OF CONVEX SETS.\label{MPCaPr}}
Let \(V_1\) and \(V_2\) be compact convex sets,
\[V_1\subset{}\mathbb{R}^{n_1},\ \  V_2\subset{}\mathbb{R}^{n_2}.\]
 Then the Cartesian product \(V_1\times{}V_2\)
is a compact convex set embedded into
the Cartesian product \(\mathbb{R}^{n_1}\times{}\mathbb{R}^{n_2}\).
Since \(\mathbb{R}^{n_1}\times{}\mathbb{R}^{n_2}\)
can be naturally identified with \(\mathbb{R}^{n_1+n_2}\),
we can consider \(V_1\times{}V_2\) as being embedded into \(\mathbb{R}^{n_1+n_2}\):
 \[V_1\times{}V_2\subset{}\mathbb{R}^{n_1+n_2}.\]
The natural question arises:\\ %
 \textit{How to express the Minkowski polynomial
\(M_{\,\,V_1\times{}V_2}^{\mathbb{R}^{n_1+n_2}}\) for the
Cartesian product \(V_1\times{}V_2\) in terms of the Minkowski
polynomials\footnote{The Minkowski polynomials
\(M_{V_1},\,M_{V_2},\,M_{V_1\times{} V_2}\) are considered with
respect to the ambient spaces
\(\mathbb{R}^{n_1},\,\mathbb{R}^{n_2},\,\mathbb{R}^{n_1+n_2}\)
respectively.}
\(M_{\,\,V_1}^{\mathbb{R}^{n_1}}\) and \(M_{\,\,V_2}^{\mathbb{R}^{n_2}}\) for the Cartesian factors \(V_1\)  and \(V_2\) ?} \\
To answer this question, we introduce a special multiplication in the set
of polynomials, the so-called \textit{M-multiplication}, which is suitable
for this goal. %
\begin{definition} %
\label{DMMul}%
The \textsf{\textsf{M-product} \(t^k\Mm{}t^l\)}
of two monomials \(t^k\) and \(t^l\)
is defined as
\begin{equation}%
\label{Mmm}%
t^k\Mm{}t^l\stackrel{\textup{\tiny{}def}}{=}
{\frac{\Gamma\big(\frac{k}{2}+1\big)\Gamma\big(\frac{l}{2}+1\big)}%
{\Gamma\big(\frac{k+l}{2}+1\big)}}\,t^{k+l}\,,\ \  k\geq{}0,\,l\geq{}0.
\end{equation}%
\textup{\small It is clear that}
\begin{equation}%
\label{PrMM}
\textup{a). }t^0\Mm{}t^k=t^k,\ \ \textup{b). } t^k\Mm{}t^l=t^l\Mm{}t^k,\ \
\textup{c). }(t^k\Mm{}t^l)\Mm{}t^m=t^k\Mm{}(t^l\Mm{t^m}).
\end{equation}%
{\small{}The M-multiplication \eqref{Mmm} of monomials can
 be extended to the multiplication
of polynomials by linearity:}
\begin{subequations}
\label{MMult}
\begin{gather}%
\label{MMult1}
\textup{For\ \ } A(t)=\sum\limits_{0\leq k\leq n_{{}_{1}}}a_kt^k,\ \ \
B(t)=\sum\limits_{0\leq l\leq n_{{}_{2}}}b_lt^l,\\
(A \,\Mm{}B)(t)=\sum\limits_{\substack{0\leq k\leq n_1,\\
0\leq l\leq n_2}}\hspace*{-1.0ex}a_kb_l(t^k\Mm{}t^l)=
{\sum\limits_{\substack{0\leq k\leq n_1,
0\leq l\leq n_2}}a_kb_l{\frac{\Gamma\big(\frac{k}{2}+1\big)\Gamma\big(\frac{l}{2}+1\big)}%
{\Gamma\big(\frac{k+l}{2}+1\big)}}\,t^{k+l}},
\notag
\end{gather}
{\small and finally,} the \textsf{M-product \(A\Mm{}B\)
of the polynomials} \(A\) and \(B\)
 is defined as
\begin{equation}%
\label{MMult2}
(A \,\Mm{}B)(t)=\sum\limits_{0\leq r\leq{}n_1+n_2}
\bigg(\sum\limits_{\substack{k\geq{}0,\,l\geq{}0,
k+l=r}}\hspace*{-1.5ex}
a_k\,b_l\,{\frac{\Gamma\big(\frac{k}{2}+1\big)\Gamma\big(\frac{l}{2}+1\big)}%
{\Gamma\big(\frac{k+l}{2}+1\big)}}\bigg)\,t^r\,.
\end{equation}%
\end{subequations}
\end{definition} %
\vspace{4.0ex}
\noindent
From (\ref{PrMM}.b) and (\ref{PrMM}.c) it follows that
\[A\Mm{}B=B\Mm{}A,\quad (A\Mm{}B)\Mm{}C=A\Mm(B\Mm{}C)
\]
for every polynomials \(A,\,B,\,C\). In particular, the `triple product'
\(A\Mm{}B\Mm{}C\) is well defined. This triple product can
be explicitly expressed in terms of the coefficients of the factors: if
\begin{equation*}
A(t)=\sum\limits_{0\leq k\leq n_1}a_kt^k,\ \
B(t)=\sum\limits_{0\leq l\leq n_2}b_lt^l,\ \
C(t)=\sum\limits_{0\leq m\leq n_3}c_mt^m,\ \
\end{equation*}
then
\begin{equation*}
(A\Mm{}B\Mm{}C)(t)=\hspace*{-3.5ex}
\sum\limits_{0\leq r\leq{}n_1+n_2+n_3}\hspace*{-1.0ex}
\bigg(\sum\limits_{\substack{k\geq{}0,\,l\geq{}0,\,m\geq{}0\\
k+l+m=r}}\hspace*{-2.5ex}
a_k\,b_l\,c_m{\frac{\Gamma\big(\frac{k}{2}+1\big)\Gamma\big(\frac{l}{2}+1\big)%
\Gamma\big(\frac{m}{2}+1\big)}%
{\Gamma\big(\frac{k+l+m}{2}+1\big)}}\bigg)\,t^r\,.
\end{equation*}
It is clear, that for every number \(\lambda\)
and for every polynomials \(A\) and \(B\),
\[(\lambda{}A)\Mm{}B=\lambda(A\Mm{}B).\]
Moreover, if
\[\mathbb{I}(t)\equiv 1,\ \ \mathbb{T}(t)\equiv{}t,\] then
\[\mathbb{I}\Mm{}A=A.\]
Thus, \textit{the polynomial \(\mathbb{I}\) is the unity
with respect to the \(M\)-Multiplication.}\\
It is worthy to mention that
\begin{equation}%
\label{Pon}%
t^{(\Mm{}k)}\stackrel{\textup{\tiny{}def}}{=}
\underbrace{t\Mm{}t\Mm{}\,\cdots\,\Mm{}t}_{k}=
\frac{(\sqrt{\pi}/2)^k}{\Gamma(\frac{k}{2}+1)}t^k.
\end{equation}%
\begin{remark} %
 The M-multiplication
by the polynomial \(\mathbb{T}\)
 is related to the transformation of the form \eqref{MTP}:
\begin{subequations}
\label{UInt}
\begin{align}
\textup{If }\hspace{22.0ex} A(t)&=\phantom{2^{-p}t^p}\sum\limits_{0\leq k\leq n} a_kt^k,\\
\textup{then }\hspace{5.5ex}%
(\underbrace{\mathbb{T}\Mm{}\,\cdots\,\Mm{}\mathbb{T}}_{p}\Mm{}A)(t)
&=2^{-p}t^p\sum\limits_{0\leq k\leq n} a_k\gamma_{\,\,k}^{(p)}t^k,
\end{align}
\end{subequations}
where the `multipliers' \(\gamma_{\,\,k}^{(p)}\) are defined by \eqref{MuSe}.\\
\end{remark}
\begin{lemma} %
\label{InReMP}
The M-product \(A\Mm{}B\) of polynomials \(A\) and \(B\) admits
the integral\footnote{The integrals in the right hand sides of \eqref{IRTP}
are Stieltjes integrals.} representation\footnote{At least, for \(t>0\).}
\begin{subequations}
\label{IRTP}
\begin{equation}
\label{IRTP1}
(A\Mm{}B)(t)=A(0)B(t)+\int\limits_{0}^{t}A\big((t^2-\tau^2)^{1/2}\big)\,dB(\tau)\,,
\end{equation}
as well as
\begin{equation}
\label{IRTP2}
(A\Mm{}B)(t)=A(t)B(0)+\int\limits_{0}^{t}B\big((t^2-\tau^2)^{1/2}\big)\,dA(\tau)\,.
\notag
\end{equation}
\end{subequations}
\end{lemma}
\noindent
\textsf{PROOF}. First of all, the expressions in the right hand sides
of \eqref{IRTP} are equal: Integrating by parts and replacing the variable
\(\tau\to(t^2-\tau^2)^{1/2}\), we obtain
\[A(0)B(t)+\int\limits_{0}^{t}A\big((t^2-\tau^2)^{1/2}\big)\,dB(\tau)=
A(t)B(0)+\int\limits_{0}^{t}B\big((t^2-\tau^2)^{1/2}\big)\,dA(\tau).\]
So, the expressions in the right hand sides
of \eqref{IRTP} which at the first glance are asymmetric with respect
 to \(A\) and \(B\) actually are symmetric. Let
 \begin{equation*}
 \label{ETC}
 A(t)=\sum\limits_{0\leq k\leq n_{{}_{1}}}a_kt^k,\ \ \
B(t)=\sum\limits_{0\leq l\leq n_{{}_{2}}}b_lt^l
\end{equation*}
be the expressions for the polynomials \(A\) and \(B\) in terms of their coefficients.
Let us substitute these polynomials into the right hand side
of \eqref{IRTP1}:
\begin{gather*}
A(0)B(t)+\int\limits_{0}^{t}A\big((t^2-\tau^2)^{1/2}\big)\,dB(\tau)=\\
a_0\sum\limits_{0\leq l\leq n_{{}_{2}}}b_lt^l+
\int\limits_0^t\bigg(\sum\limits_{0\leq k\leq n_{{}_{1}}}%
a_k(t^2-\tau^2)^{k/2}\bigg)\cdot %
\bigg(\sum\limits_{1\leq l\leq n_{{}_{2}}}l\,b_{l}\tau^{l-1}\bigg)\,d\tau=
\notag\\
\sum\limits_{0\leq l\leq n_{{}_{2}}}a_0b_lt^l+
\sum\limits_{\substack{0\leq k\leq{}n_1\\
1\leq{}l{}\leq{}n_2}}a_k{}b_l\cdot{}l\int\limits_0^t(t^2-\tau^2)^{k/2}\tau^{l-1}d\tau\,.
\notag
\end{gather*}
Changing variable \(\tau\to{}t\tau^{1/2}\), we get
\begin{multline*}%
l\int\limits_0^t(t^2-\tau^2)^{k/2}\tau^{l-1}d\tau
=t^{k+l}(l/2)\int\limits_0^1(1-\tau)^{k/2}\tau^{l/2-1}d\tau
=\\
t^{k+l}(l/2) B\Big(\frac{k}{2}+1;\frac{l}{2}\Big)\,.
\end{multline*}%
Now, according to Euler,
\[(l/2) B\Big(\frac{k}{2}+1;\frac{l}{2}\Big)=
\frac{\Gamma\big(\frac{k}{2}+1\big)\frac{l}{2}\Gamma\big(\frac{l}{2}\big)}%
{\Gamma\big(\frac{k+l}{2}+1\big)}=
\frac{\Gamma\big(\frac{k}{2}+1\big)\Gamma\big(\frac{l}{2}+1\big)}%
{\Gamma\big(\frac{k+l}{2}+1\big)}\,.\]
Thus, the right hand side of \eqref{IRTP1} can be transformed into the
right hand side of \eqref{MMult2}.
\hfill{}Q.E.D.\\
\begin{theorem}
\label{MPCarP}
Given the compact convex sets \(V_1\) and \(V_2\),
\(V_1\subset{}\mathbb{R}^{n_1},  V_2\subset{}\mathbb{R}^{n_2},\)
let \(M_{\,V_1}^{\mathbb{R}^{n_1}}(t),\,M_{\,V_2}^{\mathbb{R}^{n_2}}(t)\)
be the Minkowski polynomials for the sets
\(V_1\) and \(V_2\). Then the Minkowski polynomial
\(M_{\,\,V_1\times{}V_2}^{\mathbb{R}^{n_1+n_2}}\) of the Cartesian product \(V_1\times{}V_2\)
 is equal to the M-product
of  the polynomials \(M_{V_1}^{\mathbb{R}^{n_1}}\) and \(M_{V_2}^{\mathbb{R}^{n_12}}\):
\begin{equation}
\label{MPCP}%
M_{\,\,V_1\times{}V_2}^{\mathbb{R}^{n_1+n_2}}
=M_{\,V_1}^{\mathbb{R}^{n_1}}\Mm{}M_{\,V_2}^{\mathbb{R}^{n_2}}\,.
\end{equation}
\end{theorem}
A sketch of proof of this theorem can be found in \cite{Had},
Chapter VI, Section 6.1.9. A detailed proof is presented below.
\begin{remark} Let \(S\) be `the origin' of \(\mathbb{R}^1\), that is
the one-point set:\\
 \(S=\{t:t=0\}\).
Then \(M_S(t)=2t\), that is
\begin{equation}%
\label{MOPS}
M_S(t)=2\,\mathbb{T}(t).
\end{equation}%
Let \(V\) be a compact convex set embedded into \(R^{n}\).
 The Cartesian product \(V\times{}\underbrace{S\,\times\,\cdots\,\times\,S}_p\)
 can be identified with the
 convex set
\(V\times{}0^p,\,\, V\times{}0^p\subset\mathbb{R}^{n+p}\).
Thus,
 \[M_{V\times{}\underbrace{\scriptstyle{}S\,\times\,\cdots\,\times\,S}_p}(t)
 =M_{\,\,V\times{}0^p}^{\mathbb{R}^{n+p}}(t)\,,\]
 or
 \begin{equation}
 \label{OInt}%
 2^p(\underbrace{\mathbb{T}\Mm\,\cdots\,\Mm{}\mathbb{T}}_p)\Mm{}M_{\,V}^{\mathbb{R}^{n}}
 =M_{\,\,V\times{}0^p}^{\mathbb{R}^{n+p}}.
 \end{equation}
In view of \eqref{UInt} and \eqref{MOPS}, the equality \eqref{OInt}
 is another form of the equality \eqref{MPn2}.\\
\end{remark}
\noindent \textsf{PROOF of THEOREM \ref{MPCarP}.} Denote
\[V=V_1\times{}V_2\,.\]
 According to the
identification
\(\mathbb{R}^{n_1+n_2}=\mathbb{R}^{n_1}\times{}\mathbb{R}^{n_2}\),
we present a point \(x\in\mathbb{R}^{n_1+n_2}\) as a pair
\(x=(x_1,\,x_2)\), where \(x_1\in\mathbb{R}^{n_1},\
x_2\in\mathbb{R}^{n_2}\). It is clear that
\begin{equation}%
\label{PiThe}%
 \dist_{\mathbb{R}^{n_1+n_2}}^2(x,\,V)=
\dist_{\mathbb{R}^{n_1}}^2(x_1,\,V_1)+\dist_{\mathbb{R}^{n_2}}^2(x_2,\,V_2).
\end{equation}
For \(\tau>0,\,\tau^{\prime}>0,\,\tau^{\prime\prime}>0\), let
\(V(\tau), V_1(\tau^{\prime})\) and \(V_2(\tau^{\prime\prime})\)
be the \(\tau\)-neighborhood of \(V\) with respect to
\(\mathbb{R}^{n_1+n_2}\), the \(\tau^{\prime}\)-neighborhood of
\(V_1\) w.\,r.\,t. \(\mathbb{R}^{n_1}\) and
\(\tau^{\prime\prime}\)-neighborhood of \(V_2\) w.\,r.\,t.
\(\mathbb{R}^{n_2}\) respectively:
\begin{gather*}
V(\tau)=V+\tau{}B_{n_1+n_2},\ \
V_1(\tau^{\prime})=V_1+{\tau}^{\prime}B_{n_1},\ \ %
V_2(\tau^{\prime\prime})=V_2+\tau^{\prime\prime}B_{n_2}; \\
V,\,B_{n_1+n_2}\subset\mathbb{R}^{n_1+n_2};\quad{}
V_1,\,B_{n_1}\subset\mathbb{R}^{n_1};\quad{}
V_2,\,B_{n_2}\subset\mathbb{R}^{n_2}.
\end{gather*}

Here \(B^n\) be the Euclidean ball of the radius one in
\(\mathbb{R}^n\). (With \(n=n_1+n_2, n_1,\,n_2\) respectively.)

Given a number \(t,\,t>0\), consider the \(t\)-neighborhood
\(V(t)\) of \(V=V_{1}\times{}V_2\), and let
\[0=\tau_0<\tau_1<\,\dots\,<\tau_{N-1}<\tau_N=t\]
be a partition of the interval \([0,\,t]\). From \eqref{PiThe} it
follows that
\begin{align}%
\label{DoIn}%
\big(V_1(0)\times{}V_2(t)\big)%
\cup&\Big(\bigcup_{1\leq{}k\leq{}N}\big(V_1(\tau_k)%
\setminus{}V_1(\tau_{k-1})\big)\times{}V_2\big((t^2-\tau_{k}^2)^{1/2
}\big)\Big)
\notag\\
&\subseteq
V(t)\subseteq\\
\big(V_1(0)\times{}V_2(t)\big)%
\cup&\Big(\bigcup_{1\leq{}k\leq{}N}\big(V_1(\tau_k)
\setminus{}V_1(\tau_{k-1})\big)\times{}V_2\,%
\big((t^2-\tau_{k-1}^2)^{1/2}\big)\Big)\notag\,.
\end{align}
Since \(V_1(\tau_k)\supseteq{}V_1(\tau_{k-1})\),
\[\textup{Vol}_{n_1}\big(V_1(\tau_k)
\setminus{}V_1(\tau_{k-1}\big)=
{\textup{Vol}}_{n_1}\big(V_1(\tau_k)\big)-\textup{Vol}_{n_1}\big(V_1(\tau_{k-1})\big),\]
thus
\begin{multline*}
\textup{Vol}_{\,n_1+n_2}\Big(\big(V_1(\tau_k)
\setminus{}V_1(\tau_{k-1})\big)\times{}V_2\,%
\big((t^2-\tau_{\,l}^2)^{1/2}\big)\Big)=\\
\Big({\textup{Vol}}_{n_1}\big(V_1(\tau_k)\big)-
\textup{Vol}_{n_1}\big(V_1(\tau_{k-1})\big)\Big)\cdot
\textup{Vol}_{n_2}\Big(V_2\big((t^2-\tau_{\,l}^2)^{1/2}\big)\Big),\
\\
 l=k-1\textup{ or }l=k.
\end{multline*}
Moreover
\[\textup{Vol}_{\,n_1+n_2}\Big(V_1(0)\times{}V_2(t)\Big)=
\textup{Vol}_{\,n_1}\big(V_1(0)\big)\cdot{}\textup{Vol}_{\,n_2}\big(V_2(t)\big)\,.\]
In the notation of Minkowski polynomials, the last equalities take the form
\begin{subequations}
\label{PrFo}
\begin{gather}
\textup{Vol}_{\,n_1+n_2}\Big(\big(V_1(\tau_k)
\setminus{}V_1(\tau_{k-1})\big)\times{}%
V_2\,\big((t^2-\tau_{\,l}^2)^{1/2}\big)\Big)=\notag{}\\
\Big(M_{V_1}(\tau_k)-M_{V_1}(\tau_{k-1})\Big)\cdot{}M_{V_2}\big((t^2-\tau_{\,l}^2)^{1/2}\big),
\quad{}l=k-1\textup{ or }l=k,\\
\textup{Vol}_{\,n_1+n_2}\Big(V_1(0)\times{}V_2(t)\Big)=M_{V_1}(0)\cdot{}M_{V_2}(t)\,,
\end{gather}
and also
\begin{equation}
\textup{Vol}_{\,n_1+n_2}\big(V(t)\big)=M_{V}(t).
\end{equation}
\end{subequations}
Since the sets \(V_1(\tau_k)\setminus{}V_1(\tau_{k-1})\) for different
\(k\) do not intersect, and none of these sets intersects with the set \(V_1(0)\),
it follows from \eqref{DoIn} and \eqref{PrFo}
that
\begin{gather}
M_{V_1}(0)\cdot{}M_{V_2}(t)+
\sum\limits_{1\leq{}k\leq{}N}\Big(M_{V_1}(\tau_k)-M_{V_1}(\tau_{k-1})\Big)%
\cdot{}M_{V_2}\big((t^2-\tau_{k}^2)^{1/2}\big)\notag{}\\
\label{FBoS}%
\leq{}M_{V}(t)\leq{}\\
M_{V_1}(0)\cdot{}M_{V_2}(t)+
\sum\limits_{1\leq{}k\leq{}N}\Big(M_{V_1}(\tau_k)-M_{V_1}(\tau_{k-1})\Big)%
\cdot{}M_{V_2}\big((t^2-\tau_{k-1}^2)^{1/2}\big)\notag{}\,.
\end{gather}
Passing to the limit as \(\max{(\tau_{k}-\tau_{k-1})}\to{}0\) in
the last inequality, we express the Minkowski polynomial
\(M_V(t)\) as the Stieltjes integral
\begin{equation}
\label{IRepr}
M_V(t)=M_{V_1}(0)\cdot{}M_{V_2}(t)+%
\int\limits_{0}^tM_{V_2}\big((t^2-\tau^2)^{1/2}\big)\,dM_{V_1}(\tau)\,.
\end{equation}
According to Lemma \ref{InReMP}, the expression in the
right hand side of \eqref{IRepr}
is equal to
\(\big(M_{V_1}\Mm{}M_{V_2}\big)(t)\).\hfill{}Q.E.D.\\[1.0ex]

\section{
 PROPERTIES OF ENTIRE FUNCTIONS GENERATING\\
  THE MINKOWSKI AND WEYL
POLYNOMIALS\\ FOR THE DEGENERATED CONVEX SETS
\(\boldsymbol{B^{n+1}\times{}0^q}\). \label{LocRoot}}
In this section we investigate location of roots of the entire
functions generating the Minkowski and Weyl polynomials related to
the `degenerated' convex sets \(B^{n+1}\times{}0^q\). These are:
\begin{itemize}
\item
The entire functions \(\mathcal{M}_{B^n\times{}0^q}\) which
appears in \eqref{JPMEb}, in particular for \(q=1\) in
\eqref{LiEFMPS2} .
\item
The entire function
\(\mathcal{W}_{\partial{}(B^{\infty}\times{}0)}^{\,p},\ \
1\leq{}p<\infty\), which appears in \eqref{ReWP3}\,.
\item
The entire function
\(\mathcal{W}_{\partial{}(B^{\infty}\times{}0)}^{\,\infty}\) which
appears in \eqref{ReWP4}\,.
\end{itemize}
The functions \(\mathcal{M}_{B^n\times{}0^q}\),
\(\mathcal{W}_{\partial{}(B^{\infty}\times{}0)}^{\,p},\ \
1\leq{}p<\infty\) can not be calculated explicitly (except a very
special values of the parameters \(p\) or \(q\)), but they can be
calculated asymptotically.

 The above mentioned functions admit integrable representations:
\begin{equation}%
\label{MInt}%
\mathcal{M}_{B^n\times{}0^q}(t)=q\int\limits_{0}^{1}(1-\xi^2)^{\frac{q}{2}-1}\xi\,{}e^{\xi{}t}\,d\xi\,;
\end{equation}%
\begin{equation}%
\label{WInt}%
\mathcal{W}_{\partial{}(B^{\infty}\times{}0)}^{\,p}(t)= %
p\int\limits_{0}^{1}(1-\xi^2)^{\frac{p}{2}-1}\xi\cos{}t\xi\,d\xi\,;
\end{equation}%
These integral representation can be used for the asymptotic
calculation of the functions \(\mathcal{M}_{B^n\times{}0^q}\),
\(\mathcal{W}_{\partial{}(B^{\infty}\times{}0)}^{\,p}\).

Another way to to calculate the functions
\eqref{MInt},\,\eqref{WInt} asymptotically is to use the structure
of their Taylor series. The Taylor coefficients of each of these
functions are ratios of factorials: these functions belong to the
so-called \textit{Fox-Wright function}, \cite{CrCs4}.

\textsf{The Fox-Wright function} is defined as
\begin{equation}
\label{FoWr}%
 {\sideset{_p}{_q}\Psii}%
\left\{\genfrac{}{}{0pt}{1}%
{\alpha_1\,\,\alpha_2\,\,\,\,\,\alpha_p}%
{\beta_1\,\,\beta_2\,\,\,\,\,\beta_p};%
\genfrac{}{}{0pt}{1}%
{\rho_1\,\,\rho_2\,\,\,\,\,\rho_q}%
{\sigma_1\,\,\sigma_2\,\,\,\,\,\sigma_q};%
z\right\}%
 \stackrel{\textup{\tiny
def}}{=}\sum\limits_{0\leq{}k<\infty}%
\frac{\prod_{j=1}^{p}\Gamma{}(\alpha_{j}\,k+\beta_j)}%
{\prod_{j=1}^{q}\Gamma{}(\rho_{j}\,k+\sigma_j)}%
\cdot\frac{x^k}{k!}\,.
\end{equation}

Comparing \eqref{FoWr} with the Taylor expansions
\eqref{JPMEb},\,\eqref{ReWP3},\,\eqref{ReWP4},\,we see that
\begin{equation}%
\label{FRMq}%
 \mathcal{M}_{B^n\times{}0^q}(t)=\Gamma\left(\frac{q}{2}+1\right)\cdot%
 {\sideset{_1}{_1}\Psii}%
\left\{\genfrac{}{}{0pt}{1}%
{\frac{1}{2}}%
{1};%
\genfrac{}{}{0pt}{1}%
{\frac{1}{2}}%
{1+\frac{q}{2}};%
t\right\}\,,\ \ 1\leq{}q<\infty\,,%
\end{equation}%
\begin{equation}
\label{As1to}
\mathcal{W}^{\,p}_{\partial{}(B^{\infty}\times{}0^{1})}(t)=
\Gamma({\textstyle{}\frac{1}{2}})\cdot\Gamma({\textstyle{}\frac{p}{2}+1})\cdot%
{\sideset{_1}{_2}\Psii}%
\left\{\genfrac{}{}{0pt}{1}%
{1}%
{1};%
\genfrac{}{}{0pt}{1}%
{1\,\,\,\,\,\,\,\,1\,\,\,\,\,}%
{\frac{1}{2}\,\,\frac{p}{2}+1};%
-\frac{t^2}{4}\right\}\,,\ \ 1\leq{}p<\infty\,,%
\end{equation}
Asymptotic behavior of the functions
\({\sideset{_p}{_q}\Psii}(z)\) has been studied by E.\,Barnes,
\cite{Bar}, G.N.\,Watson, \cite{Wat}, G.\,Fox, \cite{Fox}),
E.M.\,Wright, \cite{Wr1}, \cite{Wr2}.

\paragraph{Analysis of the function
\(\mathcal{M}_{B^n\times{}0^q}(t)\):} \ We would like to
investigate for which \(q\) this function belongs to the Hurwitz
class \(\mathscr{H}\). According to \eqref{FRMq}, we may readdress
the question to the proportional function \({\sideset{_1}{_1}\Psii}%
\left\{\genfrac{}{}{0pt}{1}%
{\frac{1}{2}}%
{1};%
\genfrac{}{}{0pt}{1}%
{\frac{1}{2}}%
{1+\frac{q}{2}};%
z\right\}\). From the Taylor expansion it is clear that this
function is an entire function of exponential type. Since the
Taylor coefficients of the function are positive, its defect\,%
\footnote{We recall that the defect of an entire function \(H\) of
exponential type is defined by \eqref{DefH}.} %
 is non-negative. So, \textit{the function \(\mathcal{M}_{B^n\times{}0^q}(t)\)
 is in the Hurwitz class \(\mathscr{H}\) is and only if all roots
 of the function %
 \({\sideset{_1}{_1}\Psii}%
\left\{\genfrac{}{}{0pt}{1}%
{\frac{1}{2}}%
{1};%
\genfrac{}{}{0pt}{1}%
{\frac{1}{2}}%
{1+\frac{q}{2}};%
z\right\}\) %
are situated in the open left half plane.} To investigate the
location of roots of the last function, we use the following
asymptotic approximation, which can be derived from the results
stated in \cite{Wr1}, \cite{Wr2}:\\
For any \(\varepsilon: \ 0<\varepsilon<\dfrac{\pi}{2}\),\\[-3.0ex]
\begin{multline}
\label{Asymp1}
{\sideset{_1}{_1}\Psii}%
\left\{\genfrac{}{}{0pt}{1}%
{\frac{1}{2}}%
{1};%
\genfrac{}{}{0pt}{1}%
{\frac{1}{2}}%
{1+\frac{q}{2}};%
z\right\}=\\
=\begin{cases}
2^{\frac{q}{2}}z^{-\frac{q}{2}}e^z\left(1+r_1(z)\right)\,,&|\arg{}z|\leq{}\frac{\pi}{2}-\varepsilon\,;\\[1.0ex]
\frac{2}{\Gamma(\frac{q}{2})}z^{-2}+r_2(z),\,&|\arg{}z-\pi|\leq{}\frac{\pi}{2}-\varepsilon\,;\\[1.0ex]
2^{\frac{q}{2}}z^{-\frac{q}{2}}e^z+
\frac{2}{\Gamma(\frac{q}{2})}z^{-2}+r_3(z)\,,&|\arg{}z\mp\frac{\pi}{2}|\leq{}\varepsilon\,.
\end{cases}
\end{multline}
The reminders admit the estimates:
\begin{multline}%
\label{Rema1}%
 |r_1(z)|\leq{}C_1(\varepsilon)|z|^{-1},
|\arg{}z|\leq\frac{\pi}{2}-\varepsilon\,;\\
 |r_2(z)|\leq{}C_2(\varepsilon)|z|^{-2},
|\arg{}z-\pi|\leq\frac{\pi}{2}-\varepsilon\,;\hspace*{15.0ex}\\
|r_3(z)|\leq{}C_3\left(|z|^{-2}+|e^z||z|^{-(\frac{q}{2}+1)}\right)\,,\
\ |\arg{}z\mp\frac{\pi}{2}|\leq{}\varepsilon\,,
\end{multline}
where the values \(C_1(\varepsilon)<+\infty, \
C_2(\varepsilon)<+\infty, \ C_3(\varepsilon)<+\infty\) do not
depend on \(z\).

 From \eqref{Asymp1}, \eqref{Rema1} it follows
that for any \(\varepsilon>0\) the function %
\({\sideset{_1}{_1}\Psii}%
\left\{\genfrac{}{}{0pt}{1}%
{\frac{1}{2}}%
{1};%
\genfrac{}{}{0pt}{1}%
{\frac{1}{2}}%
{1+\frac{q}{2}};%
z\right\}\) %
has not more that finitely many roots outside the
angular domain
\(\{z:\,|\arg{}z\mp{}\frac{\pi}{2}|\leq{}\varepsilon\). Inside
this domain the analyzed function %
\({\sideset{_1}{_1}\Psii}%
\left\{\genfrac{}{}{0pt}{1}%
{\frac{1}{2}}%
{1};%
\genfrac{}{}{0pt}{1}%
{\frac{1}{2}}%
{1+\frac{q}{2}};%
z\right\}\) %
has infinitely many roots, and these roots are asymptotically
close to the roots of the approximating function \(f_q(z)\):
\[f_q(z)=2^{\frac{q}{2}}z^{-\frac{q}{2}}\left(e^z+
\frac{2^{1-\frac{q}{2}}}{\Gamma(\frac{q}{2})}z^{\frac{q}{2}-2}\right).\]
Investigating the location of roots of the approximating function \(f_q(z)\),
one should distinguish several cases: \\ %
\(\boldsymbol{q=4}\). In this case, the equation \(f_q(z)=0\) is
the equation \(e^z+\frac{1}{2}=0\), so, the roots of the
approximating function can be found explicitly: these root form an
arithmetical progression located on the straight line
\(\{z=x+iy:\,x=-\ln{}2,\,-\infty<y<\infty\}\). From this and
\eqref{Asymp1}-\eqref{Rema1} it follows that the roots of the
the analyzed function %
\({\sideset{_1}{_1}\Psii}%
\left\{\genfrac{}{}{0pt}{1}%
{\frac{1}{2}}%
{1};%
\genfrac{}{}{0pt}{1}%
{\frac{1}{2}}%
{1+\frac{q}{2}};%
z\right\}\) %
function are asymptotically close to the above appeared straight
line. Thus, \textit{for \(q=4\) all roots of the function
\({\sideset{_1}{_1}\Psii}%
\left\{\genfrac{}{}{0pt}{1}%
{\frac{1}{2}}%
{1};%
\genfrac{}{}{0pt}{1}%
{\frac{1}{2}}%
{1+\frac{q}{2}};%
z\right\}\) but finitely many are disposed in the open left half
plane.} Actually, for \(q=4\) \textsf{\textit{all}} roots of this
function are disposed in the open left half plane. To establish
this, one need the further analysis. This will be done a little
bit later.\\
\(\boldsymbol{q\not=4}\). In this case, the equation \(f_q(z)=0\)
is the equation %
\[e^z+
c_q\,z^{\frac{q}{2}-2}=0\,,\ \ \
c_q={\frac{2^{1-\frac{q}{2}}}{\Gamma(\frac{q}{2})}}\,,\] where the
exponent \(\frac{q}{2}-2\) is different from zero. The last
equation has infinitely many roots which have no finite
accumulation points and which are asymptotically close to the
`logarithmic parabola'
\begin{equation}
\label{LogPar}%
 x=({\textstyle\frac{q}{2}}-1)\ln{}(|y|+1)+ln|c|,\ -\infty<y<\infty\,,\ \
 (z=x+iy).
\end{equation}
From this and from \eqref{Asymp1}-\eqref{Rema1} it follows that
the roots of the
function \({\sideset{_1}{_1}\Psii}%
\left\{\genfrac{}{}{0pt}{1}%
{\frac{1}{2}}%
{1};%
\genfrac{}{}{0pt}{1}%
{\frac{1}{2}}%
{1+\frac{q}{2}};%
z\right\}\)  are asymptotically close to the logarithmic parabola
\eqref{LogPar}. Now we should distinguish the cases \(q<4\) and
\(q>4\).\\
\(\boldsymbol{q<4}\). In this case, the logarithmic parabola,
\eqref{LogPar} except may by its compact subset, is located inside
the left half plane. Since the roots of the analyzed function are
asymptotically close to this parabola, all roots but finitely many
are located in
the left half plane. \\
\(\boldsymbol{q>4}\). In this case, the logarithmic parabola,
\eqref{LogPar} except may by its compact subset, is located inside
the right half plane. So, all roots of the function
\({\sideset{_1}{_1}\Psii}%
\left\{\genfrac{}{}{0pt}{1}%
{\frac{1}{2}}%
{1};%
\genfrac{}{}{0pt}{1}%
{\frac{1}{2}}%
{1+\frac{q}{2}};%
z\right\}\), except finitely many, are located in the right half
plane.

Let us formulate this result as
\begin{lemma}
\label{NonHq} If \(q>4\), then the entire function
\(\mathcal{M}_{B^n\times{}0^q}\) has infinitely many roots within
the right half plane. In particular, this function does not belong
to the Hurwitz class \(\mathscr{H}\).
\end{lemma}
Claim 2 of Lemma \ref{HPEGM} is a consequence of Lemma
\ref{NonHq}.
\begin{lemma}%
\label{HurCl}%
For \(q:\,0\leq{}q\leq{}2\), the function
\(\mathcal{M}_{B^n\times{}0^q}\) belongs to the Hurwitz class
\(\mathscr{H}\).
\end{lemma}
\textsf{PROOF}. For \(q=0\), the assertion is evident: the
function in question is equal to \(e^z\). To investigate the case
\(q>0\), we use the integral representation
\begin{equation}
\label{IRPs1}
{\sideset{_1}{_1}\Psii}%
\left\{\genfrac{}{}{0pt}{1}%
{\frac{1}{2}}%
{1};%
\genfrac{}{}{0pt}{1}%
{\frac{1}{2}}%
{1+\frac{q}{2}};%
z\right\}=\frac{q}{\Gamma(\frac{q}{2}+1)}\,I_q(z)\,,
\end{equation}
where
\begin{equation}
\label{IRPs2}
I_q(z)=\int\limits_{0}^{1}(1-\xi^2)^{\frac{q}{2}-1}\xi\,{}e^{\xi{}t}\,d\xi\,.
\end{equation}
The defect of the entire function \({\sideset{_1}{_1}\Psii}%
\left\{\genfrac{}{}{0pt}{1}%
{\frac{1}{2}}%
{1};%
\genfrac{}{}{0pt}{1}%
{\frac{1}{2}}%
{1+\frac{q}{2}};%
z\right\}\) is non-negative. So it is enough to prove that this
function has no roots in the closed right half plane. The function
\(I_q(z)\) is of the form
\begin{equation}
\label{IRPs3}
I_q(z)=\int\limits_{0}^{1}\varphi_q(\xi)e^{\xi{}z}d\xi\,,
\end{equation}
where
\begin{equation}
\label{IRPs4} \varphi_q(\xi)=(1-\xi^2)^{\frac{q}{2}-1}\xi,\ \
0\leq{}\xi\leq{}1\,.
\end{equation}
The crucial circumstance is: \\
\textit{For \(q:\ 0\leq{}q\leq{}2\), the function
\(\varphi_q(\xi)\) is
positive and strictly increasing on the interval \((0,1).\)}\\
\begin{lemma}
\label{NZRHP}
\textsf{\textup{[Polya]}} %
\label{MoImHu}%
If \(\varphi(\xi)\) is a non-negative increasing function on the
interval \([0,1]\), then the entire function
\begin{equation}
\label{IRPs5}%
 I(z)=\int\limits_0^1\varphi(\xi)e^{\xi{}z}d\xi
\end{equation}
has no zeros in the closed right half plane.
\end{lemma}%
 This lemma is a continual analog of one Theorem of one theorem of
S.\,Kakeya. Proof of this Lemma and the reference to the paper of
S.\,Kakeya could be found in \cite{Po1}, \S \,1. We give another
proof. We have learned the idea of this proof  from \cite{OsPe}.
(See
Lemma 4 there.)\\[1.0ex]
\textsf{PROOF of Lemma \ref{NZRHP}.} Let \(z=x+iy\). Since
\(f(x)>0\) for  \(x\geq{}0\), \(f(x)\) has no zeros for
\(0\leq{}x<\infty\). Let us show that \(f(z)\) has no zeros for
\(0\leq{}x<\infty,\,y\not=0\). It is enough to consider the case
\(y>0\) only. We prove that \(\textup{Im}\,(e^{-z}f(z))<0\) for
\(z=x+iy,\,\,x\geq{}0,\,y>0\), thus \(f(z)\not=0\) for
\(z=x+iy,\,\,x\geq{}0,\,y>0\). To prove this, we use the integral
representation
\begin{equation}
\label{NIPa}
e^{-z}f(z)=\int\limits_0^{\infty}\psi(\xi)e^{-i\xi{}y}d\xi\,,
\end{equation}
where
\begin{equation}
\label{IRPoIP} \psi(\xi)=
\begin{cases}
\varphi(1-\xi)e^{-x\xi}, & 0\leq{}\xi\leq{}1\,,\\
0,&1<\xi<\infty\,.
\end{cases}
\end{equation}
In particular,
\begin{equation}
\label{IRPoIP1}
-(\textup{Im}\,e^{-z}f(z))=\int\limits_0^{\infty}\psi(\xi)\sin{}\xi{}y{}\,d\xi\,,
\ \ z=x+iy,\ \ x\geq{ }0,\,y>0\,,
\end{equation}
where the function \(\psi(\xi)\) is non-negative and decreasing on
\([0,\infty)\),  strictly decreasing on some non-empty open
interval, and \(\psi(\infty)=0\). Further,
\begin{multline}
\label{Fur}
\int\limits_0^{\infty}\psi(\xi)\sin{}\xi{}y{}\,d\xi=\sum\limits_{k=0}^{\infty}\,
\int\limits_{\frac{k\pi}{y}}^{\frac{(k+1)\pi}{y}}\psi(\xi)\sin{}\xi{}y{}\,d\xi=\\
\int\limits_{0}^{\frac{\pi}{y}}\left(\sum\limits_{k=0}^{\infty}
(-1)^k\psi(\xi+{\textstyle\frac{k\pi}{y}})\right)\sin{}\xi{}y{}\,d\xi>0\,:
\end{multline}
The series under the last integral is a Leibnitz type series. Thus
the sum of this series is non-negative on the interval of
integration, and is strictly positive on some subinterval.
\hfill\framebox[0.45em]{ }
\begin{lemma}%
\label{QEqFour}%
For \(q=4\), the function \(\mathcal{M}_{B^n\times{}0^q}\)
 belongs to the Hurwitz class \(\mathscr{H}\).
\end{lemma}%
\textsf{PROOF.} For \(q=4\), the integral in \eqref{IRPs2} can be
calculated explicitly:
\begin{equation}
\label{ExpCal4}%
 I_4(z)=\frac{(2z^2-6z+6)e^z+(z^2-6)}{z^4}\,\cdot
\end{equation}
Our goal is to prove that the function
\(\frac{(2z^2-6z+6)e^z+(z^2-6)}{z^4}\) has no roots in the closed
right half plane. Instead to investigate this function, we will
investigate the function
\begin{equation}
\label{AuxFunc1}%
 f(z)=(2z^2-6z+6)+(z^2-6)e^{-z}\,.
\end{equation}
We prove that the function \(f(z)\) has no other roots in the
closed right half plane than the root at the point \(z=0\) of
multiplicity four. The function \(f\) is of the form
\begin{equation}
\label{AuxFunc2}%
 f(z)=g(z)+h(z)\,, \ \textup{where} \ \ g(z)=(2z^2-6z+6),\
 h(z)=(z^2-6)e^{-z}\,.
\end{equation}
In the right half plane the function \(h\) is subordinate to the
function \(g\) in the following sense. For \(R>0\), let us
consider the contour \(\Gamma_R\) which consists of the interval
\(I_R\) of the imaginary axis and of the semicircle \(C_R\)
located in the right half plane:
\begin{equation}
\label{AuxFunc3}%
 \Gamma_R=I_R\cup{}C_R\,, \textup{ where }I_R=[-iR,iR], \
 C_R=\{z:\,|z|=R,\,\textup{Re}\,z\geq{}0\,.
\end{equation}
It is clear that \(|g(z)|\geq{}1.75\,|z|^2,\ \
|h(z)|\leq{}1.25\,|z|^2 \)  if  \(\ z\in{}C_R\) and \(R\) is large
enough. In particular, \(|g(z)|>|h(z)|\)  if \(z\in{}C_R\) and
\(R\) is large enough. On the imaginary axis,
\begin{equation}
\label{BDiCa}%
 |g(iy)|^2=36+12y^2+4y^4,\ \
|h(iy)|^2=36+12y^2+y^4\,,\ -\infty<y<\infty,
\end{equation}
 In particular,
\(|g(z)|\geq{}|h(z)|\) for \(z\in{}I_R\), and the inequality is
strict for \(z\not=0\). Thus,
\begin{equation}
\label{AuxFunc4}%
|g(z)|\geq|h(z)|\ \ \textup{ for } \ \ z\in\Gamma_R \ \ \textup{
and \(R\) is large enough.}
\end{equation}
For \(0<\varepsilon<1\), consider the function
\begin{equation}
\label{AuxFunc5}%
f_{\varepsilon}(z)=g(z)+(1-\varepsilon)h(z)\,.
\end{equation}
The function \(g\), which is a polynomial, has two simple roots:
\(z_{1,2}=\frac{3\pm{}i\sqrt{3}}{2}\). They are located in the
open right half plane. In view of \eqref{AuxFunc5} and Rouche's
theorem, for \(\varepsilon>0\) the function \(f_{\varepsilon}(z)\)
has precisely two roots \(z_1(\varepsilon),\,z_2(\varepsilon)\) in
the open right half plane. For \(\varepsilon\) positive an very
small, the roots \(z_1(\varepsilon),\,z_2(\varepsilon)\) are
located very close to the boundary point \(z=0\). This can be
shown by the asymptotic calculation. Since \(f_{\varepsilon}(z)
=6\varepsilon+{\textstyle\frac{1-\varepsilon}{4}}z^4+o(|z|^4)\) as
\(z\to{}0\),  the equation \(f_{\varepsilon}(z)=0\) has the roots
\(z_1(\varepsilon),\,z_2(\varepsilon)\) which behave as
\[z_{1,2}(\varepsilon)=\varepsilon^{\frac{1}{4}}24^{\frac{1}{4}}e^{\pm\frac{\pi}{4}}
(1+o(\varepsilon))\ \ \textup{as} \ \varepsilon\to{}+0\,.\] Since
for \({\varepsilon}>0\) the function \(f_{\varepsilon}\) has only
two roots in the open right half plane, there are no other roots
than \(z_1(\varepsilon),\,z_2(\varepsilon)\) there. Since
\(f(z)=\lim_{\varepsilon\to{}+0}f_{\varepsilon}(z)\), the function
\(f(z)\) has no roots in the right upper half plane. (We apply
Hurwitz's theorem.) From \eqref{BDiCa} it follows that
 the function \(f\) does not
vanish on the imaginary axis except the  point \(z=0\). At this
point the function \(f\) has the root of multiplicity four. Thus,
for \(q=4\) the function \(I_q(z)\) is in the Hurwitz class.
\hfill\framebox[0.45em]{ }\\

 Claim 1 of Lemma \ref{HPEGM} is a consequence of Lemma \ref{HurCl}
 and Lemma \ref{QEqFour}.

\begin{remark}
\label{Agr}%
 From \eqref{IRPs1} and \eqref{ExpCal4}, the explicit expression follows:
 \begin{equation}
 \label{ExAg}
 {\sideset{_1}{_1}\Psii}%
\left\{\genfrac{}{}{0pt}{1}%
{\frac{1}{2}}%
{1};%
\genfrac{}{}{0pt}{1}%
{\frac{1}{2}}%
{1+\frac{q}{2}};%
z\right\}=2\frac{(2z^2-6z+6)e^z+(z^2-6)}{z^4}\ \ \textup{for} \ \
q=4\,.
 \end{equation} This expression agrees with the asymptotic
 \eqref{Asymp1}.
\end{remark}
\paragraph{Analysis of the function
\(\mathcal{W}_{\partial{}(B^{\infty}\times{}0)}^{\,p}\):}  \ \ \ \
We may
calculate the function\\
\(\mathcal{W}_{\partial{}(B^{\infty}\times{}0)}^{\,p}\) \ \
asymptotically expressing it in terms of the appropriate
Fox-Wright function, \eqref{As1to}, and then refer to the
asymptotic expansion of this Fox-Wright function. However, we
derive the asymptotic of the function\\ %
\(\mathcal{W}_{\partial{}(B^{\infty}\times{}0)}^{\,p}\) \,from the
asymptotic of the function \(\mathcal{M}_{B^n\times{}0^p}(t)\).
From \eqref{MInt} and \eqref{MInt} it follows that
\begin{equation}%
\label{RelMW}%
 \mathcal{W}_{\partial{}(B^{\infty}\times{}0)}^{\,p}(t)=
\frac{1}{2}\big(\mathcal{M}_{B^n\times{}0^p}(it)+\mathcal{M}_{B^n\times{}0^p}(-it)\big)\,.
\end{equation}%
Comparing \eqref{RelMW} with \eqref{Asymp1}, we see that
\begin{multline}%
\label{AsyW}%
 \mathcal{W}_{\partial{}(B^{\infty}\times{}0)}^{\,p}(t)=\\
=%
 \begin{cases}%
\hspace*{4.0ex}2^{\frac{p}{2}}\hspace*{0.8ex}
t^{-\frac{p}{2}}\hspace*{1.0ex} \cos{(t-\frac{\pi{}p}{4}})+
\frac{2}{\Gamma(\frac{p}{2})}t^{-2}+r_1(t), &|\arg{}t|\leq\varepsilon\,,\\
2^{\frac{q}{2}}
(te^{-i\pi})^{-\frac{p}{2}}\cos{(t+\frac{\pi{}p}{4}})+
\frac{2}{\Gamma(\frac{p}{2})}t^{-2}+r_2(t),\ \ &|\arg{}t-\pi|\leq\varepsilon\,,\\
2^{\frac{p}{2}}(t\,e^{\mp{}\frac{i\pi}{2}})^{-\frac{p}{2}}\,e^{\mp{}it}(1+r_3(t))\,,
&|\arg{}t\pm{}\frac{\pi}{2}|\leq{}\frac{\pi}{2}-\varepsilon\,,
\end{cases}%
\end{multline}%
where the reminders \(r_1(t),\,r_2(t),\,r_2(t)\) admit the
estimates
\begin{subequations}
\label{EstRema}
\begin{equation}%
\label{EstRema1}%
|r_1(t)|\leq{}C_1(\varepsilon)\left(|t|^{-(1+\frac{p}{2})}e^{|\textup{Im}\,t|}+|t|^{-3}\right)\,,
\ \ \ \ \ \ \ |\arg{}t|\leq\varepsilon\,,
\end{equation}%
\begin{equation}%
\label{EstRema2}%
|r_2(t)|\leq{}C_2(\varepsilon)\left(|t|^{-(1+\frac{p}{2})}e^{|\textup{Im}\,t|}+|t|^{-3}\right)\,,
\ \    |\arg{}t-\pi|\leq\varepsilon\,,
\end{equation}%
\begin{equation}%
\label{EstRema3}%
|r_3(t)|\leq{}C_3(\varepsilon)\,|t|^{-1}|,\hspace*{9.0ex}
 |\arg{}t\mp{}\frac{\pi}{2}|\leq{}\frac{\pi}{2}-\varepsilon\,,
\end{equation}%
\end{subequations}%
 and
\(C_1(\varepsilon)<\infty,\,C_2(\varepsilon)<\infty,\,C_3(\varepsilon)<\infty\)
for every \(\varepsilon:\,0<\varepsilon<\frac{\pi}{2}.\) Moreover,
the function
\(\mathcal{W}_{\partial{}(B^{\infty}\times{}0)}^{\,p}(t)\) is even
function of \(t\), and takes real values at real \(t\).

From \eqref{AsyW}, \eqref{EstRema} it follows that for every
\(\varepsilon>0\) the function
\(\mathcal{W}_{\partial{}(B^{\infty}\times{}0)}^{\,p}(t)\) may
have not more that finitely many roots within the angles
\(\{t:\,|\arg{}(t\mp{}\frac{\pi}{2}|\leq{}\frac{\pi}{2}-\varepsilon\}\,.\)
Within the angle \(\{t:\,|\arg{}t|\leq{}\varepsilon\}\), the
function
\(\mathcal{W}_{\partial{}(B^{\infty}\times{}0)}^{\,p}(t)\) has
infinitely many roots, and these roots are asymptotically close to
the roots of the approximating function
\[f_p(t)=2^{\frac{p}{2}} t^{-\frac{p}{2}}
\cos{(t-\textstyle{\frac{\pi{}p}{4}}})+
\frac{2}{\Gamma(\frac{p}{2})}t^{-2}\,,\ \ \
|\arg{}t|\leq{}\varepsilon\,.\]
 (Since the function
\(\mathcal{W}_{\partial{}(B^{\infty}\times{}0)}^{\,p}(t)\) is
even, there is no need to study its behavior within the angle
\(\{t:\,|\arg{}t-\pi|\leq{}\varepsilon\}\).) The behavior of roots
of the approximating equation \(f_p(t)=0\), that is the equation
\begin{equation}
\label{ApprEqp}%
 \cos{(t-\textstyle{\frac{\pi{}p}{4}}})+
\frac{2^{1-\frac{p}{2}}}{\Gamma(\frac{p}{2})}\,t^{\frac{p}{2}-2}=0\,,
\ \ \ |\arg{}t|\leq{}\varepsilon\,,
\end{equation}
depends on \(p\).%
\begin{itemize}
\item[]%
 If \(0<p<4\), then all but finitely many roots of the
equation \eqref{ApprEqp} are real and simple, and these roots are
asymptotically close to the roots of the equation
\(\cos{}(t-\textstyle{\frac{\pi{}p}{4}})=0\).  %
\item[]%
If \(p=4\), then all but finitely many roots of the equation
\eqref{ApprEqp} are real and simple, and these roots are
asymptotically close to the roots of the equation
\(\cos{}(t-\pi)=0\).  %
\item[]%
If \(p>4\), then all but finitely many roots of the equation
\eqref{ApprEqp} are \textsf{non}-real and simple, they are located
symmetrically with respect to the real axis, and are
asymptotically close to the `logarithmic parabola'
\[|y|=({\textstyle\frac{p}{2}}-1)
\ln(|x|+1)+\ln{}c_p\,,\ \
c_p=\frac{2^{1-\frac{p}{2}}}{\Gamma(\frac{p}{2})}\,,\ \ \
0\leq{}x<\infty\,.\]
\end{itemize}
Thus we prove the following
\begin{lemma}
\label{AnWp}%
 For each \(p:\,0\leq{}p<\infty\), the function
 \(\mathcal{W}_{\partial{}(B^{\infty}\times{}0)}^{\,p}(t)\) has
infinitely many roots.
 All but finitely many these roots are
 simple. They lie symmetric with respect to the point \(z=0\).
 \vspace{-3.0ex}
 \begin{enumerate}
 \item
 If \(0\leq{}p\leq{}4\), then all but finitely many these
 roots are real\,;
 \item
 If \ \(4<p\), then all but finitely many these
 roots are non-real.
 In particular, the function
\(\mathcal{W}_{\partial{}(B^{\infty}\times{}0)}^{\,p}(t)\) does
not belong to the Laguerre-Polya class
\(\mathscr{L}\text{-}\mathscr{P}\).
 \end{enumerate}
\end{lemma}
\begin{lemma}
\label{WPBLP}%
 For \(0<p\leq{}2\), as well as for \(p=4\), the
function \(\mathcal{W}_{\partial{}(B^{\infty}\times{}0)}^{\,p}\)
belongs to the Laguerre-Polya class.
\end{lemma}
\textsf{PROOF}. The equality \eqref{RelMW} is a starting point of
our reasoning. If the function \(\mathcal{M}_{B^n\times{}0^p}(t)\)
is in the Hurwitz class \(\mathscr{H}\), then the function
\begin{equation}
\label{RHC} %
\Omega(t)=\mathcal{M}_{B^n\times{}0^p}(it)
\end{equation}
is in the class \(\mathscr{P}\) in the sense of of \cite{Lev1}.
\begin{definition}\textup{[B.Levin, \cite{Lev1},Chapter VII,
Section 4.] } \label{ClP} An entire function \(\Omega(t)\) of
exponential type belongs to the class
\(\mathscr{P}\) if:%
 \begin{enumerate}
 \item
\(\Omega(t)\) has no roots in the closed lower half-plane
\(\{t:\,\textup{Im}\,t\leq{}0\}.\)
 \item
 The defect \(d_{\Omega}\) of the function \(\omega\) is
 non-negative, where
 \[2d_{\Omega}=\varlimsup\limits_{r\to+\infty}\frac{\ln|\Omega(-ir)|}{r}-
\varlimsup\limits_{r\to+\infty}\frac{\ln|\Omega(ir)|}{r}\,. \]
 \end{enumerate}
\end{definition}

 In the book \cite{Lev1} of B.Ya.Levin, the following version of the Hermite-Bieler
Theorem is proved:
\begin{nonumtheorem}\textup{[\,\cite{Lev1}, Chapter VII,
Section 4, Theorem 7\,]} If an entire function \(\Omega(t)\) is in
class \(\mathscr{P}\), then its real and imaginary parts
\({}^{\mathscr{R}}\Omega(t)\) and \({}^{\mathscr{I}}\Omega(t)\):
\[{}^{\mathscr{R}}\Omega(t)=\frac{\Omega(t)+\overline{\Omega(\overline{t})}}{2},
\quad{}{}^{\mathscr{I}}\Omega(t)=\frac{\Omega(t)-\overline{\Omega(\overline{t})}}{2i},\]
possess the properties:\\[-4.0ex]
\begin{enumerate}
\item
The roots of each of the functions \({}^{\mathscr{R}}\Omega(t)\)
and \({}^{\mathscr{I}}\Omega(t)\) are real and simple;
\item
The root sets of the functions \({}^{\mathscr{R}}\Omega(t)\) and
\({}^{\mathscr{I}}\Omega(t)\) interlace.
\end{enumerate}
\end{nonumtheorem}
Let us apply this theorem to the function \(\Omega(t)\) defined by
\eqref{RHC}:
\[\Omega(t)=\mathcal{M}_{B^n\times{}0^p}(it)\,,\] taking into
account that the function \(\mathcal{M}_{B^n\times{}0^p}\) is
real:\\
\(\mathcal{M}_{B^n\times{}0^p}(t)\equiv\overline{\mathcal{M}_{B^n\times{}0^p}(\overline{t})}\),
or what is the same,
\(\mathcal{M}_{B^n\times{}0^p}(-it)\equiv\overline{\mathcal{M}_{B^n\times{}0^p}(i\overline{t})}\).
Hence,
\[{}^{\mathscr{R}}\Omega(t)=\frac{1}{2}
\big(\mathcal{M}_{B^n\times{}0^p}(it)+\mathcal{M}_{B^n\times{}0^p}(-it)\big)\,,\]
or taking \eqref{RelMW} into account,
\[{}^{\mathscr{R}}\Omega(t)= \mathcal{W}_{\partial{}(B^{\infty}\times{}0)}^{\,p}(t)\,.\]
Thus, the following result holds:
\begin{lemma}
\label{CondL}%
 If the function \(\mathcal{M}_{B^n\times{}0^p}(t)\) belongs to
 the Hurwitz class \(\mathscr{H}\), then the function
 \(\mathcal{W}_{\partial{}(B^{\infty}\times{}0)}^{\,p}(t)\)
 belongs to the Laguerre-Polya class \(\mathscr{L}\text{-}\mathscr{P}\).
\end{lemma}
Combining Lemma \ref{CondL} with Lemmas \ref{HurCl} and
\ref{QEqFour}, we obtain Lemma \ref{WPBLP}. \mbox{ \
}\hfill\framebox[0.45em]{ }
\begin{remark}%
\label{Never}%
For the function \(\Omega(t)\) of the form \eqref{RHC}, its real
part \({}^{\mathscr{R}}\Omega(t)\) \({}^{\mathscr{I}}\Omega(t)\)
has infinitely many non-real roots if \(p>4\). Nevertheless,  all
roots of the imaginary part \({}^{\mathscr{I}}\Omega(t)\) are real
for every \(p\geq{}0\).
\end{remark}%
 Indeed, according to
\eqref{MInt} and \eqref{RHC},
\begin{equation}
\label{ImPaO}%
 {}^{\mathscr{I}}\Omega(t)=p\int\limits_{0}^{1}
(1-\xi^2)^{\frac{p}{2}-1}\xi\,{}\sin{\xi{}t}\,d\xi\,.
\end{equation}
According to A.Hurwitz (see \textup{\cite{Wat}, section 15.27}),
for every \(\nu\geq{}-1\), all roots of the entire function
\[\Big(\frac{t}{2}\Big)^{-\nu}J_{\nu}(t)=
\sum\limits_{l=0}^{\infty}\frac{(-1)^l}{l!\,\Gamma(\nu+l+1)}\Big(\frac{t^2}{4}\Big)^l\]
 are real. (\(J_{\nu}(t)\) is the
Bessel function of the index \(\nu\).) For \(\nu>\frac{1}{2}\),
the function \((\frac{t}{2})^{-\nu}J_{\nu}(t)\) admits the
integral representation
\begin{equation}
\label{InReBe} \frac{1}{2}\Big(\frac{t}{2}\Big)^{-\nu}J_{\nu}(t)=
\frac{1}{\Gamma(\nu+\frac{1}{2})\,\Gamma(\frac{1}{2})}
\int\limits_0^1(1-\xi^2)^{\nu-\frac{1}{2}}\cos{t\xi}\,d\xi\,.
\end{equation}
Thus, for \(\nu>-\frac{1}{2}\) all roots of the entire function
\(\int\limits_0^1(1-\xi^2)^{\nu-\frac{1}{2}}\cos{t\xi}\,d\xi\) are
real.  If all roots of a real entire function of exponential type
are real, then all roots of its derivative are real as well. Thus,
for \(\nu>-\frac{1}{2}\) all roots of the function
\(\int\limits_0^1(1-\xi^2)^{\nu-\frac{1}{2}}\xi\sin{t\xi}\,d\xi\)
are real. However, for \(\nu=\frac{p-1}{2}\), the last function
coincides with the function
\(\frac{1}{p}\,\,{}^{\mathscr{I}}\Omega(t)\).
\hfill\framebox[0.45em]{ }

 For \(p=3\), we do not know whether the function
 \(\mathcal{W}_{\partial{}(B^{\infty}\times{}0)}^{\,p}(t)\)
 belongs to the Laguerre-Polya class or not. Our conjecture is
 that YES. Let us formulate our conjectures more precisely.
 Let us formulate our conjectures in terms of the Fox-Wright
 functions.

 \noindent
 \textsf{CONJECTURE 1.} \textit{For \(0\leq{}\lambda\leq{}2\), all roots of
 the Fox-Wright function
 \begin{equation}
 \label{FRHur}
 {\sideset{_1}{_1}\Psii}%
\left\{\genfrac{}{}{0pt}{1}%
{\frac{1}{2}}%
{1};%
\genfrac{}{}{0pt}{1}%
{\frac{1}{2}}%
{1+\lambda};%
t\right\}=\sum\limits_{0\leq{}k<\infty}
\frac{\Gamma(\frac{k}{2}+1)}{\Gamma(\frac{k}{2}+1+\lambda)}\frac{t^k}{k!}
 \end{equation}
 lie in the open left half plane}.

 We proved that the answer is affirmative for
 \(0\leq{}\lambda\leq{}1\), and for \(\lambda=2\).

 \noindent
 \textsf{CONJECTURE 2.}
 \textit{For \(0\leq{}\lambda\leq{}2\), all roots of
 the Fox-Wright function
 \begin{equation}
 \label{FRLP}
 {\sideset{_1}{_2}\Psii}%
\left\{\genfrac{}{}{0pt}{1}%
{1}%
{1};%
\genfrac{}{}{0pt}{1}%
{1}%
{\frac{1}{2}\,}\,%
\genfrac{}{}{0pt}{1}%
{1}%
{1+\lambda};%
t\right\}=\sum\limits_{0\leq{}l<\infty}
\frac{\Gamma(l+1)}{\Gamma(l+\frac{1}{2})\Gamma(l+1+\lambda)}\,\frac{t^l}{l!}
\end{equation}
 are negative and simple.}

 We proved that the answer is affirmative for
 \(0\leq{}\lambda\leq{}1\), and for \(\lambda=2\).

From Hermite-Bieler theorem it follows that if  Conjecture 1 holds
for some \(\lambda\), then for this \(\lambda\) Conjecture 2 holds
as well.

The conjectures 1 and 2 are related to some deep questions related
to  `meromorphic multiplier sequences'. (See \cite{CrCs4}, Problem
1.1.)

\section{CONCLUDING REMARKS. \label{ConRem}}
\textbf{1.} In the present paper we use geometric consideration to
a small extent. The only general geometric tool which we used was
the Alexandrov-Fenchel inequalities. We did not use the
monotonicity properties of the coefficients of the Minkowski
polynomials. If \(V_0,\,V_1,\,V_2\) are convex sets, such that
\[V_1\subseteq{}V_0\subseteq{}V_2\]
 and
 \(M_{V_0}(t),\,M_{V_1}(t),\,M_{V_2}(t)\) are their Minkowskii
 polynomials,
 \[M_{V_j}(t)=\sum\limits_{0\leq{}k\leq{}n}m_k(V_j)t^k,\quad{}j=0,\,1,\,2,\]
then for the coefficients of these polynomials the inequalities
\[m_k(V_1)\leq{}m_k(V_0)\leq{}m_k(V_2)\,,\quad{}0\leq{}k\leq{}n.\]
hold. In this connection, the use of the \textit{Kharitonov
criterion of stability} may be helpful. (Concerning the Haritonov
criterion see Chapters 5 and 7 of the book \cite{BCK} and the
literature quoted there.) The Kharitonov criterion deals with the
`interval stability' of polynomials. In its simplest form, this
criterion allow to determine whether the polynomial
\begin{equation}%
\label{IntPO}
 f(t)=\sum\limits_{0\leq{}k\leq{}n}a_kt^k
\end{equation}
with the real coefficients \(a_k\) is stable from the information
that these coefficients belongs to some intervals:
\begin{equation}%
\label{CoefInt}%
 a_k^{-}\leq{}a_k\leq{}a_k^{+}\,,\quad{}0\leq{}k\leq{}n\,.
\end{equation}
Applying this criterion, one need to construct certain polynomials
 from the given numbers \(a_k^{-},\,a_k^{+},\
0\leq{}k\leq{}n\,.\) (There are finitely many such polynomials.)
If all these polynomials are stable, then the arbitrary polynomial
\(f(t)\), \eqref{IntPO}, whose coefficients satisfy the
inequalities \eqref{CoefInt}, is stable.

\noindent%
\textbf{2.} In the example of a convex set \(V\) whose Minkowski
polynomial is not dissipative, the set \(V\) are very `flattened'
in some direction. (See Theorem \ref{NMR}.)

 \textit{What one can
say about Minkowski polynomials of those convex set \(V\) which
are `isotropic'? }

The notion of \textit{isotropy} may be defined in the following
way.
\begin{definition}
\label{DefIso} The solid convex set \(V\),
\(V\subset{}\mathbb{R}^n\), is said to be isotropic (with respect
to the point \(0\)), if the integral
\[\int\limits_{V}|\langle{}x,e\rangle|^2dv_n(x)\]
takes the same value (i.e. is constant with respect to \(e\)) for
every vector \(e\in\mathbb{R}^n\) such that
\(\langle{}e,e\rangle=1\,.\) Here \(\langle{}\,.\,,\,.\,\rangle\)
is the standard scalar product in \(\mathbb{R}^n\), and
\(dv_n(x)\) is the standard \(n\)-dimensional element on volume.
\end{definition}
\newpage

\end{document}